\documentclass[12pt]{amsart}
\usepackage{amssymb}
\usepackage{amsmath}
\usepackage{amsfonts}
\usepackage{setspace,hyperref}
\usepackage{graphicx}
\usepackage{color}

\numberwithin{equation}{section}
\theoremstyle{plain}

\input{tcilatex}

\begin{document}
\title[Testing Functional Inequalities]{Testing Functional Inequalities}
\thanks{We would like to thank an editor, an associate editor, three
anonymous referees, David Mason, and Oliver Linton for their helpful
comments on earlier versions of this paper. Lee thanks the Economic and
Social Research Council for the ESRC Centre for Microdata Methods and
Practice (RES-589-28-0001) and the European Research Council for the
research grant (ERC-2009-StG-240910-ROMETA). Whang thanks the Korea Research
Foundation for the research grant (KRF-2009-327-B00094).}
\date{14 August 2012.}
\author[Lee]{Sokbae Lee}
\address{Department of Economics, Seoul National University, 1 Gwanak-ro,
Gwanak-gu, Seoul, 151-742, Republic of Korea, and Centre for Microdata
Methods and Practice, Institute for Fiscal Studies, 7 Ridgmount Street,
London, WC1E 7AE, UK.}
\email{sokbae@gmail.com}
\author[Song]{Kyungchul Song}
\address{Department of Economics, University of British Columbia, 997 - 1873
East Mall, Vancouver, BC, V6T 1Z1, Canada}
\email{kysong@mail.ubc.ca}
\author[Whang]{Yoon-Jae Whang}
\address{Department of Economics, Seoul National University, 1 Gwanak-ro,
Gwanak-gu, Seoul, 151-742, Republic of Korea.}
\email{whang@snu.ac.kr}

\begin{abstract}
{\footnotesize This paper develops tests for inequality constraints of
nonparametric regression functions. The test statistics involve a one-sided
version of $L_p$-type functionals of kernel estimators $(1 \leq p < \infty)$%
. Drawing on the approach of Poissonization, this paper establishes that the
tests are asymptotically distribution free, admitting asymptotic normal
approximation. In particular, the tests using the standard normal critical
values have asymptotically correct size and are consistent against general
fixed alternatives. Furthermore, we establish conditions under which the
tests have nontrivial local power against Pitman local alternatives. Some
results from Monte Carlo simulations are presented. \newline
}

{\footnotesize \ }

{\footnotesize \noindent \textsc{Key words.} Conditional moment
inequalities, kernel estimation, one-sided test, local power, $L_p$ norm,
Poissonization. \newline
}

{\footnotesize \noindent \textsc{JEL Subject Classification.} C12, C14. }

\bigskip

{\footnotesize \noindent \textsc{AMS Subject Classification.} 62G10, 62G08,
62G20. }
\end{abstract}

\maketitle

\onehalfspacing

\section{Introduction}

Suppose that we observe $\{(Y_{i}^{\prime },X_{i}^{\prime })^{\prime
}\}_{i=1}^{n}$ that are i.i.d. copies from a random vector, $(Y^{\prime
},X^{\prime })^{\prime }\in \mathbf{R}^{J}\times \mathbf{R}^{d}$. Write $%
Y_{i}=(Y_{1i},\cdot \cdot \cdot ,Y_{Ji})^{\prime }\in \mathbf{R}^{J}$ and
define $m_{j}(x)\equiv \mathbf{E}[Y_{ji}|X_{i}=x],$ $j=1,2,\cdot \cdot \cdot
,J$. The notation $\equiv $ indicates definition.

This paper focuses on the problem of testing functional inequalities:%
\begin{equation}
\begin{split}
H_{0}& :m_{j}(x)\leq 0\text{ for all $(x,j)\in \mathcal{X}\times \mathcal{J}$%
,}\text{\ vs.} \\
H_{1}& :m_{j}(x)>0\text{ for some $(x,j)\in \mathcal{X}\times \mathcal{J}$},
\end{split}
\label{ineqtest}
\end{equation}%
where $\mathcal{X}\subset \mathbf{R}^{d}$ is the domain of interest and $%
\mathcal{J}\equiv \{1,\ldots ,J\}$. Our testing problem is relevant in
various applied settings. For example, in a randomized controlled trial, a
researcher observes either an outcome with treatment $(W_{1})$ or an outcome
without treatment $(W_{0})$ along with observable pre-determined
characteristics of the subjects ($X$). Let $D=1$ if the subject belongs to
the treatment group and 0 otherwise. Suppose that assignment to treatment is
random and independent of $X$ and that the assignment probability $p\equiv
P\{D=1\},\,0<p<1,$ is fixed by the experiment design. Then the average
treatment effect $\mathbf{E}(W_{1}-W_{0}|X=x)$, conditional on $X=x$, can be
written as 
\begin{equation*}
\mathbf{E}(W_{1}-W_{0}|X=x)=\mathbf{E}\left[ \frac{DW}{p}-\frac{(1-D)W}{1-p}%
\bigg|X=x\right] ,
\end{equation*}%
where $W\equiv DW_{1}+(1-D)W_{0}$. In this setup, it may be of interest to
test whether or not $m(x)\equiv \mathbf{E}(W_{1}-W_{0}|X=x)\leq 0$ for all $%
x $.

In economic theory, primitive assumptions of economic models generate
certain testable implications in the form of functional inequalities. For
example, Chiappori, Jullien, Salani\'{e}, and Salani\'{e} (2006) formulated
some testable restrictions in the study of insurance markets. Our tests are
applicable for testing their restrictions (e.g. equation (4) of Chiappori,
Jullien, Salani\'{e}, and Salani\'{e} (2006)). Furthermore, our method can
be used to test for monotone treatment response (see, e.g. Manski (1997)).
For example, testing for a decreasing demand curve for each level of price
in treatments and for each value of covariates falls within the framework of
this paper.

Our test statistic can also be used to construct confidence regions for a
parameter that is partially identified under conditional moment
inequalities. See, among many others, Andrews and Shi (2011a,b), Armstrong
(2011), Chernozhukov, Lee, and Rosen (2009), Chetverikov (2012), and
references therein for inference with conditional moment inequalities.

This paper proposes a one-sided $L_{p}$ approach in testing nonparametric
functional inequalities. While measuring the quality of an estimated
nonparametric function by its $L_{p}$-distance from the true function has
long received attention in the literature (see Devroye and Gy\"{o}rfi
(1985), for an elegant treatment of the $L_{1}$ norm of nonparametric
density estimation), the advance of this approach for general nonparametric
testing seems to have been rather slow relative to other approaches, perhaps
due to its technical complexity. 

Cs\"{o}rg\H{o} and Horv\'{a}th (1988) first established a central limit
theorem for the $L_{p}$-distance of a kernel density estimator from its
population counterpart, and Horv\'{a}th (1991) introduced a Poissonization
technique into the analysis of the $L_{p}$-distance. Beirlant and Mason
(1995) developed a different Poissonization technique and established a
central limit theorem for the $L_{p}$-distance of kernel density estimators
and regressograms from their expected values without assuming smoothness
conditions for the nonparametric functions. Gin\'{e}, Mason and Zaitsev
(2003: GMZ, hereafter) employed this technique to prove the weak convergence
of an $L_{1}$-distance process indexed by kernel functions in kernel density
estimators.

This paper builds on the contributions of Beirlant and Mason (1995) and GMZ
to develop methods for testing \eqref{ineqtest}. In particular, the tests
that we propose are studentized versions of one-sided $L_{p}$-type
functionals. We show that our proposed test statistic is distributed as
standard normal under the least favorable case of the null hypothesis. Thus,
our tests using the standard normal critical values have asymptotically
correct size. We also show that our tests are consistent against general
fixed alternatives and carry out local power analysis with Pitman
alternatives. For the latter, we establish conditions under which the tests
have nontrivial local power against Pitman local alternatives, including
some $n^{-1/2}$-converging Pitman sequences.

Our tests have the following desirable properties. First, our tests do not
require usual smoothness conditions for nonparametric functions for their
asymptotic validity and consistency. This is because we do not need
pointwise or uniform consistency of an unknown function to implement our
tests. For example, a studentized version of our statistic can be estimated
without need for controlling the bias. Second, our tests for \eqref{ineqtest}
are distribution free under the least favorable case of the null hypothesis
where $m_{j}(x)=0,$ for all $x\in \mathcal{X}$ and for all $j\in \mathcal{J}$
and at the same time have nontrivial power against some, though not all, $%
n^{-1/2}$-converging Pitman local alternatives. This is somewhat unexpected,
given that nonparametric goodness-of-fit tests that involve random vectors
of a multi-dimension and have nontrivial power against $n^{-1/2}$-converging
Pitman sequences are not often distribution free. Exceptions are tests that
use an innovation martingale approach (see, e.g., Khmaladze (1993), Stute,
Thies and Zhu (1998), Bai (2003), and Khmaladze and Koul (2004)) or some
tests of independence (or conditional independence) among random variables
(see, e.g., Blum, Kiefer, and Rosenblatt (1961), Delgado and Mora (2000) and
Song (2009)). Third, the local power calculation of our tests for %
\eqref{ineqtest} reveals an interesting contrast with other nonparametric
tests based on kernel smoothers, e.g. H\"{a}rdle and Mammen (1993) and
Horowitz and Spokoiny (2001), where the latter tests are known to have
trivial power against $n^{-1/2}$-converging Pitman local alternatives. Our
inequality tests can have nontrivial local powers against $n^{-1/2}$%
-converging Pitman local alternatives, provided that a certain integral
associated with local alternatives is strictly positive. On the other hand,
it is shown in Section \ref{sec:test-equ} that our equality tests have
trivial power against $n^{-1/2}$-converging Pitman local alternatives.
Therefore, the one-sided nature of inequality testing is the source of our
different local power results. This finding appears new in the literature to
the best of our knowledge.


The remainder of the paper is as follows. Section \ref{lit-rev} discusses
the related literature. Section \ref{sec:theory} provides an informal
description of our test statistic for a simple case, and establishes
conditions under which our tests have asymptotically valid size when the
null hypothesis is true and also are consistent against fixed alternatives.
We also obtain local power results for the leading cases when $p=1$ and $p=2$%
. In Section \ref{sec:test-equ}, we make comparison with functional equality
tests and highlight the main differences between testing inequalities and
equalities in terms of local power. In Section \ref{sec:mc}, we report
results of some Monte Carlo simulations that show that our tests perform
well in finite samples. The proofs of main theorems are contained in Section %
\ref{sec:proof}, along with a roadmap for the proof of the main theorem. 

\section{Related Literature}

\label{lit-rev}

In this section, we provide details on the related literature. The
literature on hypothesis testing involving nonparametric functions has a
long history. Many studies have focused on testing parametric or
semiparametric specifications of regression functions against nonparametric
alternatives. See, e.g., Bickel and Rosenblatt (1973), H\"{a}rdle and Mammen
(1993), Stute (1997), Delgado and Gonz\'{a}lez Manteiga (2000), Horowitz and
Spokoiny (2001), and Khmaladze and Koul (2004) among many others. The
testing problem in this paper is different from the aforementioned papers,
as the focus is on whether certain inequality (or equality) restrictions
hold, rather than on whether certain parametric specifications are plausible.

When $J=1$, our testing problem is also different from testing 
\begin{eqnarray*}
H_{0} &:&m(x)=0\text{ for all $x\in \mathcal{X}$,\thinspace\ against} \\
H_{1} &:&m(x)\geq 0\text{ for all $x\in \mathcal{X}$ with strict inequality
for some $x\in \mathcal{X}.$}
\end{eqnarray*}%
Related to this type of testing problems, see Hall, Huber, and Speckman
(1997) and Koul and Schick (1997, 2003) among others. In their setup, the
possibility that $m(x)<0$ for some $x$ is excluded, so that a consistent
test can be constructed using a linear functional of $m(x)$. On the other
hand, in our setup, negative values of $m(x)$ for some $x$ are allowed under
both $H_{0}$ and $H_{1}$. As a result, a linear functional of $m(x)$ would
not be suitable for our purpose.

There also exist some papers that consider the testing problem in (\ref%
{ineqtest}). For example, Hall and Yatchew (2005) and Andrews and Shi
(2011a,b) considered functions of the form $u\mapsto \max \{u,0\}^{p}$ to
develop tests for (\ref{ineqtest}). However, their tests are not
distribution free, although they achieve local power against some $n^{-1/2}$%
-converging sequences. See also Hall and van Keilegom (2005) for the use of
the one-sided $L_{p}$-type functionals for testing for monotone increasing
hazard rate. None of the aforementioned papers developed test statistics of
one-sided $L_{p}$-type functionals with kernel estimators like ours. See
some remarks of Ghosal, Sen, and van der Vaart (2000, p.1070) on difficulty
in dealing with one-sided $L_{p}$-type functionals with kernel estimators.

In view of Bickel and Rosenblatt (1973) who considered both $L_{2}$ and sup
tests, a one-sided sup test appears to be a natural alternative to the $%
L_{p} $-type tests studied in this paper. For example, Chernozhukov, Lee,
and Rosen (2009) considered a sup norm approach in testing inequality
constraints of nonparametric functions. Also, it may be of interest to
develop sup tests based on a one-sided version of a bootstrap uniform
confidence interval of $\hat{g}_{n}$, similar to Claeskens and van Keilegom
(2003). The sup tests typically do not have nontrivial power against any $%
n^{-1/2}$-converging alternatives, but they may have better power against
some \textquotedblleft sharp peak\textquotedblright\ type alternatives
(Liero, L\"{a}uter and Konakov, 1998).

Testing for inequality is related to testing for monotonicity since a null
hypothesis associated inequality (respectively, monotonicity) can also be
framed as that of monotonicity (respectively, convexity) of integrated
moments. For example, Durot (2003) and Delgado and Escanciano (2011, 2012)
used the least concave majorant operator to characterize their null
hypotheses and developed tests based on the isotonic regression methods.

Finally, we mention that there exist other applications of the
Poissonization method. For example, Anderson, Linton, and Whang (2012)
developed methodology for kernel estimation of a polarization measure; Lee
and Whang (2009) established asymptotic null distributions for the $L_{1}$%
-type test statistics for conditional treatment effects; and Mason (2009)
established both finite sample and asymptotic moment bounds for the $L_{p}$
risk for kernel density estimators. See also Mason and Polonik (2009) and
Biau, Cadre, Mason, and Pelletier (2009) for asymptotic distribution theory
in support estimation.

Among all the aforementioned papers, our work is most closely related to Lee
and Whang (2009), but differs substantially in several important ways.
First, we consider the case of multiple functional inequalities, in contrast
to the single inequality case of Lee and Whang (2009). This extension
requires different arguments (see, e.g. Lemma A7 in Section \ref%
{sec:appendix-lemmas}) and is necessary in order to encompass important
applications such as testing monotonic treatment response and inference with
conditional moment inequalities. Second, we extend the $L_{1}$ statistic to
the general $L_{p}$ statistic. Such an extension is not only theoretically
challenging because many of the results of GMZ apply only to the $L_{1}$
statistic (See, e.g., Lemmas A3 and Lemmas A8 in Section \ref%
{sec:appendix-lemmas}), but also useful to applied econometricians because
the $L_{p}$-type test statistics with different values of $p$ generally have
different power properties. Third, regularity conditions are weaker in this
paper than those in Lee and Whang (2009). In particular, we allow the
underlying functions to be non-smooth, which should be useful in some
contexts. We believe that none of these extensions are trivial. 
Therefore, we view these two papers as complements rather than substitutes.

The testing framework in this paper could be easily extended to testing
stochastic dominance conditional on covariates in the one-sample case or in
the program evaluation setup described in the introduction. For the latter
setup, testing conditional stochastic dominance amounts to testing $%
m(x,y)\equiv \mathbf{E}[1(W_{1} \leq y) -1(W_{0} \leq y)|X=x)\leq 0$ for all 
$(x,y) \in \mathcal{XY}$, where $\mathcal{XY}$ is the domain of the interest
and $W_1$ and $W_0$, as before, are outcomes for treatment and control
groups, respectively. Then a conditional stochastic dominance test can be
developed by combining a density weighted kernel estimator of $m(x,y)$ with
a one-sided $L_{p}$-type functional. However, it is not straightforward to
extend our framework to general two-sample cases. This is because the
propensity score $P(D=1|X=x)$ is unknown in general and has to be estimated
to implement the test. See, for example, Lee and Whang (2009), Delgado and
Escanciano (2011), and Hsu (2011) for testing conditional treatment effects,
including testing conditional stochastic dominance, in general two-sample
cases.

\section{Test Statistics and Asymptotic Properties}

\label{sec:theory}

\subsection{An Informal Description of Our Test Statistics}

Our tests are based on one-sided $L_{p}$-type functionals. For $1\leq
p<\infty ,$ let $\Lambda _{p}:\mathbf{R}\mapsto \mathbf{R}$ be such that $%
\Lambda _{p}(v)\equiv \max \{v,0\}^{p},\ v\in \mathbf{R}$. Consider the
following one-sided $L_{p}$-type functionals:%
\begin{equation*}
\varphi \mapsto \Gamma _{j}(\varphi )\equiv \int_{\mathcal{X}}\Lambda
_{p}(\varphi (x))w_{j}(x)dx,\text{ for }j\in \mathcal{J},
\end{equation*}%
where $w_{j}:\mathbf{R}^{d}\rightarrow \lbrack 0,\infty )\ $is a nonnegative
weight function. Let $f$ denote the density function of $X$ and define $%
g_{j}(x)\equiv m_{j}(x)f(x)$. 
To construct a test statistic, define%
\begin{equation*}
\hat{g}_{jn}(x)\equiv \frac{1}{nh^{d}}\sum_{i=1}^{n}Y_{ji}K\left( \frac{%
x-X_{i}}{h}\right) ,
\end{equation*}%
where $K:\mathbf{R}^{d}\mapsto \mathbf{R}$ is a kernel function and $h$ a
bandwidth parameter satisfying $h\rightarrow 0$ as $n\rightarrow \infty $.
Our test statistic is a suitably studentized of version of $\Gamma _{j}(\hat{%
g}_{jn}(x) )$'s.

Note that we focus on values of $x$ for which $\hat{g}_{jn}(x)>0$ through
the use of $\Lambda _{p}(v)$. Thus, we expect that when $H_{0}$ is true, a
suitably studentized version of $\Gamma _{j}(\hat{g}_{jn})$ is
\textquotedblleft not too large\textquotedblright\ for each $j\in \mathcal{J}
$ but that when $H_{0}$ is false, it will diverge for some $j\in \mathcal{J}$%
. This motivates the use of a weighted sum of $\Gamma_{j}(\hat{g}_{jn})$ as
a test statistic. We require that at least one component of $X$ be
continuously distributed. If some elements of $X$ are discrete, we can
modify the integral in the functional above by using some product measure
between the Lebesgue and counting measures.

We show in Section \ref{asymp-theory-validity} that under weak assumptions,
there exist nonstochastic sequences $a_{jn}\in \mathbf{R},$ $j\in \mathcal{J}
$, and $\sigma _{n}\in (0,\infty )$ such that as $n\rightarrow \infty $, 
\begin{equation}
T_{n}\equiv \frac{1}{\sigma _{n}}\sum_{j=1}^{J}\left\{
n^{p/2}h^{(p-1)d/2}\Gamma _{j}(\hat{g}_{jn})-a_{jn}\right\} \overset{d}{%
\rightarrow }N(0,1),  \label{Tn-def}
\end{equation}%
under the least favorable case of the null hypothesis, where $m_{j}(x)=0,$
for all $x\in \mathcal{X}$ and for all $j\in \mathcal{J}$. This is done
first by deriving asymptotic results for the Poissonized version of the
processes, $\{\hat{g}_{jn}(x):x\in \mathcal{X}\},$ $j\in \mathcal{J}$, and
then by translating them back into those for the original processes through
the de-Poissonization lemma of Beirlant and Mason (1995). See Appendix \ref%
{sec:appendix:roadmap} for details.

To construct a test statistic, we replace $a_{jn}$ and $\sigma _{n}$ by
appropriate estimators to obtain a feasible version of $T_{n},$ say, $\hat{T}%
_{n}$, and show that the limiting distribution remains the same under a
stronger bandwidth condition. Hence, we obtain a distribution free and
consistent test for the nonparametric functional inequality constraints.

To provide a preview of local power analysis with Pitman alternatives in
Section \ref{sec:lap}, suppose that $J=1$ and $p=1$, and the form of the
local alternatives is $g_{1}(x)= \varrho_n \delta _{1}(x)$ for some function 
$\delta _{1}(x)$, where $\varrho_n$ is a sequence of real numbers that
converges to $0$ as $n \rightarrow \infty$. Then (1) if $\int_{\mathcal{X}%
}\delta_{1}(x) w_{1}(x)dx > 0$, our test has nontrivial power against
sequences of local alternatives with $\varrho_n \propto n^{-1/2}$; (2) if $%
\int_{\mathcal{X}}\delta_{1}(x) w_{1}(x)dx = 0$, our test has nontrivial
power only against sequences of local alternatives for which $\varrho_n
\rightarrow 0$ at a rate slower than $n^{-1/2}$; and (3) if $\int_{\mathcal{X%
}}\delta_{1}(x) w_{1}(x)dx < 0$, our test is locally biased whether or not $%
\varrho_n \propto n^{-1/2}$, although our test is a consistent test against
general fixed alternatives.

An alternative statistic is a max statistic such as $\max_{j\in \mathcal{J}%
}\left\{ n^{p/2}h^{(p-1)d/2}\Gamma _{j}(\hat{g}_{jn})-a_{jn}\right\} $,
which we do not pursue in this paper since the ``max'' version of the test
is not typically asymptotically pivotal. 

\subsection{Test Statistics and Asymptotic Validity}

\label{asymp-theory-validity}

Define $\mathcal{S}_{j}\equiv \{x\in \mathcal{X}:w_{j}(x)>0\}$ for each $%
j\in \mathcal{J},$ and, given $\varepsilon >0$, let $\mathcal{S}%
_{j}^{\varepsilon }$ be an $\varepsilon $-enlargement of $\mathcal{S}_{j}$,
i.e., $\mathcal{S}_{j}^{\varepsilon }\equiv \{x+a:x\in \mathcal{S}_{j},\
a\in \lbrack -\varepsilon ,\varepsilon ]^{d}\}$. For $1\leq p<\infty $, let%
\begin{equation}
r_{j,p}(x)\equiv \mathbf{E[}|Y_{ji}|^{p}|X_{i}=x]f(x).  \label{rp}
\end{equation}%
We introduce the following assumptions.\bigskip

\noindent \textsc{Assumption 1:} (i) For each $j \in \mathcal{J}$ and for
some $\varepsilon >0$, $r_{j,2}(x)$ is bounded away from zero and $%
r_{j,2p+2}(x)$ is bounded, both uniformly in $x \in \mathcal{S}%
_{j}^{\varepsilon }$.

\noindent (ii) For each $j\in \mathcal{J},$ $w_{j}(\cdot )$ is nonnegative
on $\mathcal{X}$ and $0<\int_{\mathcal{X}}w_{j}^{s}(x)dx<\infty $, where $%
s\in \{1,2\}$.

\noindent (iii) For $\varepsilon >0$ in (i), $\mathcal{S}_{j}^{\varepsilon
}\subset \mathcal{X}$ for all $j\in \mathcal{J}$.\bigskip

\noindent \textsc{Assumption 2:} $K(u)=\Pi _{s=1}^{d}K_{s}(u_{s}),$ $%
u=(u_{1},\cdot \cdot \cdot ,u_{d})$, with each $K_{s}:\mathbf{R}\rightarrow 
\mathbf{R},$ $s=1,\cdot \cdot \cdot ,d,$\ satisfying that (a) $%
K_{s}(u_{s})=0 $ for all $u_{s}\in \mathbf{R}\backslash \lbrack -1/2,1/2]$,
(b) $K_{s}$ is of bounded variation, and (c) $||K_{s}||_{\infty }\equiv
\sup_{u_{s}\in \mathbf{R}}|K_{s}(u_{s})|<\infty $ and $\int
K_{s}(u_{s})du_{s}=1$.\bigskip

Assumption 1(i) imposes that $\inf \{ r_{j,2}(x):x\in \mathcal{S}%
_{j}^{\varepsilon }\} > 0$ and $\sup \{ r_{j,2p+2}(x):x\in \mathcal{S}%
_{j}^{\varepsilon }\} < \infty$ for each $j\in \mathcal{J}$. Assumption
1(ii) is a weak condition on the weight function. Nonnegativity is important
since we develop a sum statistic over $j$. Assumption 1(iii) is introduced
to avoid the boundary problem of kernel estimators by requiring that $w_{j}$
have support inside an $\varepsilon $-shrunk subset of $\mathcal{X}$. Note
that Assumptions 1(i) and (iii) require that $\mathcal{S}_{j}$ be a bounded
set for each $j\in \mathcal{J}$. The conditions for the kernel function in
Assumption 2 are quite flexible, except that the kernel functions have
bounded support.

Define for $j,k\in \mathcal{J}$ and $x\in \mathbf{R}^{d},$%
\begin{eqnarray*}
\rho _{jk,n}(x) &\equiv &\frac{1}{h^{d}}\mathbf{E}\left[ Y_{ji}Y_{ki}K^{2}%
\left( \frac{x-X_{i}}{h}\right) \right] , \\
\rho _{jn}^{2}(x) &\equiv &\frac{1}{h^{d}}\mathbf{E}\left[
Y_{ji}^{2}K^{2}\left( \frac{x-X_{i}}{h}\right) \right] , \\
\rho _{jk}(x) &\equiv &\mathbf{E}\left[ Y_{ji}Y_{ki}|X_{i}=x\right] f(x)\int
K^{2}(u)du,\text{ and} \\
\rho _{j}^{2}(x) &\equiv &\mathbf{E}[Y_{ji}^{2}|X_{i}=x]f(x)\int K^{2}(u)du.
\end{eqnarray*}%
Let $\mathbb{Z}_{1}$ and $\mathbb{Z}_{2}$\ denote mutually independent
standard normal random variables. We introduce the following quantities:%
\begin{eqnarray}
a_{jn} &\equiv &h^{-d/2}\int_{\mathcal{X}}\rho _{jn}^{p}(x)w_{j}(x)dx\cdot 
\mathbf{E}\Lambda _{p}(\mathbb{Z}_{1})\text{ and}  \label{an} \\
\sigma _{jk,n} &\equiv &\int_{\mathcal{X}}q_{jk,p}(x)\rho _{jn}^{p}(x)\rho
_{kn}^{p}(x)w_{j}(x)w_{k}(x)dx,  \notag
\end{eqnarray}%
where $q_{jk,p}(x)\equiv \int_{\lbrack -1,1]^{d}}Cov(\Lambda _{p}(\sqrt{%
1-t_{jk}^{2}(x,u)}\mathbb{Z}_{1}+t_{jk}(x,u)\mathbb{Z}_{2}),\ \Lambda _{p}(%
\mathbb{Z}_{2}))du$ and%
\begin{equation*}
t_{jk}(x,u)\equiv \frac{\rho _{jk}(x)}{\rho _{j}(x)\rho _{k}(x)}\cdot \frac{%
\int K\left( x\right) K\left( x+u\right) dx}{\int K^{2}\left( x\right) dx}.
\end{equation*}%
Let $\Sigma _{n}$ be a $J\times J$ matrix whose $(j,k)$-th entry is given by 
$\sigma _{jk,n}$. Later we use $\Sigma_n $ to normalize the test statistic.
The scale normalization matrix $\Sigma_n $ does not depend on $x$, and this
is not because we are assuming conditional homoskedasticity in the null
hypothesis, but because $\Sigma_n $ is constituted by covariances of random
quantities that already have $x$ integrated out. We also define $\Sigma $ to
be a $J\times J$ matrix whose $(j,k)$-th entry is given by $\sigma _{jk}$,
where%
\begin{equation*}
\sigma _{jk}\equiv \int_{\mathcal{X}}q_{jk,p}(x)\rho _{j}^{p}(x)\rho
_{k}^{p}(x)w_{j}(x)w_{k}(x)dx.
\end{equation*}%
As for $\Sigma $, we introduce the following assumption.\bigskip

\noindent \textsc{Assumption 3:} $\Sigma $ is positive definite.\bigskip

For example, Assumption 3 excludes the case where $Y_{ji}$ and $Y_{ki}$ $(j
\neq k)$ are perfectly correlated conditional on $X_{i}=x$ for almost all $x$
with $w_j \equiv w_k$.

The following theorem is the first main result of this paper.\bigskip

\noindent \textsc{Theorem 1:} \textit{Suppose that Assumptions 1-3 hold} 
\textit{and that} $h\rightarrow 0$ \textit{and }$n^{-1/2}h^{-d}\rightarrow 0$
\textit{as} $n\rightarrow \infty $. \textit{Furthermore, assume that }$%
m_{j}(x)=0$ \textit{for almost all} $x\in \mathcal{X}$ \textit{and for all }$%
j\in \mathcal{J}.$ \textit{Then } 
\begin{equation*}
T_{n}\equiv \frac{1}{\sigma _{n}}\sum_{j=1}^{J}\left\{
n^{p/2}h^{(p-1)d/2}\Gamma _{j}(\hat{g}_{jn})-a_{jn}\right\} \overset{d}{%
\rightarrow }N(0,1),
\end{equation*}%
\textit{where} $\sigma _{n}^{2}\equiv \mathbf{1}^{\prime }\Sigma _{n}\mathbf{%
1}$, and $\mathbf{1}$ is a vector of ones. \bigskip

Note that when $J=1$, $\sigma _{n}^{2}$ takes the simple form of $q_{p}\int_{%
\mathcal{X}}\rho _{1n}^{2p}(x)w_{1}^{2}(x)dx,$ where%
\begin{eqnarray*}
q_{p} &\equiv &\int_{[-1,1]^{d}}Cov(\Lambda _{p}(\sqrt{1-t^{2}(u)}\mathbb{Z}%
_{1}+t(u)\mathbb{Z}_{2}),\ \Lambda _{p}(\mathbb{Z}_{2}))du,\text{ and} \\
t(u) &\equiv &\int K\left( x\right) K\left( x+u\right) dx/\int K^{2}\left(
x\right) dx.
\end{eqnarray*}

To develop a feasible testing procedure, we construct estimators of $a_{jn}$%
's and $\sigma _{n}^{2}$ as follows. First, define%
\begin{eqnarray}
\hat{\rho}_{jk,n}(x) &\equiv &\frac{1}{nh^{d}}%
\sum_{i=1}^{n}Y_{ji}Y_{ki}K^{2}\left( \frac{x-X_{i}}{h}\right) ,\text{ and}
\label{rho2} \\
\hat{\rho}_{jn}^{2}(x) &\equiv &\frac{1}{nh^{d}}%
\sum_{i=1}^{n}Y_{ji}^{2}K^{2}\left( \frac{x-X_{i}}{h}\right) .  \notag
\end{eqnarray}%
We estimate $a_{jn}$ and $\sigma _{jk,n}$ by:%
\begin{eqnarray*}
\hat{a}_{jn} &\equiv &h^{-d/2}\int_{\mathcal{X}}\hat{\rho}%
_{jn}^{p}(x)w_{j}(x)dx\cdot \mathbf{E}\Lambda _{p}(\mathbb{Z}_{1})\text{ and}
\\
\hat{\sigma}_{jk,n} &\equiv &\int_{\mathcal{X}}\hat{q}_{jk,p}(x)\hat{\rho}%
_{jn}^{p}(x)\hat{\rho}_{kn}^{p}(x)w_{j}(x)w_{k}(x)dx,
\end{eqnarray*}%
where $\hat{q}_{jk,p}(x)\equiv \int_{\lbrack -1,1]^{d}}Cov(\Lambda _{p}(%
\sqrt{1-\hat{t}_{jk}^{2}(x,u)}\mathbb{Z}_{1}+\hat{t}_{jk}(x,u)\mathbb{Z}%
_{2}),\ \Lambda _{p}(\mathbb{Z}_{2}))du$ and%
\begin{equation*}
\hat{t}_{jk}(x,u)\equiv \frac{\hat{\rho}_{jk,n}(x)}{\hat{\rho}_{jn}(x)\hat{%
\rho}_{kn}(x)}\cdot \frac{\int K\left( x\right) K\left( x+u\right) dx}{\int
K^{2}\left( x\right) dx}.
\end{equation*}%
Note that \textbf{$\mathbf{E}$}$\Lambda _{1}(\mathbb{Z}_{1})=1/\sqrt{2\pi }%
\approx 0.39894$ and \textbf{$\mathbf{E}$}$\Lambda _{2}(\mathbb{Z}_{1})=1/2$%
. When $p$ is an integer, the covariance expression in $q_{jk,p}(x)$ can be
computed using the moment generating function of a truncated multivariate
normal distribution (Tallis, 1961). More practically, simulated draws from $%
\mathbb{Z}_{1}$ and $\mathbb{Z}_{2}$ can be used to compute the quantities 
\textbf{$\mathbf{E}$}$\Lambda _{p}(\mathbb{Z}_{1})$ and $q_{jk,p}(x)$ for
general values of $p$. The integrals appearing above can be evaluated using
methods of numerical integration. We define $\hat{\Sigma}_{n}$ to be a $%
J\times J$ matrix whose $(j,k)$-th entry is given by $\hat{\sigma}_{jk,n}$.

Let $\hat{\sigma}_{n}^{2}\equiv \mathbf{1}^{\prime }\hat{\Sigma}_{n}\mathbf{1%
}$. Our test statistic is taken to be%
\begin{equation}
\hat{T}_{n}\equiv \frac{1}{\hat{\sigma}_{n}}\sum_{j=1}^{J}\left\{
n^{p/2}h^{(p-1)d/2}\Gamma _{j}(\hat{g}_{jn})-\hat{a}_{jn}\right\} .
\label{shat}
\end{equation}%
Let $z_{1-\alpha }\equiv \Phi ^{-1}(1-\alpha ),$ where $\Phi $ denotes the
cumulative distribution function of $N(0,1)$. This paper proposes using the
following test:%
\begin{equation}
\text{Reject }H_{0}\ \text{if and only if }\hat{T}_{n}>z_{1-\alpha }\text{.}
\label{test}
\end{equation}%
The following theorem shows that the test has an asymptotically valid size.
\bigskip

\noindent \textsc{Theorem 2:} \textit{Suppose that Assumptions 1-3 hold and
that} $h\rightarrow 0$ \textit{and }$n^{-1/2}h^{-3d/2}\rightarrow 0$, 
\textit{as} $n\rightarrow \infty $. \textit{Furthermore, assume that the
kernel function} $K$ \textit{in Assumption 2 is nonnegative. Then under the
null hypothesis, we have}%
\begin{equation*}
\lim_{n\rightarrow \infty }P\{\hat{T}_{n}>z_{1-\alpha }\}\leq \alpha ,
\end{equation*}%
\textit{with equality holding if }$m_{j}(x)=0$ \textit{for almost all} $x\in 
\mathcal{X}$ \textit{and for all }$j\in \mathcal{J}.$\bigskip

Note that the probability of making an error of rejecting the true null
hypothesis is largest when $m_{j}(x)=0$ for almost all $x\in \mathcal{X}$
and for all $j\in \mathcal{J}$, namely, when we are in \textit{the least
favorable case of the null hypothesis}.

The nonparametric test does not require assumptions for $m_{j}$'s and $f$
beyond those in Assumption 1(i), even after replacing $a_{jn}$'s and $\sigma
_{n}^{2}$ by their estimators. In particular, the theory does not require
continuity or differentiability of $f $ or $m_{j}$'s. This is because we do
not need to control the bias to implement the test. This result uses the
assumption that the kernel function $K$ is nonnegative to control the size
of the test. (See the proof of Theorem 2 for details.)

The bandwidth condition for Theorem 2 is stronger than that in Theorem 1.
This is mainly due to the treatment of the estimation errors in $\hat{a}%
_{jn} $ and $\hat{\sigma}_{n}^{2}.$ For the bandwidth parameter, it suffices
to take $h=c_{1}n^{-s}$ with $0<s<1/(3d)$ for a constant $c_{1}>0.$ In
general, optimal bandwidth choice for nonparametric testing is different
from that for nonparametric estimation as we need to balance the size and
power of the test instead of the bias and variance of an estimator. For
example, Gao and Gijbels (2008) considered testing a parametric null
hypothesis against a nonparametric alternative and derived a
bandwidth-selection rule by utilizing an Edgeworth expansion of the
asymptotic distribution of the test statistic concerned. The methods of Gao
and Gijbels (2008) are not directly applicable to our tests, and it is a
challenging problem to develop a theory of optimal bandwidths for our tests.
We provide some simulation evidence regarding sensitivity to the choice of $%
h $ in Section \ref{sec:mc}.

According to Theorems 1-2, each choice of the weight functions $w_{j}$ leads
to an asymptotically valid test. The actual choice of $w_{j}$ may reflect
the relative importance of individual inequality restrictions. When it is of
little practical significance to treat individual inequality restrictions
differently, one may choose simply $w_{j}(x)=1\{x\in \mathcal{S}\}$ with
some common support $\mathcal{S}$. Perhaps more naturally, to avoid undue
influences of different scales across $Y_{ji}$'s, one may use $w_{j}(x)=%
\tilde{\sigma}_{jj,n}^{-1/2}\bar{w}(x)$, for some common nonnegative weight
function $\bar{w}(x)$, where%
\begin{equation*}
\tilde{\sigma}_{jj,n}\equiv q_{p}\int_{\mathcal{X}}\hat{\rho}_{jn}^{2p}(x)%
\bar{w}^{2}(x)dx,j\in \mathcal{J},
\end{equation*}%
where $\hat{\rho}_{jn}^{2}(x)$ is given as in (\ref{rho2}). Then $\tilde{%
\sigma}_{jj,n}$ is consistent for $\sigma _{jj,n}$ (see the proof of Theorem
2), and just as the estimation error of $\hat{\sigma}_{n}$ in (\ref{test})
leaves the limiting distribution of $T_{n}$ under the null hypothesis
intact, so does the estimation error of $\tilde{\sigma}_{jj,n}$.

The following result shows the consistency of the test in (\ref{test})
against fixed alternatives.\bigskip

\noindent \textsc{Theorem 3:} \textit{Suppose that Assumptions 1-3 hold and
that} $h\rightarrow 0$ \textit{and }$n^{-1/2}h^{-3d/2}\rightarrow 0$, 
\textit{as} $n\rightarrow \infty $. \textit{If $H_1$ is true and $\Gamma
_{j}(g_{j})>0$ for some $j\in \mathcal{J}$, then we have}%
\begin{equation*}
\lim_{n\rightarrow \infty }P\{\hat{T}_{n}>z_{1-\alpha }\}=1.
\end{equation*}
\bigskip

\subsection{Local Asymptotic Power}

\label{sec:lap}

We determine the power of the test in (\ref{test}) against some sequences of
local alternatives. Consider the following sequences of local alternatives
converging to the null hypothesis at the rate$\;n^{-1/2}$, respectively:%
\begin{equation}
H_{\delta }:g_{j}(x)=n^{-1/2}\delta _{j}(x),\text{ for each }j\in \mathcal{J}%
,  \label{ha_d}
\end{equation}%
where $\delta _{j}(\cdot )$'s are bounded real functions on $\mathbf{R}^{d}.$

The following theorem establishes a representation of the local asymptotic
power functions, when $p\in \{1,2\}$. For simplicity of notation, let us
introduce the following definition: for $s\in \{1,2\}$, $z\in \{-1,0,1\}$, a
given weight function vector $w\equiv (w_{1},\cdot \cdot \cdot ,w_{J})$, and
the direction $\delta =(\delta _{1},\cdot \cdot \cdot ,\delta _{J})^{\prime
} $, let $\eta _{s,z}(w,\delta )\equiv \sum_{j=1}^{J}\int_{\mathcal{X}%
}\delta _{j}^{s}(x)\rho _{j}^{z}(x)w_{j}(x)dx,$ and let $\sigma ^{2}\equiv 
\mathbf{1}^{\prime }\Sigma \mathbf{1}$. \bigskip

\noindent \textsc{Theorem 4:} \textit{Suppose that Assumptions 1-3 hold and
that} $h\rightarrow 0$ \textit{and }$n^{-1/2}h^{-3d/2}\rightarrow 0$, 
\textit{as} $n\rightarrow \infty $.

\noindent (i) \textit{If }$p=1$\textit{, then, under }$H_{\delta },\ $%
\textit{we have}%
\begin{equation*}
\lim_{n\rightarrow \infty }P\{\hat{T}_{n}>z_{1-\alpha }\}=1-\Phi
(z_{1-\alpha }-\eta _{1,0}(w,\delta )/2\sigma ).
\end{equation*}%
\textit{\ }

\noindent (ii) \textit{If }$p=2,$ then\textit{, under }$H_{\delta },\ $%
\textit{we have}%
\begin{equation*}
\lim_{n\rightarrow \infty }P\{\hat{T}_{n}>z_{1-\alpha }\}=1-\Phi
(z_{1-\alpha }-\eta _{1,1}(w,\delta )/(\sigma \sqrt{\pi /2})).
\end{equation*}
\bigskip

Theorem 4 gives explicit local asymptotic power functions under $H_{\delta }$%
, when $p=1$ and $p=2$. The local power of the test is greater than the size 
$\alpha $, whenever the ``non-centrality parameter'' ($\eta _{1,0}(w,\delta
)/2\sigma $ in the case of $p=1$ and $\eta _{1,1}(w,\delta )/(\sigma \sqrt{%
\pi /2})$ in the case of $p=2$) is strictly positive. For example, when $J=1$
and $p=1$ (or $p=2)$, the test is asymptotically locally strictly unbiased
as long as $\mu _{\delta }\equiv \int_{\mathcal{X}}\delta_1 (x)w_1(x)dx>0$
(or $\int_{\mathcal{X}}\delta_1 (x)\rho_1 (x)w_1(x)dx>0$). Notice that $\mu
_{\delta }$ can be strictly positive even if $\delta _{1}(x)$ takes negative
values for some $x\in \mathcal{X}$. Therefore, our test has nontrivial local
power against some, though not all, $n^{-1/2}$-local alternatives.

On the other hand, if the noncentrality parameter is zero, the test still
has nontrivial power against local alternatives converging to the null at
the $n^{-1/2}h^{-d/4}$ rate, which is slower than $n^{-1/2}$. To show this,
consider the following local alternatives: 
\begin{equation*}
H_{\delta }^{\ast }:g_{j}(x)=n^{-1/2}h^{-d/4}\delta _{j}(x),\text{ for each }%
j\in \mathcal{J},
\end{equation*}%
where $\delta _{j}(\cdot )$'s are bounded real functions as before$.$
Theorem 4* gives the local asymptotic power functions against $H_{\delta
}^{\ast }.$ \bigskip

\textsc{Theorem 4*:} \textit{Suppose that Assumptions 1-3 hold and that} $%
h\rightarrow 0$ \textit{and }$n^{-1/2}h^{-3d/2}\rightarrow 0$, \textit{as} $%
n\rightarrow \infty $.

\noindent (i) \textit{If }$p=1$\textit{\ and }$\eta _{1,0}(w,\delta )=0$%
\textit{, then, under }$H_{\delta }^{\ast },$\textit{\ we have\ }%
\begin{equation*}
\lim_{n\rightarrow \infty }P\{\hat{T}_{n}>z_{1-\alpha }\}=1-\Phi
(z_{1-\alpha }-\eta _{2,-1}(w,\delta )/\sqrt{8\pi }\sigma ).
\end{equation*}%
\textit{\ }

\noindent (ii) \textit{If }$p=2$\textit{\ and }$\eta _{1,1}(w,\delta )=0,$%
\textit{\ then, under }$H_{\delta }^{\ast },\ $\textit{we have}%
\begin{equation*}
\lim_{n\rightarrow \infty }P\{\hat{T}_{n}>z_{1-\alpha }\}=1-\Phi
(z_{1-\alpha }-\eta _{2,0}(w,\delta )/2\sigma ).
\end{equation*}
\bigskip

If $\eta _{1,0}(w,\delta )=0$ in the case of $p=1$ or $\eta _{1,1}(w,\delta
)=0$ in the case of $p=2$, then the local power of the test is greater than
the size $\alpha $ because the new noncentrality parameter in Theorem 4* is
strictly positive. For example, when $J=1,$ we have $\eta _{2,-1}(w,\delta
)=\int_{\mathcal{X}}\delta _{1}^{2}(x)\rho _{1}^{-1}(x)w_{1}(x)dx>0$ (and $%
\eta _{2,0}(w,\delta )=\int_{\mathcal{X}}\delta _{1}^{2}(x)w_{1}(x)dx>0)$
for all $\delta _{1}.$ Therefore, when $\eta _{1,0}(w,\delta )=0$ or $\eta
_{1,1}(w,\delta )=0,$ Theorem 4$^{\ast }$ implies that our test is strictly
locally unbiased against the $n^{-1/2}h^{-d/4}$ local alternatives $%
H_{\delta }^{\ast }$, though it has only trivial local power $(=\alpha )$
against the $n^{-1/2}$ local alternatives $H_{\delta }$.

To explain the results of Theorems 4 and 4$^{\ast }$ more intuitively,
consider the test statistic $T_{n}$ with $J=1,\ p=2$ and $d=1$. For
simplicity, take $w(\cdot )=1.$\ Let $\sigma \equiv q_{2} \int_{\mathcal{X}%
}\rho_1 ^{4}(x)dx\ $and $a_{n}\equiv h^{-1/2}\int_{\mathcal{X}}\rho_1
^{2}(x)dx\cdot \mathbf{E}\Lambda _{2}(\mathbb{Z}_{1}).$ Let the alternative
hypothesis be given by 
\begin{equation*}
H_{\delta }^{\ast }:g(x)=n^{-1/2}h^{-b}\delta_1 (x),
\end{equation*}%
where $b=0$ or $1/4$. Consider the statistic $\hat{T}_{n}$ with $\hat{\sigma}%
_{n}\ $and $\hat{a}_{n}\ $replaced by their population analogues $\sigma
_{n} $ and $a_{n},$ respectively, i.e., 
\begin{eqnarray}
T_{n} &\equiv &\frac{1}{\sigma _{n}}\left\{ nh^{1/2}\int_{\mathcal{X}%
}\Lambda _{2}\left( \hat{g}_{n}(x)\right) dx-a_{n}\right\}  \notag \\
&=&\frac{nh^{1/2}}{\sigma _{n}}\left\{ \int_{\mathcal{X}}\Lambda _{2}\left( 
\hat{g}_{n}(x)\right) dx-\int_{\mathcal{X}}E\Lambda _{2}\left( \hat{g}%
_{n}(x)\right) dx\right\}  \label{rr1} \\
&&+\frac{nh^{1/2}}{\sigma _{n}}\left\{ \int_{\mathcal{X}}E\Lambda _{2}\left( 
\hat{g}_{n}(x)\right) dx-\frac{a_{n}}{nh^{1/2} } \right\} .  \notag
\label{rr2}
\end{eqnarray}%
It is easy to see that $T_{n}$ has the same asymptotic distribution as $\hat{%
T}_{n}\ $under the local alternative hypothesis. The first term on the right
hand side of (\ref{rr1}) converges in distribution to the standard normal
distribution by the arguments similar to those used to prove Theorem 1.
Consider the second term in (\ref{rr1}). We can approximate it by%
\begin{align}
&\frac{1}{\sigma _{n}}\left\{ nh^{1/2}\int_{\mathcal{X}}E\Lambda _{2}\left( 
\hat{g}_{n}(x)\right) dx-a_{n} \right\}  \notag \\
&=\frac{1}{\sigma _{n}}\left\{ \int_{\mathcal{X}}E\Lambda _{2}\left(
h^{-1/4}\rho_1 (x)\frac{\sqrt{nh}\left[ \hat{g}_{n}(x)-E\hat{g}_{n}(x)\right]
}{\rho_1 (x)}+n^{1/2}h^{1/4}E\hat{g}_{n}(x)\right) dx-a_{n} \right\}  \notag
\\
&\simeq \sigma^{-1} \int_{\mathcal{X}}E\Lambda _{2}\left( h^{-1/4}\rho_1 (x)%
\mathbb{Z}_{1}+n^{1/2}h^{1/4}E\hat{g}_{n}(x)\right) dx- \sigma^{-1} a_{n}
\label{r3_2} \\
&\simeq \sigma^{-1} \int_{\mathcal{X}}\left\{ E\Lambda _{2}\left(
h^{-1/4}\rho_1 (x)\mathbb{Z}_{1}+h^{1/4-b}\delta_1 (x)\right) -E\Lambda
_{2}\left( h^{-1/4}\rho_1 (x)\mathbb{Z}_{1}\right) \right\} dx  \label{r3_3}
\\
&\simeq h^{-b}\left( \frac{2\phi (0)}{\sigma} \int_{\mathcal{X}}\delta_1
(x)\rho_1 (x)dx\right) +h^{1/2-2b}\left( \frac{1}{2\sigma}\int_{\mathcal{X}%
}\delta_1 ^{2}(x)dx\right) ,  \label{r3_4}
\end{align}%
where (\ref{r3_2}) follows from the Poissonization argument, (\ref{r3_3})
holds by $n^{1/2}h^{1/4}E\hat{g}_{n}(x)=h^{1/4-b}\int \delta_1
(x-uh)K(u)du\simeq h^{1/4-b}\delta_1 (x),$ and (\ref{r3_4}) uses a Taylor
expansion $\mathbf{E}\Lambda _{2}(\gamma \mathbb{Z}_{1}+\mu )-\mathbf{E}%
\Lambda _{2}(\gamma \mathbb{Z}_{1})\simeq 2\phi (0)\mu \gamma +\Phi (0)\mu
^{2}\ $with $\gamma =h^{-1/4}\rho_1 (x)$ and $\mu =h^{1/4-b}\delta_1 (x)$,
where $\phi(\cdot)$ and $\Phi(\cdot)$, respectively, denote the pdf and cdf
of the standard normal distribution. This approximation tells us that if $%
\int_{\mathcal{X}}\delta_1 (x)\rho_1 (x)dx >0$, we can take $b=0$ so that it
can achieve nontrivial power against $n^{-1/2}$ alternatives, while if $%
\int_{\mathcal{X}}\delta_1 (x)\rho_1 (x)dx$ $=0,$ then we should take $b=1/4$
so that it has nontrivial local power against $n^{-1/2}h^{-1/4}$ local
alternatives. Notice that, in the latter case, $\int_{\mathcal{X}}\delta_1
^{2}(x)dx$ is always positive.

It would also be interesting to compare local power properties of our test
with that of Andrews and Shi (2011a). Unlike our test, the test of Andrews
and Shi (2011a, Theorem 4(b)) does not require $\int_{\mathcal{X}}\delta_1
(x)\rho_1 (x)dx >0$, but excludes some $n^{-1/2}$-local alternatives. An
analytical and unambiguous comparison between the two approaches is not
straightforward, because the test of Andrews and Shi (2011a) is not
asymptotically distribution free, meaning that the local power function may
depend on the underlying data generating process in a complicated way.
However, we do compare the two approaches in our simulation studies.

When $J=1,$ thanks to Theorem 4, we can compute an optimal weight function
that maximizes the local power against a given direction $\delta $. See
Stute (1997) for related results of optimal directional tests, and Tripathi
and Kitamura (1997) for results of optimal directional and average tests
based on smoothed empirical likelihoods.

Define $\sigma _{p}^{2}(w_{1})\equiv q_{p}\int_{\mathcal{X}}\rho
_{1n}^{2p}(x)w_{1}^{2}(x)dx$ for $J=1$. The optimal weight function (denoted
by $w_{p}^{\ast }$) is taken to be a maximizer of the drift term $\eta
_{1,0}(w_{1},\delta _{1})/\sigma _{1}(w_{1})$ (in the case of $p=1$) or $%
\eta _{1,1}(w_{1},\delta _{1})/\sigma _{2}(w_{1})$ (in the case of $p=2$)
with respect to $w_{1}$ under the constraint that $w_{1}\geq 0$ and $\int_{%
\mathcal{X}}w_{1}(x)\rho ^{2p}(x)dx=1$. The latter condition is for a scale
normalization. Let $\delta _{1}^{+}=\max \{\delta _{1},0\}$. Since $\rho
_{1} $ and $w_{1}$ are nonnegative, the Cauchy-Schwarz inequality suggests
that the optimal weight function is given by%
\begin{equation}
w_{p}^{\ast }(x)=\left\{ 
\begin{tabular}{ll}
$\frac{\delta _{1}^{+}(x)\rho _{1}^{-2}(x)}{\sqrt{\int_{\mathcal{X}}(\delta
_{1}^{+})^{2}(x)\rho _{1}^{-2}(x)dx}},$ & $\text{ if }p=1,\text{ and}$ \\ 
$\frac{\delta _{1}^{+}(x)\rho _{1}^{-3}(x)}{\sqrt{\int_{\mathcal{X}}(\delta
_{1}^{+})^{2}(x)\rho _{1}^{-2}(x)dx}},$ & $\text{ if }p=2.$%
\end{tabular}%
\right.  \label{wf}
\end{equation}%
To satisfy Assumption 1(iii), we assume that the support of $\delta _{1}$ is
contained in an $\varepsilon $-shrunk subset of $\mathcal{X}$. With this
choice of an optimal weight function, the local power function becomes:%
\begin{equation*}
\begin{tabular}{ll}
$1-\Phi \left( z_{1-\alpha }-\frac{\sqrt{\int_{\mathcal{X}}(\delta
_{1}^{+})^{2}(x)\rho _{1}^{-2}(x)dx}}{2\sqrt{q_{1}}}\right) ,$ & $\text{ if }%
p=1,\text{ and}$ \\ 
$1-\Phi \left( z_{1-\alpha }-\frac{\sqrt{\int_{\mathcal{X}}(\delta
_{1}^{+})^{2}(x)\rho _{1}^{-2}(x)dx}}{\sqrt{q_{2}\pi /2}}\right) ,$ & $\text{
if }p=2.$%
\end{tabular}%
\end{equation*}

\section{Comparison with Testing Functional Equalities}

\label{sec:test-equ}

It is straightforward to follow the proofs of Theorems 1-3 to develop a test
for equality restrictions:%
\begin{eqnarray}
H_{0} &:&m_{j}(x)=0\text{ for all $(x,j)\in \mathcal{X}\times \mathcal{J}$,}%
\text{ vs.}  \label{eqtest} \\
H_{1} &:&m_{j}(x)\neq 0\text{ for some $(x,j)\in \mathcal{X}\times \mathcal{J%
}$}.  \notag
\end{eqnarray}%
For this test, we redefine $\Lambda _{p}(v)=|v|^{p}$ and, using this,
redefine $\hat{T}_{n}$ in (\ref{shat}) and $\sigma ^{2}$. Then under the
null hypothesis,%
\begin{equation*}
\hat{T}_{n}\overset{d}{\rightarrow }N(0,1).
\end{equation*}%
Therefore, we can take a critical value in the same way as before. The
asymptotic validity of this test under the null hypothesis in \eqref{eqtest}
follows under precisely the same conditions as in Theorem 2. However, the
convergence rates of the inequality tests and the equality tests under local
alternatives are different, as we shall see now.

Consider the local alternatives converging to the null hypothesis at the rate%
$\;n^{-1/2}h^{-d/4}$:%
\begin{equation}
H_{\delta }^{\ast }:g_{j}(x)=n^{-1/2}h^{-d/4}\delta _{j}(x),\text{ for each }%
j\in \mathcal{J},  \label{hastar}
\end{equation}%
where $\delta _{j}(\cdot )$'s are again bounded real functions on $\mathbf{R}%
^{d}$. The following theorem establishes the local asymptotic power
functions of the test based on $\hat{T}_{n}.$\bigskip

\noindent \textsc{Theorem 5:} \textit{Suppose that Assumptions 1-3 hold and
that} $h\rightarrow 0$ \textit{and }$n^{-1/2}h^{-3d/2}\rightarrow 0$, 
\textit{as} $n\rightarrow \infty $.

\noindent (i) \textit{If }$p=1$\textit{, then} \textit{under} $H_{\delta
}^{\ast },\ $\textit{we have}%
\begin{equation*}
\lim_{n\rightarrow \infty }P\{\hat{T}_{n}>z_{1-\alpha }\}=1-\Phi
(z_{1-\alpha }-\eta _{2,-1}(w,\delta )/(\sqrt{2\pi }\sigma )).
\end{equation*}%
\noindent (ii) \textit{If }$p=2$\textit{, then under} $H_{\delta }^{\ast },\ 
$\textit{we have}%
\begin{equation*}
\lim_{n\rightarrow \infty }P\{\hat{T}_{n}>z_{1-\alpha }\}=1-\Phi
(z_{1-\alpha }-\eta _{2,0}(w,\delta )/\sigma ).
\end{equation*}%
\bigskip

Theorem 5 shows that the equality tests (on (\ref{eqtest})), in contrast to
the inequality tests (on (\ref{ineqtest})), have nontrivial local power
against alternatives converging to the null at rate $n^{-1/2}h^{-d/4},$
which is slower than $n^{-1/2}.$ This phenomenon of different convergence
rates arises because $\Lambda _{p}$ is symmetric around zero in the case of
equality tests, and it is not in the case of inequality tests. To see this
closely, observe that in the case of $p=1$, the power comparison between the
equality test and the inequality test is reduced to comparison between $%
\mathbf{E}|\mathbb{Z}_{1}+\mu |-\mathbf{E}|\mathbb{Z}_{1}|$ and $\mathbf{E}%
\max \{\mathbb{Z}_{1}+\mu ,0\}-\mathbf{E}\max \{\mathbb{Z}_{1},0\}$ for $\mu 
$ close to zero, where $\mathbb{Z}_{1}$ follows a standard normal
distribution with $\phi $ denoting its density. Note that we can approximate 
$\mathbf{E}|\mathbb{Z}_{1}+\mu |-\mathbf{E}|\mathbb{Z}_{1}|$ by $\{\phi
^{\prime \prime }(0)+2\phi (0)\}\mu ^{2}$ for $\mu $ close to zero, and
approximate $\mathbf{E}\max \{\mathbb{Z}_{1}+\mu ,0\}-\mathbf{E}\max \{%
\mathbb{Z}_{1},0\}$ by $\Phi (0)\mu $ for $\mu $ close to zero. The smaller
scale $\mu ^{2}$ in the former case arises because the leading term in the
expansion of $\mathbf{E}|\mathbb{Z}_{1}+\mu |-\mathbf{E}|\mathbb{Z}_{1}|$
around $\mu =0$ disappears due to the symmetry of the absolute value
function $|\cdot |$. Therefore, the different rate of convergence arises due
to our symmetric treatment of the alternative hypotheses (positive or
negative) in the equality test, in contrast to the asymmetric treatment in
the inequality test.

Since $\eta _{2,-1}(w,\delta )$ and $\eta _{2,0}(w,\delta )$ are always
nonnegative, the equality tests are \textit{locally asymptotically unbiased}
against any local alternatives. In contrast, the terms $\eta _{1,0}(w,\delta
)$ and $\eta _{1,1}(w,\delta )$ in the local asymptotic power functions of
the inequality tests in Theorem 4 can take negative values for some local
alternatives, implying that the inequality tests might be asymptotically
biased against such local alternatives. This feature is not due to the form
of our proposed inequality test, but is rather a common feature in testing
moment inequalities. It is because the null hypothesis is given by a
composite hypothesis and most of the powerful tests are not similar on the
boundary and hence biased against some local alternatives. In principle, one
can construct a test that is asymptotically similar on the boundary, but
such a test has typically poor power. See Andrews (2011) for details.

The test in Theorem 5 shares some features common in nonparametric tests
that are known to detect some smooth local alternatives that have narrow
peaks as the sample size increases. See e.g. Fan and Li (2000) and
references therein. To see this closely, consider a sequence of non-Pitman
local alternatives of type:%
\begin{equation*}
H_{\delta _{n}}^{\ast }:g_{j}(x)=\gamma _{n}\delta _{j,n}(x),\text{ for each 
}j\in \mathcal{J},
\end{equation*}%
where $\gamma _{n}$ is a deterministic sequence and $\delta _{j,n}(x)$ is
now allowed to change over $n$. For example, one may consider $\delta
_{j,n}(x)$ to be a function with a single peak that becomes sharper as $n$
becomes large, e.g. $\delta _{j,n}(x)=L_{j}\left( (x-x_{0})/\zeta
_{n}\right) ,$ where $L_{j}\left( \cdot \right) $ is a bounded function, $%
x_{0}\in \mathbf{R}^{d}$ is a fixed point, and $\zeta _{n}\rightarrow 0$ as $%
n\rightarrow \infty $. By using the same arguments as in the proof of
Theorem 5, we can show that the two-sided version of our test has nontrivial
power against such local alternatives provided $\lim_{n\rightarrow \infty
}nh^{d/2}\gamma _{n}^{2}\eta _{2,-1}(w,\delta _{j,n})\neq 0$ (for $p=1$) or $%
\lim_{n\rightarrow \infty }nh^{d/2}\gamma _{n}^{2}\eta _{2,0}(w,\delta
_{j,n})\neq 0$ (for $p=2$). However, since our main interest lies in testing
functional inequalities, we will not pursue further local power properties
of the equality test. On the other hand, it would also be interesting to see
whether it would give an adaptive, rate-optimal test to take the supremum of
our two-sided version of our test over a set of bandwidths, as in Horowitz
and Spokoiny (2001). However, the latter study is beyond of the scope of
this paper.

As in Section \ref{sec:lap}, when $J=1$, an optimal directional test under (%
\ref{hastar}) can also be obtained by following the arguments leading up to (%
\ref{wf}) so that%
\begin{equation*}
w_{p}^{\ast }(x)=\left\{ 
\begin{tabular}{ll}
$\frac{\delta _{1}^{2}(x)\rho _{1}^{-3}(x)}{\sqrt{\int_{\mathcal{X}}\delta
_{1}^{4}(x)\rho _{1}^{-4}(x)dx}},$ & $\text{ if }p=1,\text{ and}$ \\ 
$\frac{\delta _{1}^{2}(x)\rho _{1}^{-4}(x)}{\sqrt{\int_{\mathcal{X}}\delta
_{1}^{4}(x)\rho _{1}^{-4}(x)dx}},$ & $\text{ if }p=2.$%
\end{tabular}%
\right.
\end{equation*}%
Similarly as before, let the support of $\delta _{1}$ be contained in an $%
\varepsilon $-shrunk subset of $\mathcal{X}$. The optimal weight function
yields the following local power functions:%
\begin{equation*}
\begin{tabular}{ll}
$1-\Phi \left( z_{1-\alpha }-\frac{\sqrt{\int \delta _{1}^{4}(x)\rho
_{1}^{-4}(x)dx}}{\sqrt{2\pi \bar{q}_{1}}}\right) ,$ & $\text{ if }p=1,\text{
and}$ \\ 
$1-\Phi \left( z_{1-\alpha }-\frac{\sqrt{\int \delta _{1}^{4}(x)\rho
_{1}^{-4}(x)dx}}{\sqrt{\bar{q}_{2}}}\right) ,$ & $\text{ if }p=2,$%
\end{tabular}%
\end{equation*}%
where $\bar{q}_{p}\equiv \int_{\lbrack -1,1]^{d}}Cov(|\sqrt{1-t^{2}(u)}%
\mathbb{Z}_{1}+t(u)\mathbb{Z}_{2}|^{p},\ |\mathbb{Z}_{2}|^{p})du$, for $p\in
\{1,2\}$.

\section{Monte Carlo Experiments}

\label{sec:mc}

This section reports the finite-sample performance of the one-sided $L_{1}$-
and $L_{2}$-type tests from a Monte Carlo study. In the experiments, $n$
observations of a pair of random variables $(Y,X)$ were generated from $%
Y=m(X)+\sigma (X)U$, where $X\sim \text{Unif}[0,1]$ and $U\sim N(0,1)$ and $%
X $ and $U$ are independent. In all the experiments, we set $\mathcal{X}%
=[0.05,0.95]$.

To evaluate the finite-sample size of the tests, we first set $m(x)\equiv 0$%
. We call this case DGP0. In addition, we consider the following alternative
model 
\begin{equation}
m(x)=x(1-x)-c_{m}  \label{alt-model}
\end{equation}%
where $c_{m}\in \{0.25, 0.20, 0.15, 0.10, 0.05\}$. We call these 5 cases
DGPs 1-5. When $c_{m}=0.25$ (DGP1), we have $m(x)<0$ for all $x\neq 0.5$ and 
$m(x)=0$ with $x=0.5$. Hence, this case corresponds to the \textquotedblleft
interior\textquotedblright\ of the null hypothesis. In view of asymptotic
theory, we expect the empirical probability of rejecting $H_{0}$ to converge
to zero as $n$ gets large. When $c_{m} < 0.25$ (DGPs 2-5), we have $m(x)>0$
for some $x$. Therefore, these four cases are considered to see the
finite-sample power of our tests. Two different functions of $\sigma (x)$
are considered: $\sigma (x)\equiv 1$ (homoskedastic error) and $\sigma (x)=x$
(heteroskedastic error).

The experiments use sample sizes of $n=50,200,1000$ and the nominal level of 
$\alpha =0.05$. We performed 1000 Monte Carlo replications in each
experiment. In implementing both $L_{1}$ and $L_{2}$-type tests, we used $%
K(u)=(3/2)(1-(2u)^{2})I(|u|\leq 1/2)$ and $h=c_{h}\times \hat{s}_{X}\times
n^{-1/5}$, where $I(A)$ is the usual indicator function that has value one
if $A$ is true and zero otherwise, $c_{h}$ is a constant and $\hat{s}_{X}$
is the sample standard deviation of $X$. To check the sensitivity to the
choice of the bandwidth, eight different values of $c_{h}$ are considered: $%
\{ 0.75, 1, 1.25, 1.5, 1.75, 2, 2.25, 2.50 \}$. Finally, we considered the
uniform weight function: $w(x)=1$ and the inverse standard error weight
function: $w(x) = 1/\rho_n(x)$.

To evaluate the relative performance of our test, we have also implemented
one of test statistics proposed by Andrews and Shi (2011a), specifically
their Cram\'{e}r-von Mises-type (CvM) statistic with both plug-in asymptotic
(PA/Asy) and asymptotic generalized moment selection (GMS/Asy) critical
values. Specifically, countable hypercubes are used as instrument functions,
and tuning parameters were chosen, following suggestions as in Section 9 of
Andrews and Shi (2011a).

\begin{figure}[htbp]
\caption{Results of Monte Carlo Experiments: $L_1$ test and $\protect\sigma%
(x) \equiv 1$}
\label{figure1}
\begin{center}
\makebox{
\includegraphics[origin=bl,scale=.8,angle=180]{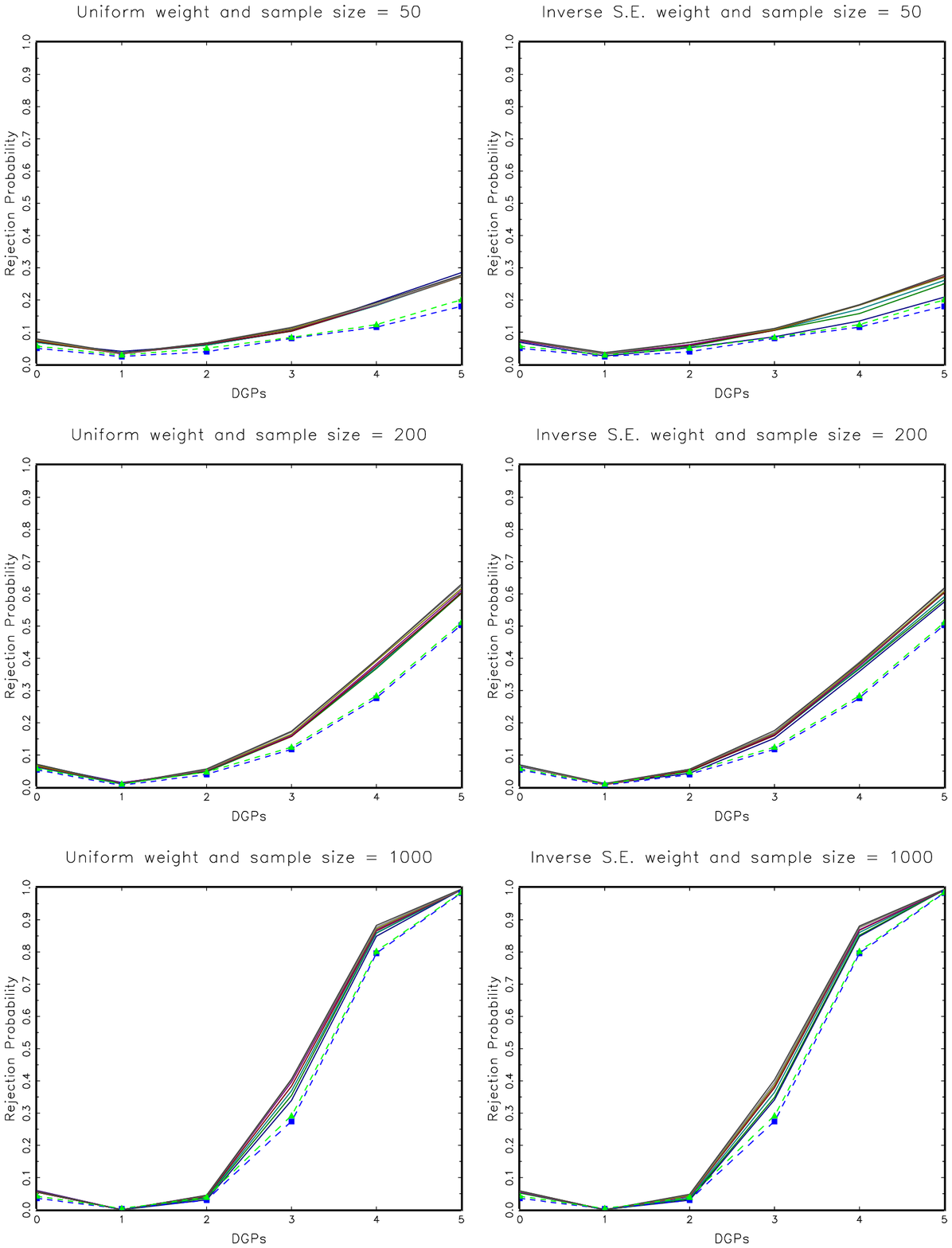}
}
\end{center}
\par
\parbox{5in}{Notes:  8 different solid lines in each panel correspond to our test with 8 different bandwidth values.
2 dotted lines correspond to the test of Andrews and Shi (2011a) with PA and GMS critical values.
The nominal level for each test is $\alpha = 0.05$. There are 1000 Monte Carlo replications in each
experiment.}
\end{figure}

Empirical rejection probabilities are plotted in Figures \ref{figure1}-\ref%
{figure4}. 8 different solid lines in each panel correspond to our test with
8 different bandwidth values. 2 dotted lines correspond to the test of
Andrews and Shi (2011a) with PA and GMS critical values. For each case, the
test with the GMS critical value gives slightly higher rejection
probabilities than that with the PA critical value. When $H_{0}$ is true and 
$m(x)\equiv 0$ (DGP0), the differences between the nominal and empirical
rejection probabilities are small. When $H_{0}$ is true and $m(x)$ is (\ref%
{alt-model}) with $c_{m}=0.25$ (the interior case DGP1), the empirical
rejection probabilities are smaller than the nominal level and become almost
zero for $n=1000$.

When $H_{0}$ is false and the correct model is (\ref{alt-model}) with $c_{m}
< 0.25$ (DGPs 2-5), the power of both the $L_{1}$ and $L_{2}$ tests is
increasing as $c_m$ gets smaller. This finding is consistent with asymptotic
theory since it is likely that our test will be more powerful when $\int_{%
\mathcal{X}}m(x)w(x)dx$ is larger. Note that in DGPs 3-5, ($c_{m}=0.15,
0.10, 0.05$), the rejection probabilities increase as $n$ gets large. This
is in line with the asymptotic theory in the preceding sections, for our
test is consistent for these values of $c_{m}$. However, the rejection
probabilities are quite small even with $n=1000$ for $c_m = 0.20$ (DGP 2).
This is not surprising given that our test can be biased, as shown in
Section \ref{sec:lap}. To further investigate the issue of bias associated
with $\int_{\mathcal{X}} m(x)w(x)dx$, we carried out an additional
simulation with $m(x) = \sin(2 \pi x)$. It turns out that rejection
probabilities were almost one across different values of the bandwidth for
both weight functions and for both homoskedastic and heteroskedastic errors.
This seems to be consistent with Theorem 4* in Section \ref{sec:lap}. We do
not report full details of additional simulation results for brevity.

Simulation results for the CvM statistics are similar to our test
statistics. More precisely, in Figure \ref{figure1} (the homoskedasticity
case), the $L_1$ test with both weight functions seems to be more powerful
than Andrews and Shi's test, whereas in Figure \ref{figure4}, their test
appears to be more powerful than the $L_2$ test with the uniform weight.
However, for most cases, power performances are comparable between each
other. Note further that there is little difference between PA and GMS
critical values for the CvM statistic of Andrews and Shi (2011a). This is
due to the fact that $m(x)$ is either flat or has a maximum at a single
point. We note that the results are not very sensitive to the bandwidth
choice for our tests. Finally, regarding the choice of the weight function,
we would like to recommend the inverse standard error weight since it seems
to perform better than the uniform weight in simulations.

\section{Proofs}

\label{sec:proof}

This section begins with a roadmap for the proof, where the roles of
technical lemmas and main difficulties are explained. Then we state the
lemmas and present the proofs of the theorems. 

\subsection{The Roadmap for the Proof of Theorem 1}

\label{sec:appendix:roadmap}

The proof of Theorem 1 follows the structure of the proof of the
finite-dimensional convergence in Theorem 1.1 of GMZ.

Under the condition of Theorem 1 that $m_{j}(x)=0$ for almost all $x\in 
\mathcal{X}$ and for all $j\in \mathcal{J}$, we can show that $\mathbf{E}%
\hat{g}_{jn}(x)=0$ for almost all $x$ in the support of $w_{j}$ from some
large $n$ on. This means that by letting $v_{jn}(x)\equiv \hat{g}_{jn}(x)-%
\mathbf{E}\hat{g}_{jn}(x)$ and $\zeta _{n}(A)\equiv
\sum_{j=1}^{J}\int_{A}\Lambda _{p}(v_{jn}(x))w_{j}(x)dx$ with some $A\subset 
\mathcal{X}$, we can write $T_{n}$ as%
\begin{eqnarray}
&&\frac{n^{p/2}h^{(p-1)d/2}}{\sigma _{n}}\left\{ \zeta _{n}(\mathcal{X}%
\backslash A)-\mathbf{E}\zeta _{n}(\mathcal{X}\backslash A)\right\}
\label{dec} \\
&&+\frac{n^{p/2}h^{(p-1)d/2}}{\sigma _{n}}\left\{ \zeta _{n}(A)-\mathbf{E}%
\zeta _{n}(A)\right\} .  \notag
\end{eqnarray}%
The main part of the proof of Theorem 1 establishes asymptotic normality for
the second term and asymptotic negligibility for the first term when $A$ is
chosen to nearly cover $\mathcal{X}$. The proof of asymptotic normality
employs the Poissonization method of GMZ which prevents us from choosing $A$
to cover $\mathcal{X}$ entirely. This makes the proof intricate. The
asymptotic arguments for both terms of (\ref{dec}) require that $\sigma _{n}$
is an asymptotically stable quantity. Hence we begin by dealing with $\sigma
_{n}$.\bigskip

\noindent \textit{Step 1:} In Lemma A7, we show that given appropriate $%
A\subset \mathcal{X}$, $\sigma _{n}(A)\rightarrow \sigma (A)>0$ as $%
n\rightarrow \infty $, for some $\sigma (A)>0$, where $\sigma _{n}(A)$ is $%
\sigma _{n}$ except that the integral domains of $\sigma _{jk,n}$ are
restricted to $A$. To prove the convergence, we choose the domain $A$ to be
such that nonparametric functions $\rho _{jk,n}$ that constitute $\sigma
_{n} $ are continuous and uniformly convergent on this domain. That we can
choose such $A$ to be large enough is ensured by Lemma A1. The proof is
lengthy, the main step being the approximation of covariances of Poissonized
sums. For this approximation, we use a type of a Berry-Esseen bound for sums
of independent random variables due to Sweeting (1977). This bound is
restated in Lemma A2. Since the bound involves various moments of random
quantities, we prepare these moment bounds in Lemmas A4 and A5.\bigskip

\noindent \textit{Step 2:} We establish that the second term in (\ref{dec})
is asymptotically standard normal when $A$ nearly covers $\mathcal{X}$.
First, we use Lemma A6 to show that the second component in (\ref{dec}) is
asymptotically equivalent to 
\begin{equation}
\frac{n^{p/2}h^{(p-1)d/2}}{\sigma _{n}}\left\{ \zeta _{n}(A)-\mathbf{E}\zeta
_{N}(A)\right\} ,  \label{ns}
\end{equation}%
where $\zeta _{N}(A)\equiv \sum_{j=1}^{J}\int_{A}\Lambda
_{p}(v_{jN}(x))w_{j}(x)dx$, $v_{jN}(x)\equiv \hat{g}_{jN}(x)-\mathbf{E}\hat{g%
}_{jn}(x)$,%
\begin{equation}
\hat{g}_{jN}(x)\equiv \frac{1}{nh^{d}}\sum_{i=1}^{N}Y_{ji}K\left( \frac{%
x-X_{i}}{h}\right) ,  \label{gN}
\end{equation}%
and $N$ is a Poisson random variable with mean $n$ and independent of all
the other random variables. Then consider%
\begin{equation}
S_{n}(A)\equiv \frac{n^{p/2}h^{(p-1)d/2}\{\zeta _{N}(A)-\mathbf{E}\zeta
_{N}(A)\}}{\sigma _{n}(A)},  \label{Sn}
\end{equation}%
where $\sigma _{n}^{2}(A)\equiv \sum_{j=1}^{J}\sum_{k=1}^{J}\sigma
_{jk,n}(A) $ and $\sigma _{jk,n}(A)$ is $\sigma _{jk,n}$ with the integral
domain restricted to $A$. Note that the numerator of $S_{n}(A)$ is based on
the Poissonized version $v_{jN}(x)$ so that when we cut the integral in $%
\zeta _{N}(A)$ into integrals on small disjoint domains and sum them, this
latter sum behaves like a sum of independent random variables. In Lemma A9,
we construct this sum and apply the CLT to obtain asymptotic normality for $%
S_{n}(A)$. Then in Lemma A10, using the de-Poissonization lemma of Beirlant
and Mason (1995), we deduce that the conditional distribution of $S_{n}(A)$
given $N=n$ converges to a standard normal distribution. (This lemma
requires the set $\mathcal{X}\backslash A$ to stay nonempty.) Since this
conditional distribution is nothing but the distribution of (\ref{ns}), we
conclude that the second term in (\ref{dec}) is asymptotically standard
normal. However, this sequence of arguments so far presumes that $\sigma
_{n} $ is an asymptotically right scale, which means that $\sigma _{n}$
should be based on the Poissonized version $v_{jN}(x)$ not on the original
one $v_{jn}(x)$.\bigskip

\noindent \textit{Step 3:} It remains to deal with the first term in (\ref%
{dec}). Since $\sigma _{n}(A)$ is close to $\sigma (A)>0$ by Step 1 for
large samples, it suffices to show that the quantity 
\begin{equation*}
n^{p/2}h^{(p-1)d/2}\left\{ \zeta _{n}(\mathcal{X}\backslash A)-\mathbf{E}%
\zeta _{n}(\mathcal{X}\backslash A)\right\}
\end{equation*}%
is asymptotically negligible for large $n$ and large set $A$. This is
accomplished by Lemma A8, which again uses moment bounds of Lemmas A4 and
A5. Since $w_{j}$ is square integrable, if we can take $A\subset \mathcal{X}$
such that $\int_{\mathcal{X}/A}w_{j}^{2}(x)dx$ is small, the asymptotic
negligibility of the first component in (\ref{dec}) follows by Lemma A8.
Lemma 8 extends to Lemma 6.2 of GMZ from $p=1$ to $p\geq 1$. This
generalization is necessary since the majorization inequality of Pinelis
(1994) used in GMZ is not directly applicable in the general case with $p
\geq 1$. \bigskip

\noindent \textit{Step 4:} Finally, we approximate $\mathbf{E}\zeta _{n}(A)$
in the second component in (\ref{dec}) by an estimable quantity, $\sum_{j\in 
\mathcal{J}}a_{jn}$ in Theorem 1. This step is done through Lemma A6. The
lemma is adapted from Lemma 6.3 of GMZ, but unlike their case of $L_{1}$%
-norm, our case involves the one-sided $L_{p}$-norm with $p\geq 1$. For this
modification, we use the algebraic inequality of Lemma A3. This closes the
proof of Theorem 1.

\subsection{Technical Lemmas and the Proof of Theorem 1}

\label{sec:appendix-lemmas}

We begin with technical lemmas. The lemmas are ordered so that lemmas that
come later rely on their preceding lemmas. 

The first statement of the lemma below is a special case of Theorem 2(b) of
Stein (1970) on pages 62 and 63. The second statement is an extension of
Lemma 6.1 of GMZ.\bigskip

\noindent \textsc{Lemma A1:}\textbf{\ }\textit{Let} $J(\cdot ):\mathbf{R}%
^{d}\rightarrow \mathbf{R}$ \textit{be a Lebesgue integrable bounded function%
}\ \textit{and} $H:\mathbf{R}^{d}\rightarrow \mathbf{R}$ \textit{be a
bounded function with compact support} $S$. \textit{Then}, \textit{for
almost every }$y\in \mathbf{R}^{d},$%
\begin{equation*}
\int_{\mathbf{R}^{d}}J(x)H_{h}\left( y-x\right) dx\rightarrow
J(y)\int_{S}H(x)dx,\text{ \textit{as} }h\rightarrow 0\text{,}
\end{equation*}%
\textit{where} $H_{h}(x)\equiv H(x/h)/h^{d}.$

\textit{Furthermore}, \textit{suppose that} $\bar{J}\equiv \int |J(z)|dz>0$. 
\textit{Then for all }$0<\varepsilon <\bar{J}\equiv \int |J(z)|dz$, \textit{%
there exist} $M>0,$ $v >0$ \textit{and a Borel set} $B$ \textit{of finite
Lebesgue measure} $m(B)$ \textit{such that} $B\subset \lbrack -M+v,M-v]^{d}$%
, $\alpha \equiv \int_{\mathbf{R}^{d}\backslash \lbrack -M,M]^{d}}|J(z)|dz>0$%
, $\int_{B}|J(z)|dz>\bar{J}-\varepsilon $, $J$ \textit{is continuous on} $B,$
\textit{and}%
\begin{equation*}
\sup_{y\in B}\left\vert \int_{\mathbf{R}^{d}}J(x)H_{h}\left( y-x\right)
dx-J(y)\int_{S}H(x)dx\right\vert \rightarrow 0,\text{ \textit{as} }%
h\rightarrow 0\text{.}
\end{equation*}

\noindent \textsc{Proof:} The first statement is a special case of Theorem
2(b) of Stein (1970) on pages 62 and 63. The second statement can be proved
following the proof of Lemma 6.1 of GMZ. Since $J$ is Lebesgue integrable,
the integral $\int_{\mathbf{R}^{d}\backslash \lbrack -M,M]^{d}}|J(z)|dz$ is
continuous in $M$ and converges to zero as $M\rightarrow \infty .$ We can
find $M>0$ and $v>0$ such that 
\begin{equation*}
\int_{\mathbf{R}^{d}\backslash \lbrack -M,M]^{d}}|J(z)|dz=\varepsilon /8%
\text{ and }\int_{\mathbf{R}^{d}\backslash \lbrack
-M+v,M-v]^{d}}|J(z)|dz=\varepsilon /4.
\end{equation*}%
The construction of the desired set $B\subset \lbrack -M+v,M-v]^{d}$ can be
done using the arguments in the proof of Lemma 6.1 of GMZ. $\blacksquare $%
\bigskip

The following result is a special case of Theorem 1 of Sweeting (1977) with $%
g(x)=\min (x,1)$ (in his notation). See also Fact 6.1 of GMZ and Fact 4 of
Mason (2009) for applications of Theorem 1 of Sweeting (1977). \bigskip

\noindent \textsc{Lemma A2} \textsc{(Sweeting (1977)):} \textit{Let} $%
\mathbb{Z}\in \mathbf{R}^{k}$ \textit{be a mean zero normal random vector
with covariance matrix }$I$ \textit{and }$\{W_{i}\}_{i=1}^{n}$ \textit{is a
set of i.i.d. random vectors in} $\mathbf{R}^{k}$ \textit{such that }$%
\mathbf{E}W_{i}=0,$ $\mathbf{E}W_{i}W_{i}^{\prime }=I,$ \textit{and} $%
\mathbf{E}||W_{i}||^{r}<\infty $, $r\geq 3.$ \textit{Then for any Borel
measurable function }$\varphi :\mathbf{R}^{k}\rightarrow \mathbf{R}$ \textit{%
such that}%
\begin{equation*}
\sup_{x\in \mathbf{R}^{k}}\frac{\left\vert \varphi (x)-\varphi
(0)\right\vert }{1+||x||^{r}\min (||x||,1)}<\infty ,
\end{equation*}%
\textit{we have}%
\begin{eqnarray*}
&&\left\vert \mathbf{E}\left[ \varphi \left( \frac{1}{\sqrt{n}}%
\sum_{i=1}^{n}W_{i}\right) \right] -\mathbf{E}\left[ \varphi (\mathbb{Z})%
\right] \right\vert \\
&\leq &c_{1}\left( \sup_{x\in \mathbf{R}^{k}}\frac{\left\vert \varphi
(x)-\varphi (0)\right\vert }{1+||x||^{r}\min (||x||,1)}\right) \left\{ \frac{%
1}{\sqrt{n}}\mathbf{E}||W_{i}||^{3}+\frac{1}{n^{(r-2)/2}}\mathbf{E}%
||W_{i}||^{r}\right\} \\
&&+c_{2}\mathbf{E}\left[ \omega _{\varphi }\left( \mathbb{Z};\frac{c_{3}}{%
\sqrt{n}}\mathbf{E}||W_{i}||^{3}\right) \right] ,
\end{eqnarray*}%
\textit{where }$c_{1},\ c_{2}$ \textit{and }$c_{3}$\textbf{\ }\textit{are
positive constants that depend only on }$k$ \textit{and} $r$ \textit{and}%
\begin{equation*}
\omega _{\varphi }(x;\varepsilon )\equiv \sup \left\{ |\varphi (x)-\varphi
(y)|:y\in \mathbf{R}^{k},||x-y||\leq \varepsilon \right\} .
\end{equation*}

The following algebraic inequality is used frequently throughout the
proofs.\bigskip

\noindent \textsc{Lemma A3:} \textit{For any }$a,b\in \mathbf{R}$,\textit{\
let }$a_{+}=\max (a,0)$\textit{\ and} $b_{+}=\max (b,0)$\textit{.
Furthermore, for any real }$a\geq 0$\textit{, if }$a=0$,\textit{\ we define }%
$\lceil a\rceil =1$, \textit{and if }$a>0,$ \textit{we define} $\lceil
a\rceil $\textit{\ to be the smallest integer greater than or equal to} $a$%
\textit{. Then for any }$p\geq 1,$%
\begin{eqnarray*}
\max \left\{ |a_{+}^{p}-b_{+}^{p}|,||a|^{p}-|b|^{p}|\right\} &\leq
&2p|a-b|\left( \sum_{k=0}^{\lceil p-1\rceil }\frac{\lceil p-1\rceil !}{k!}%
|a-b|^{\lceil p-1\rceil -k}|b|^{k}\right) ^{(p-1)/\lceil p-1\rceil } \\
&\leq &C\sum_{k=0}^{\lceil p-1\rceil }|a-b|^{p-\frac{(p-1)k}{\lceil
p-1\rceil }}|b|^{k},
\end{eqnarray*}%
\textit{for some }$C>0$\textit{\ that depends only on} $p$.\bigskip

\noindent \textsc{Proof} \textsc{:} First, we show the inequality for the
case where $p$ is a positive integer. We prove first that $||a|^{p}-|b|^{p}|$
has the desired bound. Note that in this case of $p$ being a positive
integer, the bound takes the following form:%
\begin{equation*}
2\sum_{k=0}^{p-1}\frac{p!}{k!}|a-b|^{p-k}|b|^{k}.
\end{equation*}%
When $p=1$, the bound is trivially obtained. Suppose now that the inequality
holds for a positive integer $q$. First, note that using the mean-value
theorem, convexity of the function $f(x)=|x|^{q}$ for $q\geq 1$, and the
triangular inequality,%
\begin{eqnarray*}
||a|^{q+1}-|b|^{q+1}| &\leq &(q+1)|a-b|\text{sup}_{\alpha \in \lbrack
0,1]}\left( \alpha |a|+(1-\alpha )|b|\right) ^{q} \\
&\leq &(q+1)|a-b|\text{sup}_{\alpha \in \lbrack 0,1]}\left( \alpha
|a|^{q}+(1-\alpha )|b|^{q}\right) \\
&\leq &(q+1)|a-b|\left( ||a|^{q}-|b|^{q}|+2|b|^{q}\right) .
\end{eqnarray*}%
As for $||a|^{q}-|b|^{q}|$, we apply the inequality to bound the last term by%
\begin{eqnarray*}
&&(q+1)|a-b|\left( 2\sum_{k=0}^{q-1}\frac{q!}{k!}|a-b|^{q-k}|b|^{k}+|b|^{q}%
\right) \\
&=&2\sum_{k=0}^{q}\frac{(q+1)!}{k!}|a-b|^{q-k+1}|b|^{k}.
\end{eqnarray*}%
Therefore, by the principle of mathematical induction, the desired bound in
the case of $p$ being a positive integer follows.

Certainly, we obtain the same bound for $|a_{+}^{p}-b_{+}^{p}|$ when $p=1$.
When $p>1$, we observe that by the mean-value theorem,%
\begin{eqnarray*}
|a_{+}^{p}-b_{+}^{p}| &\leq &p|a-b|\left( |a|^{p-1}+|b|^{p-1}\right) \\
&\leq &p|a-b|\left( ||a|^{p-1}-|b|^{p-1}|+2|b|^{p-1}\right) .
\end{eqnarray*}%
By applying the previous inequality to $||a|^{p-1}-|b|^{p-1}|$, we obtain
the desired bound for $|a_{+}^{p}-b_{+}^{p}|$ when $p$ is a positive integer.

Since the bound holds for any positive integer $p$, let us consider the case
where $p$ is a real number strictly larger than $1$. Again, we first show
that $||a|^{p}-|b|^{p}|$ has the desired bound. Using the mean-value theorem
as before and the fact that $|a+b|\leq 2^{1-1/s}\left(
|a|^{s}+|b|^{s}\right) ^{1/s}$ for all $s\in \lbrack 1,\infty )$ and all $%
a,b\in \mathbf{R}$, we find that for $u\equiv \lceil p-1\rceil ,$%
\begin{eqnarray*}
||a|^{p}-|b|^{p}| &\leq &p\left\vert a-b\right\vert (|a|^{p-1}+|b|^{p-1}) \\
&\leq &p\left\vert a-b\right\vert 2^{1-(p-1)/u}\left( |a|^{u}+|b|^{u}\right)
^{(p-1)/u} \\
&\leq &p\left\vert a-b\right\vert 2^{1-(p-1)/u}\left( \left\vert
|a|^{u}-|b|^{u}\right\vert +2|b|^{u}\right) ^{(p-1)/u}.
\end{eqnarray*}%
Since $u$ is a positive integer, using the previous bound, we bound the
right-hand side by%
\begin{equation*}
p\left\vert a-b\right\vert 2^{1-(p-1)/u}\left( 2\sum_{k=0}^{u-1}\frac{u!}{k!}%
|a-b|^{u-k}|b|^{k}+2|b|^{u}\right) ^{(p-1)/u}.
\end{equation*}%
Consolidating the sum in the parentheses, we obtain the wanted bound.

As for the second inequality, observe that%
\begin{eqnarray*}
&&2p|a-b|\left( \sum_{k=0}^{\lceil p-1\rceil }\frac{\lceil p-1\rceil !}{k!}%
|a-b|^{\lceil p-1\rceil -k}|b|^{k}\right) ^{(p-1)/\lceil p-1\rceil } \\
&\leq &C\max_{k\in \{0,1,\cdot \cdot \cdot ,\lceil p-1\rceil
\}}|a-b|^{p-k\{(p-1)/\lceil p-1\rceil \}}|b|^{k}\leq C\sum_{k=0}^{\lceil
p-1\rceil }|a-b|^{p-k\{(p-1)/\lceil p-1\rceil \}}|b|^{k},
\end{eqnarray*}%
for some $C>0$ that depends only on $p$. We can obtain the same bound for $%
|a_{+}^{p}-b_{+}^{p}|$ by noting that $|a_{+}^{p}-b_{+}^{p}|\leq p\left\vert
a-b\right\vert \left( |a|^{p-1}+|b|^{p-1}\right) $ and following the same
arguments afterwards as before. $\blacksquare $\bigskip

Define for $j\in \mathcal{J},$%
\begin{equation}
k_{jn,r}(x)\equiv h^{-d}\mathbf{E}\left[ \left\vert Y_{ji}K\left( \frac{%
x-X_{i}}{h}\right) \right\vert ^{r}\right] ,\ r\geq 1.  \label{kr}
\end{equation}

\noindent \textsc{Lemma A4:}\textbf{\ }\textit{Suppose that Assumptions
1(i)(iii) and 2 hold and }$h\rightarrow 0$ \textit{as} $n\rightarrow \infty $%
. \textit{Then for }$\varepsilon >0$\textit{\ in Assumption 1(i), there
exist positive integer }$n_{0}$ \textit{and constants} $c_{1},c_{2}>0$ 
\textit{such that for all }$n\geq n_{0},$ \textit{all} $r\in \lbrack
1,2p+2], $ \textit{and all} $j\in \mathcal{J},$%
\begin{gather*}
0<c_{1}\leq \text{inf}_{x\in \mathcal{S}_{j}^{\varepsilon /2}}\rho
_{jn}^{2}(x)\text{ \textit{and}} \\
\text{sup}_{x\in \mathcal{S}_{j}^{\varepsilon /2}}k_{jn,r}(x)\leq
c_{2}<\infty .
\end{gather*}

\noindent \textsc{Proof:}\textbf{\ } Since $h\rightarrow 0$ as $n\rightarrow
\infty ,$ we apply change of variables to find that from large $n$ on,%
\begin{eqnarray*}
\inf_{x\in \mathcal{S}_{j}^{\varepsilon /2}}\rho _{jn}^{2}(x) &=&\inf_{x\in 
\mathcal{S}_{j}^{\varepsilon /2}}\frac{1}{h^{d}}\mathbf{E}\left[
Y_{ji}^{2}K^{2}\left( \frac{x-X_{i}}{h}\right) \right] \\
&\geq &\inf_{x\in \mathcal{S}_{j}^{\varepsilon }}\mathbf{E}\left[
Y_{ji}^{2}|X=x\right] f(x)\int_{[-1/2,1/2]^{d}}K^{2}\left( u\right) du>c_{1},
\end{eqnarray*}%
for some $c_{1}>0$ by Assumptions 1(i) and 2. Similarly, from some large $n$
on,%
\begin{equation*}
\sup_{x\in \mathcal{S}_{j}^{\varepsilon /2}}k_{jn,r}(x)\leq \sup_{x\in 
\mathcal{S}_{j}^{\varepsilon }}\mathbf{E}\left[ |Y_{ji}|^{r}|X=x\right]
f(x)\int |K\left( u\right) |^{r}du<\infty ,
\end{equation*}%
by Assumptions 1(i) and 2. $\blacksquare $\bigskip

Define for each $j\in \mathcal{J},$%
\begin{equation*}
\hat{g}_{jN}(x)\equiv \frac{1}{nh^{d}}\sum_{i=1}^{N}Y_{ji}K\left( \frac{%
x-X_{i}}{h}\right) ,\text{ }x\in \mathcal{X}\text{,}
\end{equation*}%
where $N$ is a Poisson random variable that is common across $j\in \mathcal{J%
}$, has mean $n$, and is independent of $\{(Y_{ji},X_{i}):j\in \mathcal{J}%
\}_{i=1}^{\infty }$. Let for each $j\in \mathcal{J},$%
\begin{equation*}
v_{jn}(x)\equiv \hat{g}_{jn}(x)-\mathbf{E}\hat{g}_{jn}(x)\text{, and\ }%
v_{jN}(x)\equiv \hat{g}_{jN}(x)-\mathbf{E}\hat{g}_{jn}(x).
\end{equation*}%
We define, for each $j\in \mathcal{J},$%
\begin{eqnarray}
\xi _{jn}(x) &\equiv &\frac{\sqrt{nh^{d}}v_{jN}(x)}{\rho _{jn}(x)}\text{ and}
\label{eq2} \\
V_{jn}(x) &\equiv &\frac{\sum_{i\leq N_{1}}\left\{ Y_{ji}K\left(
(x-X_{i})/h\right) -\mathbf{E}\left( Y_{ji}K\left( (x-X_{i})/h\right)
\right) \right\} }{\sqrt{\mathbf{E}[Y_{ji}^{2}K^{2}\left( (x-X_{i})/h\right)
]}},  \notag
\end{eqnarray}%
where $N_{1}$ denotes a Poisson random variable with mean $1$ that is
independent of $\{(Y_{ji},X_{i}):j\in \mathcal{J}\}_{i=1}^{\infty }$. Then, $%
Var(V_{jn}(x))=1.$ Let $V_{jn}^{(i)}(x),\ i=1,\cdot \cdot \cdot ,n,$ be
i.i.d. copies of $V_{jn}(x)$ so that%
\begin{equation}
\xi _{jn}(x)\overset{d}{=}\frac{1}{\sqrt{n}}\sum_{i=1}^{n}V_{jn}^{(i)}(x).
\label{eq3}
\end{equation}

\noindent \textsc{Lemma A5:}\textbf{\ }\textit{Suppose that Assumptions
1(i)(iii) and 2 hold and }$h\rightarrow 0$ \textit{as} $n\rightarrow \infty $
\textit{and}%
\begin{equation*}
\text{limsup}_{n\rightarrow \infty }n^{-r/2+1}h^{(1-r/2)d}<C,
\end{equation*}%
\textit{for some constant }$C>0$ \textit{and for }$r\in \lbrack 2,2p+2]$. 
\textit{Then, for }$\varepsilon >0$ \textit{in Assumption 1(iii)},%
\begin{equation*}
\sup_{x\in \mathcal{S}_{j}^{\varepsilon /2}}\mathbf{E}\left[ |V_{jn}(x)|^{r}%
\right] \leq C_{1}h^{(1-r/2)d}\text{ \textit{and\ }}\sup_{x\in \mathcal{S}%
_{j}^{\varepsilon /2}}\mathbf{E}\left[ |\xi _{jn}(x)|^{r}\right] \leq C_{2},
\end{equation*}%
\textit{where }$C_{1}$ \textit{and} $C_{2}$\textit{\ are constants that
depend only on }$r$.\bigskip

\noindent \textsc{Proof} \textsc{:}\textbf{\ }For all $x\in \mathcal{S}%
_{j}^{\varepsilon /2}$,\ $\mathbf{E}[V_{jn}^{2}(x)]=1.$ Recall the
definition of $k_{jn,r}(x)$ in (\ref{kr}). Then for some $C_{0},C_{1}>0,$%
\begin{equation}
\sup_{x\in \mathcal{S}_{j}^{\varepsilon /2}}\mathbf{E}\left\vert
V_{jn}(x)\right\vert ^{r}\leq C_{0}\sup_{x\in \mathcal{S}_{j}^{\varepsilon
/2}}\frac{h^{d}k_{jn,r}(x)}{h^{rd/2}\rho _{jn}^{r}(x)}\leq C_{1}h^{(1-r/2)d},
\label{bd6}
\end{equation}%
by Lemma A4, completing the proof of the first statement.

As for the second statement, using (\ref{eq3}) and applying Rosenthal's
inequality (e.g. (2.3) of GMZ), we deduce that for positive constants $%
C_{3},C_{4}$ and $C_{5}$ that depend only on $r$,%
\begin{eqnarray*}
\sup_{x\in \mathcal{S}_{j}^{\varepsilon /2}}\mathbf{E}\left[ |\xi
_{jn}(x)|^{r}\right] &\leq &C_{3}\sup_{x\in \mathcal{S}_{j}^{\varepsilon
/2}}\max \{(\mathbf{E}V_{jn}^{2}(x))^{r/2},n^{-r/2+1}\mathbf{E}%
|V_{jn}(x)|^{r}\} \\
&\leq &C_{4}\max \left\{ 1,C_{5}n^{-r/2+1}h^{(1-r/2)d}\right\}
\end{eqnarray*}%
by (\ref{bd6}). By the condition that limsup$_{n\rightarrow \infty
}n^{-r/2+1}h^{(1-r/2)d}<C$, the desired result follows. $\blacksquare $%
\bigskip

The following lemma is adapted from Lemma 6.3 of GMZ. The result is obtained
by combining Lemmas A2-A5.\bigskip

\noindent \textsc{Lemma A6:}\textbf{\ }\textit{Suppose that Assumptions 1
and 2 hold and }$h\rightarrow 0$ \textit{and} $n^{-1/2}h^{-d}\rightarrow 0$ 
\textit{as} $n\rightarrow \infty $. \textit{Then for any Borel set} $%
A\subset \mathbf{R}^{d}$ \textit{and for any} $j\in \mathcal{J},$%
\begin{eqnarray*}
\int_{A}\left\{ n^{p/2}h^{(p-1)d/2}\mathbf{E}\Lambda
_{p}(v_{jN}(x))-h^{-d/2}\rho _{jn}^{p}(x)\mathbf{E}\Lambda _{p}(\mathbb{Z}%
_{1})\right\} w_{j}(x)dx &\rightarrow &0, \\
\int_{A}\left\{ n^{p/2}h^{(p-1)d/2}\mathbf{E}\Lambda
_{p}(v_{jn}(x))-h^{-d/2}\rho _{jn}^{p}(x)\mathbf{E}\Lambda _{p}(\mathbb{Z}%
_{1})\right\} w_{j}(x)dx &\rightarrow &0.
\end{eqnarray*}

\noindent \textsc{Proof} \textsc{:}\textbf{\ }Recall the definition of $\xi
_{jn}(x)$ in (\ref{eq2}) and write%
\begin{eqnarray*}
&&n^{p/2}h^{(p-1)d/2}\mathbf{E}\Lambda _{p}(v_{jN}(x))-h^{-d/2}\rho
_{jn}^{p}(x)\mathbf{E}\Lambda _{p}(\mathbb{Z}_{1}) \\
&=&h^{-d/2}\rho _{jn}^{p}(x)\left\{ \mathbf{E}\Lambda _{p}(\xi _{jn}(x))-%
\mathbf{E}\Lambda _{p}(\mathbb{Z}_{1})\right\} .
\end{eqnarray*}%
In view of Lemma A4 and Assumption 1(ii), we find that it suffices for the
first statement of the lemma to show that%
\begin{equation}
\sup_{x\in \mathcal{S}_{j}}\left\vert \mathbf{E}\Lambda _{p}(\xi _{jn}(x))-%
\mathbf{E}\Lambda _{p}(\mathbb{Z}_{1})\right\vert =o(h^{d/2}).  \label{claim}
\end{equation}%
By Lemma A5,\ $\sup_{x\in \mathcal{S}_{j}}\mathbf{E}\left\vert
V_{jn}(x)\right\vert ^{3}\leq Ch^{-d/2}$ for some $C>0$. Using Lemma A2 and
taking $r=\max \{p,3\}$ and $V_{jn}^{(i)}(x)=W_{i}$, and $\Lambda _{p}(\cdot
)=\varphi (\cdot ),$ we deduce that%
\begin{eqnarray}
&&\sup_{x\in \mathcal{S}_{j}}\left\vert \mathbf{E}\Lambda _{p}(\xi _{jn}(x))-%
\mathbf{E}\Lambda _{p}(\mathbb{Z}_{1})\right\vert  \label{arg} \\
&\leq &C_{1}n^{-1/2}\sup_{x\in \mathcal{S}_{j}}\mathbf{E}\left\vert
V_{jn}(x)\right\vert ^{3}+C_{2}n^{-(r-2)/2}\sup_{x\in \mathcal{S}_{j}}%
\mathbf{E}\left\vert V_{jn}(x)\right\vert ^{r}  \notag \\
&&+C_{3}\sup_{x\in \mathcal{S}_{j}}\mathbf{E}\left[ \omega _{\Lambda
_{p}}\left( \mathbb{Z}_{1};\frac{C_{4}}{\sqrt{n}}\mathbf{E}\left\vert
V_{jn}(x)\right\vert ^{3}\right) \right] ,  \notag
\end{eqnarray}%
for some constants $C_{s}>0,$ $s=1,2,3.$ The first two terms are $o(h^{d/2})$%
. As for the last expectation, observe that by Lemma A3,%
\begin{equation*}
\mathbf{E}\left[ \omega _{\Lambda _{p}}\left( \mathbb{Z}_{1};\frac{C_{4}}{%
\sqrt{n}}\mathbf{E}\left\vert V_{jn}(x)\right\vert ^{3}\right) \right] \leq
C\sum_{k=0}^{\lceil p-1\rceil }\left( \frac{C_{4}}{\sqrt{n}}\mathbf{E}%
\left\vert V_{jn}(x)\right\vert ^{3}\right) ^{p-\frac{(p-1)k}{\lceil
p-1\rceil }}\mathbf{E}|\mathbb{Z}_{1}|^{k}.
\end{equation*}%
The last sum is $O(n^{-1/2}h^{-d/2})=o(h^{d/2})$ uniformly over $x\in 
\mathcal{S}_{j}$, completing the proof of (\ref{claim}).

We consider the second statement. Let $\bar{V}_{jn}^{(k)}(x),\ k=1,\cdot
\cdot \cdot ,n,$ be i.i.d. copies of%
\begin{equation*}
\frac{Y_{ji}K\left( \frac{x-X_{i}}{h}\right) -\mathbf{E}\left( Y_{ji}K\left( 
\frac{x-X_{i}}{h}\right) \right) }{\sqrt{\mathbf{E}\left[ Y_{ji}^{2}K^{2}%
\left( \frac{x-X_{i}}{h}\right) \right] -\left( \mathbf{E}\left[
Y_{ji}K\left( \frac{x-X_{i}}{h}\right) \right] \right) ^{2}}}
\end{equation*}%
so that $Var(\bar{V}_{jn}^{(k)}(x))=1$. Observe that for some constants $%
C_{1},C_{2}>0$,%
\begin{equation}
\sup_{x\in \mathcal{S}_{j}}\mathbf{E}\left\vert \bar{V}_{jn}^{(k)}(x)\right%
\vert ^{3}\leq Ch^{-d/2}\sup_{x\in \mathcal{S}_{j}}\frac{k_{jn,3}(x)}{\left(
\rho _{jn}^{2}(x)-h^{d}b_{jn}^{2}(x)\right) ^{3/2}}\leq C_{2}h^{-d/2},
\label{bd}
\end{equation}%
where $b_{jn}(x)\equiv h^{-d}\mathbf{E}\left[ Y_{ji}K\left(
(x-X_{i})/h\right) \right] $. The last inequality follows by Lemma A4. Define%
\begin{equation*}
\bar{\xi}_{jn}(x)\equiv \frac{\sqrt{nh^{d}}v_{jn}(x)}{\tilde{\rho}_{jn}(x)},
\end{equation*}%
where $\tilde{\rho}_{jn}^{2}(x)\equiv nh^{d}Var(v_{jn}(x))$. Then $\bar{\xi}%
_{jn}(x)\overset{d}{=}\frac{1}{\sqrt{n}}\sum_{k=1}^{n}\bar{V}_{jn}^{(k)}(x)$%
. Using Lemma A2 and following the arguments in (\ref{arg}) analogously, we
deduce that%
\begin{equation*}
\sup_{x\in \mathcal{S}_{j}}\left\vert \mathbf{E}\Lambda _{p}\left( \bar{\xi}%
_{jn}(x)\right) -\mathbf{E}\Lambda _{p}(\mathbb{Z}_{1})\right\vert
=o(h^{d/2}).
\end{equation*}%
This leads us to conclude that%
\begin{equation*}
\int_{A}\left\{ n^{p/2}h^{(p-1)d/2}\mathbf{E}\Lambda _{p}(v_{jn}(x))-h^{-d/2}%
\tilde{\rho}_{jn}^{p}(x)\mathbf{E}\Lambda _{p}(\mathbb{Z}_{1})\right\}
w_{j}(x)dx=o(1).
\end{equation*}%
Now, there exists $n_{0}$ such that for all $n>n_{0}$, sup$_{x\in \mathcal{S}%
_{j}}h^{d}b_{jn}^{2}(x)<c_{1}/2$, where $c_{1}>0$ is the constant in Lemma
A3. Observe that for all $n>n_{0},$%
\begin{eqnarray*}
&&\sup_{x\in \mathcal{S}_{j}}h^{-d/2}\left\vert \tilde{\rho}%
_{jn}^{p}(x)-\rho _{jn}^{p}(x)\right\vert \\
&=&\sup_{x\in \mathcal{S}_{j}}h^{-d/2}\left\vert (\rho
_{jn}^{2}(x)-h^{d}b_{jn}^{2}(x))^{p/2}-(\rho _{jn}^{2}(x))^{p/2}\right\vert
\\
&\leq &\sup_{x\in \mathcal{S}_{j}}\frac{p}{2}h^{d/2}b_{jn}^{2}(x)\cdot \max
\left\{ \left( \rho _{jn}^{2}(x)+c_{1}/2\right) ^{p/2-1},\left( \rho
_{jn}^{2}(x)-c_{1}/2\right) ^{p/2-1}\right\} \text{.}
\end{eqnarray*}%
By Lemma A4, the last term is $O(h^{d/2})=o(1)$. This completes the proof. $%
\blacksquare $\bigskip

Recall the definition: $\rho _{j}^{2}(x)\equiv \mathbf{E}%
[Y_{ji}^{2}|X_{i}=x]f(x)\int K^{2}(u)du$. Let%
\begin{align}  \label{def3}
\begin{split}
\sigma _{jk,n}(A) &\equiv n^{p}h^{(p-1)d}\int_{A}\int_{A}Cov\left( \Lambda
_{p}(v_{jN}(x)),\Lambda _{p}(v_{kN}(z))\right) w_{j}(x)w_{k}(z)dxdz, \text{
and} \\
\sigma _{jk}(A) &\equiv \int_{A}q_{jk,p}(x)\rho _{j}^{p}(x)\rho
_{k}^{p}(x)w_{j}(x)w_{k}(x)dx,
\end{split}%
\end{align}%
where we recall the definition:%
\begin{equation*}
q_{jk,p}(x)\equiv \int_{\lbrack -1,1]^{d}}Cov\left( \Lambda _{p}(\sqrt{%
1-t_{jk}^{2}(x,u)}\mathbb{Z}_{1}+t_{jk}(x,u)\mathbb{Z}_{2}),\Lambda _{p}(%
\mathbb{Z}_{2})\right) du.
\end{equation*}%
Now, let $(Z_{1n}(x),Z_{2n}(z))\in \mathbf{R}^{2}$ be a jointly normal
centered random vector whose covariance matrix is the same as that of $(\xi
_{jn}(x),\xi _{kn}(z))$ for all $x,z\in \mathbf{R}^{d}$. We define%
\begin{equation*}
\tau _{jk,n}(A)\equiv \int_{A}\int_{[-1,1]^{d}}g_{jk,n}(x,u)\lambda
_{jk,n}(x,x+uh)dudx,
\end{equation*}%
where%
\begin{eqnarray*}
\lambda _{jk,n}(x,z) &\equiv &\rho _{jn}^{p}(x)\rho
_{kn}^{p}(z)w_{j}(x)w_{k}(z)1_{A}(x)1_{A}(z),\ \text{and} \\
g_{jk,n}(x,u) &\equiv &Cov\left( \Lambda _{p}(Z_{1n}(x)),\Lambda
_{p}(Z_{2n}(x+uh))\right) .
\end{eqnarray*}

The following result generalizes Lemma 6.5 of GMZ from a univariate $X$ to a
multivariate $X$. The truncation arguments in their proof on pages 752 and
753 do not apply in the case of multivariate $X$. The proof of the following
lemma employs a different approach for this part.\bigskip

\noindent \textsc{Lemma A7:}\textbf{\ }\textit{Suppose that Assumptions 1
and 2 hold and let }$h\rightarrow 0$ \textit{as} $n\rightarrow \infty $ 
\textit{satisfying }limsup$_{n\rightarrow \infty }n^{-r/2+1}h^{(1-r/2)d}<C$ 
\textit{for any} $r\in \lbrack 2,2p+2]$ \textit{for some} $C>0$.

\noindent (i) \textit{Suppose that }$A\subset \mathcal{S}_{j}\cap \mathcal{S}%
_{k}$ \textit{is any Borel set.\ Then }%
\begin{equation*}
\sigma _{jk,n}(A)=\tau _{jk,n}(A)+o(1).
\end{equation*}

\noindent (ii) \textit{Suppose further that }$A$\textit{\ has a finite
Lebesgue measure, }$\rho _{j}(\cdot )\rho _{k}(\cdot )$ \textit{and }$%
w_{j}(\cdot )w_{k}(\cdot )$\textit{\ are continuous and bounded on }$A$, 
\textit{and}%
\begin{equation}
\sup_{x\in A}\left\vert \rho _{l,n}(x)-\rho _{l}(x)\right\vert \rightarrow 0,%
\text{ \textit{as }}n\rightarrow \infty ,\text{ \textit{for} }l\in \{j,k\}.
\label{conv4}
\end{equation}%
\textit{Then, as\ }$n\rightarrow \infty ,$\textit{\ }$\tau _{jk,n}(A)=\sigma
_{jk}(A)+o(1),$ \textit{and hence from }(i),%
\begin{equation*}
\sigma _{jk,n}(A)\rightarrow \sigma _{jk}(A).
\end{equation*}

\noindent \textsc{Proof:} (i) By change of variables, we write $\sigma
_{jk,n}(A)=\tilde{\tau}_{jk,n}(A)$, where%
\begin{equation*}
\tilde{\tau}_{jk,n}(A)\equiv \int_{A}\int_{[-1,1]^{d}}Cov\left( \Lambda
_{p}(\xi _{jn}(x)),\Lambda _{p}(\xi _{kn}(x+uh))\right) \lambda
_{jk,n}(x,x+uh)dudx.
\end{equation*}%
Fix $\varepsilon _{1}\in (0,1]$ and let $c(\varepsilon _{1})=(1+\varepsilon
_{1})^{2}-1$. Let $\eta _{1}$ and $\eta _{2}$ be two independent random
variables that are independent of $(\{Y_{ji},X_{i}:j\in \mathcal{J}%
\}_{i=1}^{\infty },N)$, each having a two-point distribution that gives two
points, $\{\sqrt{c(\varepsilon _{1})}\}$ and $\{-\sqrt{c(\varepsilon _{1})}%
\},$ the equal mass of 1/2, so that $\mathbf{E}\eta _{1}=\mathbf{E}\eta
_{2}=0$ and $Var(\eta _{1})=Var(\eta _{2})=c(\varepsilon _{1})$.
Furthermore, observe that for any $r\geq 1,$%
\begin{equation}
\mathbf{E}|\eta _{1}|^{r}=\frac{1}{2}|c(\varepsilon _{1})|^{r/2}+\frac{1}{2}%
|c(\varepsilon _{1})|^{r/2}\leq C\varepsilon _{1}^{r/2},  \label{bd8}
\end{equation}%
for some constant $C>0$ that depends only on $r$. Define%
\begin{equation*}
\xi _{jn,1}^{\eta }(x)\equiv \frac{\xi _{jn}(x)+\eta _{1}}{1+\varepsilon _{1}%
}\text{ and }\xi _{kn,2}^{\eta }(x+uh)\equiv \frac{\xi _{kn}(x+uh)+\eta _{2}%
}{1+\varepsilon _{1}}.
\end{equation*}%
Note that $Var(\xi _{jn,1}^{\eta }(x))=Var(\xi _{kn,2}^{\eta }(x+uh))=1$.
Let $(Z_{1n}^{\eta }(x),Z_{2n}^{\eta }(x+uh))$ be a jointly normal centered
random vector whose covariance matrix is the same as that of $(\xi
_{jn,1}^{\eta }(x),\xi _{kn,2}^{\eta }(x+uh))$ for all $(x,u)\in \mathbf{R}%
^{d}\times \lbrack -1,1]^{d}$. Define%
\begin{eqnarray*}
\tilde{\tau}_{jk,n}^{\eta }(A) &\equiv &\int_{A}\int_{[-1,1]^{d}}Cov\left(
\Lambda _{p}(\xi _{jn,1}^{\eta }(x)),\Lambda _{p}(\xi _{kn,2}^{\eta
}(x+uh))\right) \lambda _{jk,n}(x,x+uh)dudx, \\
\tau _{jk,n}^{\eta }(A) &\equiv &\int_{A}\int_{[-1,1]^{d}}Cov\left( \Lambda
_{p}(Z_{1n}^{\eta }(x)),\Lambda _{p}(Z_{2n}^{\eta }(x+uh))\right) \lambda
_{jk,n}(x,x+uh)dudx.
\end{eqnarray*}%
Then first observe that%
\begin{eqnarray*}
\left\vert \tilde{\tau}_{jk,n}(A)-\tilde{\tau}_{jk,n}^{\eta }(A)\right\vert
&\leq &\int_{A}\int_{[-1,1]^{d}}|\Delta _{jk,n,1}^{\eta }(x,u)|\lambda
_{jk,n}(x,x+uh)dudx \\
&&+\int_{A}\int_{[-1,1]^{d}}|\Delta _{jk,n,2}^{\eta }(x,u)|\lambda
_{jk,n}(x,x+uh)dudx,
\end{eqnarray*}%
where%
\begin{eqnarray*}
\Delta _{jk,n,1}^{\eta }(x,u) &\equiv &\mathbf{E}\Lambda _{p}(\xi _{jn}(x))%
\mathbf{E}\Lambda _{p}(\xi _{kn}(x+uh)) \\
&&-\mathbf{E}\Lambda _{p}(\xi _{jn,1}^{\eta }(x))\mathbf{E}\Lambda _{p}(\xi
_{kn,2}^{\eta }(x+uh))\text{ and} \\
\Delta _{jk,n,2}^{\eta }(x,u) &\equiv &\mathbf{E}\Lambda _{p}(\xi
_{jn}(x))\Lambda _{p}(\xi _{kn}(x+uh)) \\
&&-\mathbf{E}\Lambda _{p}(\xi _{jn,1}^{\eta }(x))\Lambda _{p}(\xi
_{kn,2}^{\eta }(x+uh)).
\end{eqnarray*}%
Since for any $a,b\in \mathbf{R}$, $|a_{+}^{p}-b_{+}^{p}|\leq p|a-b|\left(
|a|^{p-1}+|b|^{p-1}\right) $, we bound $|\Delta _{jk,n,2}^{\eta }(x,u)|$ by%
\begin{eqnarray*}
&&p\mathbf{E}\left[ |\xi _{jn}(x)-\xi _{jn,1}^{\eta }(x)|\left( |\xi
_{jn}(x)|^{p-1}+|\xi _{jn,1}^{\eta }(x)|^{p-1}\right) |\xi _{kn}(x+uh)|^{p}%
\right] \\
&&+p\mathbf{E}\left[ |\xi _{kn}(x+uh)-\xi _{kn,2}^{\eta }(x+uh)|\left( |\xi
_{kn}(x+uh)|^{p-1}+|\xi _{kn,2}^{\eta }(x+uh)|^{p-1}\right) |\xi
_{jn,1}^{\eta }(x)|^{p}\right] \\
&\equiv &A_{1n}(x,u)+A_{2n}(x,u),\text{ say.}
\end{eqnarray*}%
As for $A_{1n}(x,u)$,%
\begin{eqnarray*}
A_{1n}(x,u) &\leq &p\left( \mathbf{E}\left[ |\xi _{jn}(x)-\xi _{jn,1}^{\eta
}(x)|^{2}(|\xi _{jn}(x)|^{p-1}+|\xi _{jn,1}^{\eta }(x)|^{p-1})^{2}\right]
\right) ^{1/2} \\
&&\times \left( \mathbf{E}\left[ |\xi _{kn}(x+uh)|^{2p}\right] \right)
^{1/2}.
\end{eqnarray*}%
Define%
\begin{equation*}
s\equiv \left\{ 
\begin{array}{c}
(p+1)/(p-1)\text{ if }p>1 \\ 
2\text{ if }p=1,%
\end{array}%
\right.
\end{equation*}%
and $q\equiv (1-1/s)^{-1}$. Note that by H\"{o}lder inequality,%
\begin{eqnarray*}
&&\mathbf{E}\left[ |\xi _{jn}(x)-\xi _{jn,1}^{\eta }(x)|^{2}(|\xi
_{jn}(x)|^{p-1}+|\xi _{jn,1}^{\eta }(x)|^{p-1})^{2}\right] \\
&\leq &\left( \mathbf{E}\left[ |\xi _{jn}(x)-\xi _{jn,1}^{\eta }(x)|^{2q}%
\right] \right) ^{1/q}\left( \mathbf{E}\left[ (|\xi _{jn}(x)|^{p-1}+|\xi
_{jn,1}^{\eta }(x)|^{p-1})^{2s}\right] \right) ^{1/s}.
\end{eqnarray*}%
Now,%
\begin{eqnarray*}
&&\mathbf{E}\left[ |\xi _{jn}(x)-\xi _{jn,1}^{\eta }(x)|^{2q}\right]
=(1+\varepsilon _{1})^{-2q}\mathbf{E}\left[ |\varepsilon _{1}\xi
_{jn}(x)-\eta _{1}|^{2q}\right] \\
&\leq &2^{2q-1}(1+\varepsilon _{1})^{-2q}\left\{ \varepsilon _{1}^{2q}%
\mathbf{E}\left[ |\xi _{jn}(x)|^{2q}\right] +\mathbf{E}\left[ |\eta
_{1}|^{2q}\right] \right\} .
\end{eqnarray*}%
Applying Lemma A5 and (\ref{bd8}) to the last bound, we conclude that%
\begin{equation*}
\sup_{x\in \mathcal{S}_{j}}\mathbf{E}\left[ |\xi _{jn}(x)-\xi _{jn,1}^{\eta
}(x)|^{2q}\right] \leq \frac{C_{1}(\varepsilon _{1}^{2q}+\varepsilon
_{1}^{q})}{(1+\varepsilon _{1})^{2q}}\leq C_{2}\varepsilon _{1}^{q},
\end{equation*}%
for some constants $C_{1},C_{2}>0$. Using Lemma A5, we can also see that for
some constants $C_{3},C_{4}>0$,%
\begin{equation*}
\sup_{x\in \mathcal{S}_{j}}\mathbf{E}\left[ (|\xi _{jn}(x)|^{p-1}+|\xi
_{jn,1}^{\eta }(x)|^{p-1})^{2s}\right] \leq C_{3}
\end{equation*}%
and from some large $n$ on,%
\begin{equation*}
\sup_{u\in \lbrack -1,1]^{d}}\sup_{x\in \mathcal{S}_{k}}\mathbf{E}\left[
|\xi _{kn}(x+uh)|^{2p}\right] \leq \sup_{x\in \mathcal{S}_{k}^{\varepsilon
/2}}\mathbf{E}\left[ |\xi _{kn}(x)|^{2p}\right] \leq C_{4},
\end{equation*}%
for $\varepsilon >0$ in Assumption 1(iii).\ Therefore, for some constant $%
C>0,$%
\begin{equation*}
\sup_{u\in \lbrack -1,1]^{d}}\sup_{x\in \mathcal{S}_{j}\cap \mathcal{S}%
_{k}}A_{1n}(x,u)\leq C\sqrt{\varepsilon _{1}}.
\end{equation*}%
Using similar arguments for $A_{2n}(x,u)$, we deduce that for some constant $%
C>0,$%
\begin{equation}
\sup_{u\in \lbrack -1,1]^{d}}\sup_{x\in \mathcal{S}_{j}\cap \mathcal{S}%
_{k}}|\Delta _{jk,n,2}^{\eta }(x,u)|\leq C\sqrt{\varepsilon _{1}}.
\label{one}
\end{equation}

Let us turn to $\Delta _{jk,n,1}^{\eta }(x,u)$. We bound $|\Delta
_{jk,n,1}^{\eta }(x,u)|$ by%
\begin{eqnarray*}
&&p\mathbf{E}\left[ |\xi _{jn}(x)-\xi _{jn,1}^{\eta }(x)|\left( |\xi
_{jn}(x)|^{p-1}+|\xi _{jn,1}^{\eta }(x)|^{p-1}\right) \right] \mathbf{E}%
\left[ |\xi _{kn}(x+uh)|^{p}\right] \\
&&+p\mathbf{E}\left[ |\xi _{kn}(x+uh)-\xi _{kn,2}^{\eta }(x+uh)|\left( |\xi
_{kn}(x+uh)|^{p-1}+|\xi _{kn,1}^{\eta }(x+uh)|^{p-1}\right) \right] \mathbf{E%
}[|\xi _{jn,1}^{\eta }(x)|^{p}].
\end{eqnarray*}%
Using similar arguments for $\Delta _{jk,n,2}^{\eta }(x,u)$, we find that
for some constant $C>0,$%
\begin{equation}
\sup_{u\in \lbrack -1,1]^{d}}\sup_{x\in \mathcal{S}_{j}\cap \mathcal{S}%
_{k}}|\Delta _{jk,n,1}^{\eta }(x,u)|\leq C\sqrt{\varepsilon _{1}}.
\label{two}
\end{equation}%
By Lemma A4 and Assumption 1(ii), there exist $n_{0}>0$ and $C_{1},C_{2}>0$
such that for all $n\geq n_{0},$%
\begin{eqnarray}
&&\int_{A}\int_{[-1,1]^{d}}\lambda _{jk,n}(x,x+uh)dudx  \label{lambdd} \\
&\leq &C_{1}\int_{A}\int_{[-1,1]^{d}}w_{j}(x)w_{k}(x+uh)dudx  \notag \\
&\leq &C_{2}\sqrt{\int_{A}w_{j}^{2}(x)dx}\sqrt{\int_{A}%
\int_{[-1,1]^{d}}w_{k}^{2}(x+uh)dudx}<\infty .  \notag
\end{eqnarray}%
Hence 
\begin{equation*}
\left\vert \tilde{\tau}_{jk,n}(A)-\tilde{\tau}_{jk,n}^{\eta }(A)\right\vert
\leq C_{5}\sqrt{\varepsilon _{1}}\int_{A}\int_{[-1,1]^{d}}\lambda
_{jk,n}(x,x+uh)dudx\leq C_{6}\sqrt{\varepsilon _{1}},
\end{equation*}%
for some constants $C_{5}>0$ and $C_{6}>0$.

Since the choice of $\varepsilon _{1}>0$ was arbitrary, it remains for the
proof of Lemma A7(i) to prove that%
\begin{equation}
\left\vert \tilde{\tau}_{jk,n}^{\eta }(A)-\tau _{jk,n}(A)\right\vert =o(1),
\label{tcv}
\end{equation}%
as $n\rightarrow \infty $ and then $\varepsilon _{1}\rightarrow 0$. For any $%
x\in \mathcal{S}_{j}\cap \mathcal{S}_{k}$,%
\begin{equation*}
\left( \xi _{jn,1}^{\eta }(x),\xi _{kn,2}^{\eta }(x+uh)\right) ^{\prime }%
\overset{d}{=}\frac{1}{\sqrt{n}}\sum_{i=1}^{n}W_{n}^{(i)}(x,u),
\end{equation*}%
where $W_{n}^{(i)}(x,u)$'s are i.i.d. copies of $W_{n}(x,u)\equiv
(q_{jn}(x),q_{kn}(x+uh))^{\prime }$ with 
\begin{equation*}
q_{jn}(x)\equiv \left\{ \frac{\sum_{i\leq N_{1}}Y_{ji}K\left(
(x-X_{i})/h\right) -\mathbf{E}\left[ Y_{ji}K\left( (x-X_{i})/h\right) \right]
}{h^{d/2}\rho _{jn}(x)}+\eta _{1}\right\} /(1+\varepsilon _{1})\text{.}
\end{equation*}%
Using the same arguments as in the proof of Lemma A5, we find that for $j\in 
\mathcal{J},$%
\begin{equation}
\sup_{x\in \mathcal{S}_{j}^{\varepsilon }}\mathbf{E}\left[ |q_{jn}(x)|^{3}%
\right] \leq Ch^{-d/2},\text{ for some }C>0\text{.}  \label{bd0}
\end{equation}%
Let $\Sigma _{1n}$ be the $2\times 2\ $covariance matrix of $(\xi
_{jn,1}^{\eta }(x),\xi _{kn,2}^{\eta }(x+uh))^{\prime }$. Define 
\begin{equation*}
\tilde{\Lambda}_{n,p}(v)\equiv \Lambda _{p}([\Sigma
_{1n}^{1/2}v]_{1})\Lambda _{p}([\Sigma _{1n}^{1/2}v]_{2})\text{,\ }v\in 
\mathbf{R}^{2}\text{,}
\end{equation*}%
where $[a]_{j}$ of a vector $a\in \mathbf{R}^{2}$ indicates its $j$-th
entry. There exists some $C>0$ such that for all $n$, 
\begin{eqnarray}
\sup_{v\in \mathbf{R}^{2}}\frac{\left\vert \tilde{\Lambda}_{n,p}(v)-\tilde{%
\Lambda}_{n,p}(0)\right\vert }{1+||v||^{2p+2}\min \{||v||,1\}} &\leq &C\text{
and}  \label{bds} \\
\int \sup_{u\in \mathbf{R}^{2}:||z-u||\leq \delta }\left\vert \tilde{\Lambda}%
_{n,p}(z)-\tilde{\Lambda}_{n,p}(u)\right\vert d\Phi (z) &\leq &C\delta \text{
for all }\delta \in (0,1].  \notag
\end{eqnarray}%
The correlation between $\xi _{jn,1}^{\eta }(x)$ and $\xi _{kn,2}^{\eta
}(x+uh)$ is equal to%
\begin{equation*}
\mathbf{E}\left[ \xi _{jn,1}^{\eta }(x)\xi _{kn,2}^{\eta }(x+uh)\right] =%
\frac{\mathbf{E}[\xi _{jn}(x)\xi _{kn}(x+uh)]}{(1+\varepsilon _{1})^{2}}\in %
\left[ -(1+\varepsilon _{1})^{-2},(1+\varepsilon _{1})^{-2}\right] .
\end{equation*}%
Hence, as for $\tilde{W}_{n}^{(i)}(x,u)\equiv \Sigma
_{1n}^{-1/2}W_{n}^{(i)}(x,u),$ by (\ref{bd0}),%
\begin{eqnarray}
&&\sup_{x\in \mathcal{S}_{j}\cap \mathcal{S}_{k}}\mathbf{E}||\tilde{W}%
_{n}^{(i)}(x,u)||^{3}  \label{arg51} \\
&\leq &C_{1}(1-(\mathbf{E}[\xi _{jn,1}^{\eta }(x)\xi _{kn,2}^{\eta
}(x+uh)])^{2})^{-3/2}\left\{ \text{sup}_{x\in \mathcal{S}_{j}^{\varepsilon }}%
\mathbf{E}[|q_{jn}(x)|^{3}]+\text{sup}_{x\in \mathcal{S}_{k}^{\varepsilon }}%
\mathbf{E}[|q_{kn}(x)|^{3}]\right\}  \notag \\
&\leq &C_{1}(1-(1+\varepsilon _{1})^{-4})^{-3/2}\left\{ \text{sup}_{x\in 
\mathcal{S}_{j}^{\varepsilon }}\mathbf{E}[|q_{jn}(x)|^{3}]+\text{sup}_{x\in 
\mathcal{S}_{k}^{\varepsilon }}\mathbf{E}[|q_{kn}(x)|^{3}]\right\}  \notag \\
&\leq &C_{2}(1-(1+\varepsilon _{1})^{-4})^{-3/2}h^{-d/2},\text{ for some }%
C_{1},C_{2}>0,  \notag
\end{eqnarray}%
so that $n^{-1/2}\sup_{x\in \mathcal{S}_{j}\cap \mathcal{S}_{k}}\mathbf{E}||%
\tilde{W}_{n}^{(i)}(x,u)||^{3}=O(n^{-1/2}h^{-d/2}).$ By Lemma A2 and
following the arguments in (\ref{arg}) analogously,%
\begin{eqnarray*}
&&\sup_{x\in \mathcal{S}_{j}\cap \mathcal{S}_{k}}\left\vert \mathbf{E}\tilde{%
\Lambda}_{n,p}\left( \frac{1}{\sqrt{n}}\sum_{i=1}^{n}\tilde{W}%
_{n}^{(i)}(x,u)\right) -\mathbf{E}\tilde{\Lambda}_{n,p}\left( \tilde{Z}%
_{n}^{\eta }(x,u)\right) \right\vert \\
&=&O\left( n^{-1/2}h^{-d/2}\right) =o(1),
\end{eqnarray*}%
where $\tilde{Z}_{n}^{\eta }(x,u)\equiv \Sigma _{1n}^{-1/2}(Z_{1n}^{\eta
}(x),Z_{2n}^{\eta }(x+uh))^{\prime }$. Certainly by (\ref{bd8}) and Lemma A5,%
\begin{eqnarray*}
&&Cov(\Lambda _{p}(Z_{1n}^{\eta }(x)),\Lambda _{p}(Z_{2n}^{\eta }(x+uh))) \\
&\leq &\sqrt{\mathbf{E}|Z_{1n}^{\eta }(x)|^{2p}}\sqrt{\mathbf{E}%
|Z_{2n}^{\eta }(x+uh)|^{2p}}<C,
\end{eqnarray*}%
for some $C>0$ that does not depend on $\varepsilon _{1}$. Using (\ref%
{lambdd}), we apply the dominated convergence theorem to obtain that%
\begin{equation}
\left\vert \tau _{jk,n}^{\eta }(A)-\tilde{\tau}_{jk,n}^{\eta }(A)\right\vert
=o(1)  \label{cv5}
\end{equation}%
as $n\rightarrow \infty $ for each $\varepsilon _{1}>0$.

Finally, note from (\ref{one}) and (\ref{two}) that, for all $x\in A$ and
all $u\in \lbrack -1,1]^{d},$%
\begin{eqnarray*}
&&Cov(\Lambda _{p}(Z_{1n}^{\eta }(x)),\Lambda _{p}(Z_{2n}^{\eta }(x+uh))) \\
&=&Cov(\Lambda _{p}(Z_{1n}(x)),\Lambda _{p}(Z_{2n}(x+uh)))+o(1),
\end{eqnarray*}%
where the $o(1)$ term is one that converges to zero as $n\rightarrow \infty $
and then $\varepsilon _{1}\rightarrow 0$. Therefore, by the dominated
convergence theorem,%
\begin{equation*}
\left\vert \tau _{jk,n}^{\eta }(A)-\tau _{jk,n}(A)\right\vert =o(1),
\end{equation*}%
as $n\rightarrow \infty $ and then $\varepsilon _{1}\rightarrow 0$. In view
of (\ref{cv5}), this completes the proof of (\ref{tcv}) and, as a
consequence, that of (i).\bigskip

\noindent (ii) Define $t_{jk,n}(x,u)\equiv \mathbf{E(}\xi _{jn}(x)\xi
_{kn}(x+uh)),$%
\begin{eqnarray*}
e_{jk,n}(x,u) &\equiv &\frac{1}{h^{d}}\mathbf{E}\left[ Y_{ji}Y_{ki}K\left( 
\frac{x-X_{i}}{h}\right) K\left( \frac{x-X_{i}}{h}+u\right) \right] \text{
and} \\
e_{jk}(x,u) &\equiv &\rho _{jk}(x)\frac{\int K(z)K(z+u)dz}{\int K^{2}(u)du}.
\end{eqnarray*}%
By Assumption 1(i), and Lemma A4, for almost every $x\in A$ and for each $%
u\in \lbrack -1,1]^{d}$,%
\begin{align}  \label{arg5}
\begin{split}
t_{jk,n}(x,u) &= \frac{1}{\rho _{jn}(x)\rho _{kn}(x+uh)}\frac{1}{h^{d}}%
\mathbf{E}\left[ Y_{ji}Y_{ki}K\left( \frac{x-X_{i}}{h}\right) K\left( \frac{%
x-X_{i}}{h}+u\right) \right] \\
&= \frac{e_{jk,n}(x,u)}{\rho _{jn}(x)\rho _{kn}(x+uh)}=\frac{e_{jk}(x,u)}{%
\rho _{j}(x)\rho _{k}(x+uh)}+o(1)=t_{jk}(x,u)+o(1),
\end{split}%
\end{align}%
where we recall that $t_{jk}(x,u)=e_{jk}(x,u)/(\rho _{j}(x)\rho _{k}(x))$ by
the definition of $t_{jk}(\cdot ,\cdot )$.

By (\ref{conv4}),%
\begin{equation*}
\tau _{jk,n}(A)=\int_{A}\int_{[-1,1]^{d}}g_{jk,n}(x,u)\lambda
_{jk}(x,x+uh)dudx+o(1),
\end{equation*}%
where $\lambda _{jk}(x,z)\equiv \rho _{j}^{p}(x)\rho
_{k}^{p}(z)w_{j}(x)w_{k}(z)1_{A}(x)1_{A}(z).$ By (\ref{arg5}), for almost
every $x\in A$ and for each $u\in \lbrack -1,1]^{d}$,%
\begin{equation*}
g_{jk,n}(x,u)\rightarrow g_{jk}(x,u),\text{ as }n\rightarrow \infty \text{,}
\end{equation*}%
where $g_{jk}(x,u)\equiv Cov(\Lambda _{p}(\sqrt{1-t_{jk}^{2}(x,u)}\mathbb{Z}%
_{1}+t_{jk}(x,u)\mathbb{Z}_{2}),\Lambda _{p}(\mathbb{Z}_{2})).$ Furthermore,
since $\rho _{j}(\cdot )\rho _{k}(\cdot )$ and $w_{j}(\cdot )w_{k}(\cdot )$
are continuous on $A$ and $A$ has a finite Lebesgue measure, we follow the
proof of Lemma 6.4 of GMZ to find that $g_{jk,n}(x,u)\lambda _{jk}(x,x+uh)$
converges in measure to $g_{jk}(x,u)\lambda _{jk}(x,x)$ on $A\times \lbrack
-1,1]^{d}$, as $n\rightarrow \infty $. Using the bounded convergence
theorem, we deduce the desired result. $\blacksquare $\bigskip

The following lemma is a generalization of Lemma 6.2 of GMZ from $p=1$ to $%
p\geq 1.$ The proof of GMZ does not carry over to this general case because
the majorization inequality of Pinelis (1994) used in GMZ does not apply
here. (Note that (4) in Pinelis (1994) does not apply when $p>1$.)\bigskip

\noindent \textsc{Lemma A8:}\textbf{\ }\textit{Suppose that Assumptions 1
and 2 hold. Furthermore, assume that as }$n\rightarrow \infty ,$ $%
h\rightarrow 0$, $n^{-1/2}h^{-d}\rightarrow 0.$\textit{\ Then there exists a
constant }$C>0$ \textit{such that for any Borel set }$A\subset \mathbf{R}%
^{d} $ \textit{and for all} $j\in \mathcal{J},$%
\begin{eqnarray*}
&&\text{limsup}_{n\rightarrow \infty }\mathbf{E}\left[ \left\vert
n^{p/2}h^{(p-1)d/2}\int_{A}\left\{ \Lambda _{p}(v_{jn}(x))-\mathbf{E}\left[
\Lambda _{p}(v_{jn}(x))\right] \right\} w_{j}(x)dx\right\vert \right] \\
&\leq &C\int_{A}w_{j}(x)dx+C\sqrt{\int_{A}w_{j}^{2}(x)dx}.
\end{eqnarray*}

\noindent \textsc{Proof} \textsc{:}\textbf{\ }It suffices to show that there
exists $C>0$ such that for any Borel set $A\subset \mathbf{R}^{d},$\bigskip

\noindent \textsc{Step 1:} $\mathbf{E}\left[ \left\vert
n^{p/2}h^{(p-1)d/2}\int_{A}\left( \Lambda _{p}(v_{jn}(x))-\Lambda
_{p}(v_{jN}(x))\right) w_{j}(x)dx\right\vert \right] \leq
C\int_{A}w_{j}(x)dx,$\bigskip

\noindent \textsc{Step 2:} $\mathbf{E}\left[ \left\vert
n^{p/2}h^{(p-1)d/2}\int_{A}\left( \Lambda _{p}(v_{jN}(x))-\mathbf{E}\left[
\Lambda _{p}(v_{jN}(x))\right] \right) w_{j}(x)dx\right\vert \right] \leq C%
\sqrt{\int_{A}w_{j}^{2}(x)dx},$ and\bigskip

\noindent \textsc{Step 3:} $n^{p/2}h^{(p-1)d/2}\left\vert \int_{A}\left( 
\mathbf{E}\Lambda _{p}(v_{jN}(x))-\mathbf{E}\left[ \Lambda _{p}(v_{jn}(x))%
\right] \right) w_{j}(x)dx\right\vert \rightarrow 0$ as $n\rightarrow \infty 
$.\bigskip

Indeed, by chaining Steps 1, 2 and 3, we obtain the desired result.\bigskip

\noindent \textsc{Proof of Step 1:} For simplicity, let%
\begin{eqnarray*}
u_{jn}^{2}(x) &\equiv &\mathbf{E}\left[ Y_{ji}^{2}K^{2}\left( \frac{x-X_{i}}{%
h}\right) \right] -\left( \mathbf{E}\left[ Y_{ji}K\left( \frac{x-X_{i}}{h}%
\right) \right] \right) ^{2}\text{ and} \\
\bar{V}_{n,ji}(x) &\equiv &\frac{1}{u_{jn}(x)}\left\{ Y_{ji}K\left( \frac{%
x-X_{i}}{h}\right) -\mathbf{E}\left[ Y_{ji}K\left( \frac{x-X_{i}}{h}\right) %
\right] \right\} .
\end{eqnarray*}%
We write, if $N=n,\ \sum_{i=N+1}^{n}=0,$ and if $N>n$, $\sum_{i=N+1}^{n}=-%
\sum_{i=n+1}^{N}$. Using this notation, write%
\begin{equation*}
v_{jn}(x)=\frac{1}{nh^{d}}\sum_{i=1}^{N}\bar{V}_{n,ji}(x)u_{jn}(x)+\frac{1}{%
nh^{d}}\sum_{i=N+1}^{n}\bar{V}_{n,ji}(x)u_{jn}(x).
\end{equation*}%
Now, observe that%
\begin{eqnarray*}
\frac{1}{\sqrt{nh^{d}}}\sum_{i=1}^{N}\bar{V}_{n,ji}(x)u_{jn}(x) &=&\frac{1}{%
\sqrt{nh^{d}}}\sum_{i=1}^{N}\left\{ Y_{ji}K\left( \frac{x-X_{i}}{h}\right) -%
\mathbf{E}\left[ Y_{ji}K\left( \frac{x-X_{i}}{h}\right) \right] \right\} \\
&=&\sqrt{nh^{d}}\left\{ \hat{g}_{jN}(x)-\frac{1}{h^{d}}\mathbf{E}\left[
Y_{ji}K\left( \frac{x-X_{i}}{h}\right) \right] \right\} \\
&&+\sqrt{nh^{d}}\left( \frac{n-N}{n}\right) \cdot \frac{1}{h^{d}}\cdot 
\mathbf{E}\left[ Y_{ji}K\left( \frac{x-X_{i}}{h}\right) \right] \\
&=&\sqrt{nh^{d}}v_{jN}(x)+\sqrt{nh^{d}}\left( \frac{n-N}{n}\right) \cdot 
\frac{1}{h^{d}}\cdot \mathbf{E}\left[ Y_{ji}K\left( \frac{x-X_{i}}{h}\right) %
\right] .
\end{eqnarray*}%
Letting 
\begin{eqnarray*}
\eta _{jn}(x) &\equiv &\sqrt{n}\left( \frac{n-N}{n}\right) \cdot \frac{1}{%
h^{d}}\cdot \mathbf{E}\left[ Y_{ji}K\left( \frac{x-X_{i}}{h}\right) \right] 
\text{ and} \\
s_{jn}(x) &\equiv &\frac{1}{\sqrt{nh^{d}}}\sum_{i=N+1}^{n}\bar{V}%
_{n,ji}(x)u_{jn}(x),
\end{eqnarray*}%
we can write%
\begin{equation}
\sqrt{nh^{d}}v_{jn}(x)=\sqrt{nh^{d}}v_{jN}(x)+(\sqrt{h^{d}}\eta
_{jn}(x)+s_{jn}(x)).  \label{dec5}
\end{equation}%
First, note that for some constant $C>0$, 
\begin{equation}
\sup_{x\in \mathcal{S}_{j}}u_{jn}^{2}(x)\leq Ch^{d},  \label{bd5}
\end{equation}%
from some large $n$ on, by Lemma A4. Recall the definition of $\tilde{\rho}%
_{jn}(x):$ $\tilde{\rho}_{jn}(x)\equiv \sqrt{nh^{d}Var(v_{jn}(x))}$ and note
that%
\begin{equation*}
\tilde{\rho}_{jn}^{2}(x)=\rho
_{jn}^{2}(x)-h^{d}b_{jn}^{2}(x)=h^{-d}u_{jn}^{2}(x).
\end{equation*}%
As in the proof of Lemma A5, there exist $n_{0}$, $C_{1}>0$ and $C_{2}>0$
such that for all $n\geq n_{0},$%
\begin{eqnarray}
C_{1} &>&\sup_{x\in \mathcal{S}_{j}}\sqrt{\rho
_{jn}^{2}(x)-h^{d}b_{jn}^{2}(x)}=\sup_{x\in \mathcal{S}_{j}}\tilde{\rho}%
_{jn}(x)  \label{recall} \\
&\geq &\inf_{x\in \mathcal{S}_{j}}\tilde{\rho}_{jn}(x)\geq \inf_{x\in 
\mathcal{S}_{j}}\sqrt{\rho _{jn}^{2}(x)-h^{d}b_{jn}^{2}(x)}>C_{2}.  \notag
\end{eqnarray}%
Using (\ref{bd5}), (\ref{dec5}), and (\ref{recall}), we deduce that for some 
$C_{1},C_{2},C_{3},$ and $C_{4}>0,$%
\begin{eqnarray*}
\lefteqn{\left\vert n^{p/2}h^{(p-1)d/2}\int_{A}\left( \Lambda
_{p}(v_{jn}(x))-\Lambda _{p}(v_{jN}(x))\right) w_{j}(x)dx\right\vert } \\
&\leq &C_{1}h^{-d/2}\left\vert \int_{A}\left( \Lambda _{p}\left( \sqrt{nh^{d}%
}\frac{v_{jn}(x)}{\tilde{\rho}_{jn}(x)}\right) -\Lambda _{p}\left( \sqrt{%
nh^{d}}\frac{v_{jN}(x)}{\tilde{\rho}_{jn}(x)}\right) \right)
w_{j}(x)dx\right\vert \\
&\leq &C_{2}\int_{A}\left\vert \frac{\eta _{jn}(x)}{\tilde{\rho}_{jn}(x)}%
\right\vert \left( \left\vert \sqrt{nh^{d}}\frac{v_{jn}(x)}{\tilde{\rho}%
_{jn}(x)}\right\vert ^{p-1}+\left\vert \sqrt{nh^{d}}\frac{v_{jN}(x)}{\tilde{%
\rho}_{jn}(x)}\right\vert ^{p-1}\right) w_{j}(x)dx \\
&&+C_{3}h^{-d/2}\int_{A}\left\vert \frac{s_{jn}(x)}{\tilde{\rho}_{jn}(x)}%
\right\vert \left( \left\vert \sqrt{nh^{d}}\frac{v_{jn}(x)}{\tilde{\rho}%
_{jn}(x)}\right\vert ^{p-1}+\left\vert \sqrt{nh^{d}}\frac{v_{jN}(x)}{\tilde{%
\rho}_{jn}(x)}\right\vert ^{p-1}\right) w_{j}(x)dx \\
&\leq &C_{4}\int_{A}\left\vert \eta _{jn}(x)\right\vert \left( \left\vert 
\sqrt{nh^{d}}\frac{v_{jn}(x)}{\tilde{\rho}_{jn}(x)}\right\vert
^{p-1}+\left\vert \sqrt{nh^{d}}\frac{v_{jN}(x)}{\tilde{\rho}_{jn}(x)}%
\right\vert ^{p-1}\right) w_{j}(x)dx \\
&&+C_{3}\int_{A}\left\vert \frac{s_{jn}(x)}{u_{jn}(x)}\right\vert \left(
\left\vert \sqrt{nh^{d}}\frac{v_{jn}(x)}{\tilde{\rho}_{jn}(x)}\right\vert
^{p-1}+\left\vert \sqrt{nh^{d}}\frac{v_{jN}(x)}{\tilde{\rho}_{jn}(x)}%
\right\vert ^{p-1}\right) w_{j}(x)dx \\
&=&A_{1n}+A_{2n},\text{ say.}
\end{eqnarray*}%
To deal with $A_{1n}$ and $A_{2n}$, we first show the following:\bigskip

\noindent \textsc{Claim 1:} $\sup_{x\in \mathcal{S}_{j}}\mathbf{E}[\eta
_{jn}^{2}(x)]=O(1).$\bigskip

\noindent \textsc{Claim 2:} $\sup_{x\in \mathcal{S}_{j}}\mathbf{E}%
[\left\vert s_{jn}(x)/u_{jn}(x)\right\vert ^{2}]=o(1)$.\bigskip

\noindent \textsc{Claim 3:} $\sup_{x\in \mathcal{S}_{j}}\mathbf{E}[|\sqrt{%
nh^{d}}v_{jN}(x)/\tilde{\rho}_{jn}(x)|^{2p-2}]=O(1)$.\bigskip

\noindent \textsc{Proof of Claim 1:} By Lemma A4 and the fact that $\mathbf{E%
}|n^{-1/2}(n-N)|^{2}=O(1),$%
\begin{equation*}
\sup_{x\in \mathcal{S}_{j}}\mathbf{E}\left[ \eta _{jn}^{2}(x)\right] \leq 
\mathbf{E}\left\vert \sqrt{n}\left( \frac{n-N}{n}\right) \right\vert
^{2}\cdot \sup_{x\in \mathcal{S}_{j}}\left\vert \frac{1}{h^{d}}\cdot \mathbf{%
E}\left[ Y_{ji}K\left( \frac{x-X_{i}}{h}\right) \right] \right\vert
^{2}=O(1).
\end{equation*}

\noindent \textsc{Proof of Claim 2:} Note that%
\begin{equation}
\left\vert \sqrt{nh^{d}}\frac{s_{jn}(x)}{u_{jn}(x)}\right\vert =\left\vert
\sum_{i=n+1}^{N}\bar{V}_{n,ji}(x)\right\vert .  \label{eq5}
\end{equation}%
Certainly $Var(\bar{V}_{n,ji}(x))=1.$ As seen in (\ref{bd}), sup$_{x\in 
\mathcal{S}_{j}}\mathbf{E}\left\vert \bar{V}_{n,ji}(x)\right\vert ^{3}\leq
Ch^{-d/2}$ for some $C>0$. Similarly,%
\begin{equation*}
\sup_{x\in \mathcal{S}_{j}}\mathbf{E}\left\vert \bar{V}_{n,ji}(x)\right\vert
^{4}\leq \frac{h^{d}k_{jn,4}(x)}{h^{2d}\left( \rho
_{jn}^{2}(x)-h^{d}b_{jn}^{2}(x)\right) ^{2}}\leq Ch^{-d},
\end{equation*}%
for some $C>0$. Hence by Lemma 1(i) of Horv\'{a}th (1991), for some $C>0,$%
\begin{eqnarray*}
\mathbf{E}\left( \sqrt{nh^{d}}\frac{s_{jn}(x)}{u_{jn}(x)}\right) ^{2} &\leq &%
\mathbf{E}|N-n|\mathbf{E}|\mathbb{Z}_{1}|^{2} \\
&&+C\left\{ \mathbf{E}\left\vert N-n\right\vert ^{1/2}\mathbf{E}\left\vert 
\bar{V}_{n,ji}(x)\right\vert ^{3}+\mathbf{E}\left\vert \bar{V}%
_{n,ji}(x)\right\vert ^{4}\right\} .
\end{eqnarray*}
Note that $\mathbf{E}|N-n|=O(n^{1/2})$ and $\mathbf{E}\left\vert
N-n\right\vert ^{1/2}=O(n^{1/4})$ (e.g. (2.21) and (2.22) of Horv\'{a}th
(1991)). Therefore, there exists $C>0$ such that%
\begin{equation*}
\sup_{x\in \mathcal{S}_{j}}\mathbf{E}\left( \sqrt{nh^{d}}\frac{s_{jn}(x)}{%
u_{jn}(x)}\right) ^{2}\leq C\left\{ n^{1/2}+n^{1/4}h^{-d/2}+h^{-d}\right\} .
\end{equation*}%
Since $n^{-1/2}h^{-d}\rightarrow 0$, $\sup_{x\in \mathcal{S}_{j}}\mathbf{E}%
[(s_{jn}(x)/u_{jn}(x))^{2}]=o(1)$.\bigskip

\noindent \textsc{Proof of Claim 3:} By (\ref{bd6}), Lemmas A3-A4, and (\ref%
{recall}), we have%
\begin{equation*}
\sup_{x\in \mathcal{S}_{j}}\mathbf{E}\left[ \left\vert \sqrt{nh^{d}}\frac{%
v_{jN}(x)}{\tilde{\rho}_{jn}(x)}\right\vert ^{2p-2}\right] =\sup_{x\in 
\mathcal{S}_{j}}\left\vert \frac{\rho _{jn}(x)}{\tilde{\rho}_{jn}(x)}%
\right\vert ^{2p-2}\mathbf{E}\left( \left\vert \frac{\sqrt{nh^{d}}v_{jN}(x)}{%
\rho _{jn}(x)}\right\vert ^{2p-2}\right) \leq C,
\end{equation*}%
for some $C>0$. This completes the proof of Claim 3.\bigskip

Now, using Claims 1-3, we prove Step 1. Let $\mu _{j}(A)\equiv
\int_{A}w_{j}(x)dx$. Since $h^{(p-1)d/2}=O(1)$ when $p=1,$ and $\sqrt{a+b}%
\leq \sqrt{a}+\sqrt{b}$ for any $a\geq 0$ and $b\geq 0,$%
\begin{eqnarray*}
\mathbf{E}\left[ A_{1n}\right] &\leq &C\int_{A}\mathbf{E}\left[ \left\vert
\eta _{jn}(x)\right\vert \left( \left\vert \sqrt{nh^{d}}\frac{v_{jn}(x)}{%
\tilde{\rho}_{jn}(x)}\right\vert ^{p-1}+\left\vert \sqrt{nh^{d}}\frac{%
v_{jN}(x)}{\tilde{\rho}_{jn}(x)}\right\vert ^{p-1}\right) \right] w_{j}(x)dx
\\
&\leq &C\mu _{j}(A)\sup_{x\in \mathcal{S}_{j}}\mathbf{E}\left[ \left\vert
\eta _{jn}(x)\right\vert \left( \left\vert \sqrt{nh^{d}}\frac{v_{jn}(x)}{%
\tilde{\rho}_{jn}(x)}\right\vert ^{p-1}+\left\vert \sqrt{nh^{d}}\frac{%
v_{jN}(x)}{\tilde{\rho}_{jn}(x)}\right\vert ^{p-1}\right) \right] \\
&\leq &2C\mu _{j}(A)\times \left( \sup_{x\in \mathcal{S}_{j}}\mathbf{E}\left[
\eta _{jn}^{2}(x)\right] \right) ^{1/2} \\
&&\times \left( \left( \sup_{x\in \mathcal{S}_{j}}\mathbf{E}\left[
\left\vert \sqrt{nh^{d}}\frac{v_{jn}(x)}{\tilde{\rho}_{jn}(x)}\right\vert
^{2p-2}\right] \right) ^{1/2}+\left( \sup_{x\in \mathcal{S}_{j}}\mathbf{E}%
\left[ \left\vert \sqrt{nh^{d}}\frac{v_{jN}(x)}{\tilde{\rho}_{jn}(x)}%
\right\vert ^{2p-2}\right] \right) ^{1/2}\right) .
\end{eqnarray*}%
Certainly, as in the proof of Lemma A5,%
\begin{equation}
\sup_{x\in \mathcal{S}_{j}}\mathbf{E}\left[ \left\vert \sqrt{nh^{d}}\frac{%
v_{jn}(x)}{\tilde{\rho}_{jn}(x)}\right\vert ^{2p-2}\right] \leq C,
\label{bd4}
\end{equation}%
for some constant $C>0$. Hence using Claims 1 and 3, we conclude that $%
\mathbf{E}\left[ A_{1n}\right] \leq C\mu _{j}(A)$ for some $C>0$. As for $%
A_{2n}$, similarly, we obtain that for some $C>0,$%
\begin{eqnarray*}
\mathbf{E}\left[ A_{2n}\right] &\leq &C\int_{A}\mathbf{E}\left[ \left\vert 
\frac{s_{jn}(x)}{u_{jn}(x)}\right\vert \left( \left\vert \sqrt{nh^{d}}\frac{%
v_{jn}(x)}{\tilde{\rho}_{jn}(x)}\right\vert ^{p-1}+\left\vert \sqrt{nh^{d}}%
\frac{v_{jN}(x)}{\tilde{\rho}_{jn}(x)}\right\vert ^{p-1}\right) \right]
w_{j}(x)dx \\
&\leq &2C\mu _{j}(A)\times \left( \sup_{x\in \mathcal{S}_{j}}\mathbf{E}\left[
\left\vert \frac{s_{jn}(x)}{u_{jn}(x)}\right\vert ^{2}\right] \right) ^{1/2}
\\
&&\times \left( \left( \sup_{x\in \mathcal{S}_{j}}\mathbf{E}\left[
\left\vert \sqrt{nh^{d}}\frac{v_{jn}(x)}{\tilde{\rho}_{jn}(x)}\right\vert
^{2p-2}\right] \right) ^{1/2}+\left( \sup_{x\in \mathcal{S}_{j}}\mathbf{E}%
\left[ \left\vert \sqrt{nh^{d}}\frac{v_{jN}(x)}{\tilde{\rho}_{jn}(x)}%
\right\vert ^{2p-2}\right] \right) ^{1/2}\right) .
\end{eqnarray*}%
By Claims 2 and 3 and (\ref{bd4}), $\mathbf{E}\left[ A_{2n}\right] =o(1)$.
Hence the proof of Step 1 is completed.\bigskip

\noindent \textsc{Proof of Step 2:} We can follow the proof of Lemma A7(i)
to show that%
\begin{equation*}
\mathbf{E}\left[ n^{p/2}h^{(p-1)d/2}\int_{A}\left( \left\vert
v_{jN}(x)\right\vert ^{p}-\mathbf{E}\left[ \left\vert v_{jN}(x)\right\vert
^{p}\right] \right) w_{j}(x)dx\right] ^{2}=\kappa _{jn}(A)+o(1),
\end{equation*}%
where $\kappa _{jn}(A)\equiv \int_{A}\int_{[-1,1]^{d}}r_{jn}(x,u)\lambda
_{jn}(x,x+uh)dudx,$%
\begin{eqnarray*}
\lambda _{jn}(x,z) &\equiv &\rho _{jn}^{p}(x)\rho
_{jn}^{p}(z)w_{j}(x)w_{j}(z)1_{A\cap \mathcal{S}_{j}}(x)1_{A\cap \mathcal{S}%
_{j}}(z)\text{ and} \\
r_{jn}(x,u) &\equiv &Cov\left( |Z_{jn,A}(x)|^{p},|Z_{jn,B}(x+uh)|^{p}\right)
,
\end{eqnarray*}%
with $(Z_{jn,A}(x),Z_{jn,B}(x+uh))^{\prime }\in \mathbf{R}^{2}$ denoting a
centered normal random vector whose covariance matrix is equal to that of $%
(\xi _{jn}(x),\xi _{jn}(x+uh))^{\prime }$. By Cauchy-Schwarz inequality and
Lemma A5,%
\begin{equation*}
\sup_{x\in \mathcal{S}_{j}}r_{jn}(x,u)\leq \sup_{x\in \mathcal{S}_{j}}\sqrt{%
\mathbf{E}\left\vert Z_{jn,A}(x)\right\vert ^{2p}\mathbf{E}\left\vert
Z_{jn,B}(x+uh)\right\vert ^{2p}}<\infty .
\end{equation*}%
Furthermore, for each $u\in \lbrack -1,1]^{d},$%
\begin{equation*}
\int_{A}\lambda _{jn}(x,x+uh)dx\leq \sqrt{\int_{A}w_{j}^{2}(x)dx}\sqrt{%
\int_{A+uh}w_{j}^{2}(x)dx}.
\end{equation*}%
Since $\int_{\mathcal{S}_{j}^{\varepsilon }}w_{j}^{2}(x)dx<\infty $ for some 
$\varepsilon >0$ (Assumption 1(ii)), we find that as $h\rightarrow 0$, the
last term converges to $\int_{A}w_{j}^{2}(x)dx$. We obtain the desired
result of Step 2.$\bigskip $

\noindent \textsc{Proof of Step 3:} The convergence above follows from the
proof of Lemma A6. $\blacksquare $\bigskip

Let $\mathcal{C}\subset \mathbf{R}^{d}$ be a bounded Borel set such that%
\begin{equation*}
\alpha \equiv P\left\{ X\in \mathbf{R}^{d}\backslash \mathcal{C}\right\} >0.
\end{equation*}%
For any Borel set $A\subset \mathcal{C}$, let%
\begin{eqnarray*}
\zeta _{n}(A) &\equiv &\sum_{j=1}^{J}\int_{A}\Lambda
_{p}(v_{jn}(x))w_{j}(x)dx\text{ and\ } \\
\zeta _{N}(A) &\equiv &\sum_{j=1}^{J}\int_{A}\Lambda
_{p}(v_{jN}(x))w_{j}(x)dx.
\end{eqnarray*}%
We also let $\sigma _{n}^{2}(A)\equiv \sum_{j=1}^{J}\sum_{k=1}^{J}\sigma
_{jk,n}(A),$ and $\sigma ^{2}(A)\equiv \sum_{j=1}^{J}\sum_{k=1}^{J}\sigma
_{jk}(A)$. We define%
\begin{equation*}
S_{n}(A)\equiv \frac{n^{p/2}h^{(p-1)d/2}\{\zeta _{N}(A)-\mathbf{E}\zeta
_{N}(A)\}}{\sigma _{n}(A)},
\end{equation*}%
where%
\begin{eqnarray*}
U_{n} &\equiv &\frac{1}{\sqrt{n}}\left\{ \sum_{i=1}^{N}1\left\{ X_{i}\in 
\mathcal{C}\right\} -nP\left\{ X\in \mathcal{C}\right\} \right\} \text{, and}
\\
V_{n} &\equiv &\frac{1}{\sqrt{n}}\left\{ \sum_{i=1}^{N}1\{X_{i}\in \mathbf{R}%
^{d}\backslash \mathcal{C}\}-nP\{X\in \mathbf{R}^{d}\backslash \mathcal{C}%
\}\right\} .
\end{eqnarray*}

\noindent \textsc{Lemma A9:}\textbf{\ }\textit{Suppose that Assumptions 1
and 2 hold. Furthermore, assume that as} $n\rightarrow \infty ,$ $%
h\rightarrow 0$, \textit{and} $n^{-1/2}h^{-d}\rightarrow 0$. \textit{Let }$%
A\subset \mathcal{C}$ \textit{be such that }$\sigma ^{2}(A)>0$, $\alpha
\equiv P\{X\in \mathbf{R}^{d}\backslash \mathcal{C}\}>0,$ $\rho _{j}(\cdot )$%
\textit{'s and }$w_{j}(\cdot )$\textit{'s\ are continuous and bounded on }$%
A, $ \textit{and condition in }(\ref{conv4}) \textit{is satisfied for all }$%
l=1,\cdot \cdot \cdot ,J$. \textit{Then,}%
\begin{equation*}
(S_{n}(A),U_{n})\overset{d}{\rightarrow }(\mathbb{Z}_{1},\sqrt{1-\alpha }%
\mathbb{Z}_{2}).
\end{equation*}

\noindent \textsc{Proof} \textsc{:}\textbf{\ }First, we show that%
\begin{equation}
Cov\left( S_{n}(A),U_{n}\right) \rightarrow 0.  \label{SU}
\end{equation}%
Write%
\begin{equation*}
Cov\left( S_{n}(A),U_{n}\right) =\frac{n^{p/2}h^{(p-1)d/2}}{\sigma _{n}(A)}%
\sum_{j=1}^{J}\int_{A}Cov\left( \Lambda _{p}(v_{jN}(x)),U_{n}\right)
w_{j}(x)dx.
\end{equation*}%
It suffices for (\ref{SU}) to show that%
\begin{equation}
Cov\left( n^{p/2}h^{pd/2}\{\zeta _{N}(A)-\mathbf{E}\zeta
_{N}(A)\},U_{n}\right) =o(h^{d/2}),  \label{SU2}
\end{equation}%
since $\sigma _{n}^{2}(A)\rightarrow \sigma ^{2}(A)\equiv
\sum_{j=1}^{J}\sum_{k=1}^{J}\sigma _{jk}(A)>0$ by Lemma A7. For any $x\in 
\mathcal{S}_{j}$,%
\begin{equation*}
\left( \frac{\sqrt{nh^{d}}v_{jN}(x)}{\rho _{jn}(x)},\frac{U_{n}}{\sqrt{%
P\left\{ X\in \mathcal{C}\right\} }}\right) \overset{d}{=}\left( \frac{1}{%
\sqrt{n}}\sum_{k=1}^{n}Q_{n}^{(k)}(x),\frac{1}{\sqrt{n}}%
\sum_{k=1}^{n}U^{(k)}\right) ,
\end{equation*}%
where $(Q_{n}^{(k)}(x),U^{(k)})$'s are i.i.d. copies of $(Q_{n}(x),U)$ with%
\begin{eqnarray*}
Q_{n}(x) &\equiv &\frac{1}{h^{d/2}\rho _{jn}(x)}\left\{ \sum_{i\leq
N_{1}}Y_{ji}K\left( \frac{x-X_{i}}{h}\right) -\mathbf{E}\left[ Y_{ji}K\left( 
\frac{x-X_{i}}{h}\right) \right] \right\} \text{ and} \\
U &\equiv &\frac{\sum_{i\leq N_{1}}1\left\{ X_{i}\in \mathcal{C}\right\}
-P\left\{ X\in \mathcal{C}\right\} }{\sqrt{P\left\{ X\in \mathcal{C}\right\} 
}}.
\end{eqnarray*}%
Uniformly over $x\in \mathcal{S}_{j},$%
\begin{equation}
r_{n}(x)\equiv \mathbf{E}\left[ Q_{n}(x)U\right] =O(h^{d/2})=o(1),
\label{rate5}
\end{equation}%
by Lemma A4. Let $(Z_{1n},Z_{2n})^{\prime }$ be a centered normal random
vector with the same covariance matrix as that of $(Q_{n}(x),U)^{\prime }$.
Let the 2 by 2 covariance matrix be $\Sigma _{n,2}$.

Since $\frac{1}{\sqrt{n}}\sum_{k=1}^{n}U^{(k)}$ and $Z_{2n}$ have mean zero,
we write%
\begin{eqnarray*}
&&Cov\left( \Lambda _{p}\left( \frac{1}{\sqrt{n}}%
\sum_{k=1}^{n}Q_{n}^{(k)}(x)\right) ,\frac{1}{\sqrt{n}}\sum_{k=1}^{n}U^{(k)}%
\right) -Cov\left( \Lambda _{p}\left( Z_{1n}\right) ,Z_{2n}\right) \\
&=&\mathbf{E}\left[ \Lambda _{p}\left( \frac{1}{\sqrt{n}}%
\sum_{k=1}^{n}Q_{n}^{(k)}(x)\right) \left( \frac{1}{\sqrt{n}}%
\sum_{k=1}^{n}U^{(k)}\right) \right] -\mathbf{E}\left[ \Lambda _{p}\left(
Z_{1n}\right) Z_{2n}\right] \equiv A_{n}(x)\text{, say.}
\end{eqnarray*}%
Define $\bar{\Lambda}_{n,p}(v)\equiv \Lambda _{p}([\Sigma
_{n,2}^{1/2}v]_{1})[\Sigma _{n,2}^{1/2}v]_{2}$,\ $v\in \mathbf{R}^{2}$.
There exists some $C>0$ such that for all $n\geq 1,$%
\begin{eqnarray*}
\sup_{v\in \mathbf{R}^{2}}\frac{\left\vert \bar{\Lambda}_{n,p}(v)-\bar{%
\Lambda}_{n,p}(0)\right\vert }{1+||v||^{p+1}\min \{||v||,1\}} &\leq &C\text{
and} \\
\int \sup_{u\in \mathbf{R}^{2}:||z-u||\leq \delta }\left\vert \bar{\Lambda}%
_{n,p}(z)-\bar{\Lambda}_{n,p}(u)\right\vert d\Phi (z) &\leq &C\delta ,\text{
for all }\delta \in (0,1].
\end{eqnarray*}%
Letting $W_{n}^{(k)}(x)\equiv \Sigma _{n,2}^{-1/2}\cdot
(Q_{n}^{(k)}(x),U^{(k)})^{\prime },$ observe that using (\ref{rate5}) and
following the arguments in (\ref{arg51}), from some large $n$ on, for some $%
C>0,$%
\begin{eqnarray*}
\mathbf{E}||W_{n}^{(k)}(x)||^{3} &=&\mathbf{E}||\Sigma
_{n,2}^{-1/2}(Q_{n}^{(k)}(x),U^{(k)})^{\prime }||^{3} \\
&=&\mathbf{E}[\{tr(\Sigma _{n,2}^{-1/2}(Q_{n}^{(k)}(x),U^{(k)})^{\prime
}(Q_{n}^{(k)}(x),U^{(k)})\Sigma _{n,2}^{-1/2})\}^{3/2}] \\
&\leq &C(1-r_{n}^{2}(x))^{-3/2}\mathbf{E}\left[ |Q_{n}(x)|^{3}+|U|^{3}\right]
\leq Ch^{-d/2}.
\end{eqnarray*}%
Hence, by Lemma A2,%
\begin{eqnarray*}
&&\sup_{x\in \mathcal{S}_{j}}\left\vert A_{n}(x)\right\vert =\sup_{x\in 
\mathcal{S}_{j}}\left\vert \mathbf{E}\bar{\Lambda}_{n,p}\left( \frac{1}{%
\sqrt{n}}\sum_{i=1}^{n}W_{n}^{(i)}(x)\right) -\mathbf{E}\bar{\Lambda}%
_{n,p}\left( \tilde{Z}_{n}\right) \right\vert \\
&=&O\left( n^{-1/2}h^{-d/2}\right) =o(h^{d/2}),
\end{eqnarray*}%
where $\tilde{Z}_{n}\equiv \Sigma _{n,2}^{-1/2}(Z_{1n},Z_{2n})^{\prime }$.
This completes the proof of (\ref{SU2}) and hence that of (\ref{SU}).

Now, define%
\begin{equation*}
\Delta _{n}(x)\equiv n^{p/2}h^{(p-1)d/2}\sum_{j=1}^{J}\{\Lambda
_{p}(v_{jN}(x))-\mathbf{E}[\Lambda _{p}(v_{jN}(x))]\}w_{j}(x).
\end{equation*}%
Following Mason and Polonik (2009), we slice the integral $\int_{\mathcal{X}%
}\Delta _{n}(x)dx$ into a sum of a $1$-dependent random field. Let $\mathcal{%
C}$ be as given in the lemma. Let $\mathbb{Z}^{d}$ be the set of $d$-tupes
of integers, and let $\{R_{n,\mathbf{i}}:\mathbf{i}\in \mathbb{Z}^{d}\}$ be
the collection of rectangles in $\mathbf{R}^{d}$ such that $R_{n,\mathbf{i}%
}=[a_{n,\mathbf{i}_{1}},b_{n,\mathbf{i}_{1}}]\times \cdot \cdot \cdot \cdot
\times \lbrack a_{n,\mathbf{i}_{d}},b_{n,\mathbf{i}_{d}}]$, where $\mathbf{i}%
_{j}$ is the $j$-th entry of $\mathbf{i}$, and $h\leq b_{n,\mathbf{i}%
_{j}}-a_{n,\mathbf{i}_{j}}\leq 2h$, for all $j=1,\cdot \cdot \cdot ,d$, and
two different rectangles $R_{n,\mathbf{i}}$ and $R_{n,\mathbf{j}}$ do not
have intersection with nonempty interior, and the union of the rectangles $%
R_{n,\mathbf{i}}$, $\mathbf{i}\in \mathbb{Z}_{n}^{d}$, cover $\mathcal{C}$,
from some sufficiently large $n$ on, where $\mathbb{Z}_{n}^{d}$ be the set
of $d$-tuples of integers whose absolute values less than or equal to $n$.

We let $B_{n,\mathbf{i}}=R_{n,\mathbf{i}}\cap \mathcal{C}$ and $\mathcal{I}%
_{n}\equiv \{\mathbf{i}\in \mathbb{Z}_{n}^{d}:B_{n,\mathbf{i}}\neq
\varnothing \}$. Then $B_{n,\mathbf{i}}$ has Lebesgue measure $m(B_{n,%
\mathbf{i}})$ bounded by $C_{1}h^{d}$ and the cardinality of the set $%
\mathcal{I}_{n}$ is bounded by $C_{2}h^{-d}$ for some positive constants $%
C_{1}$ and $C_{2}$. Define%
\begin{eqnarray*}
\alpha _{\mathbf{i},n} &\equiv &\frac{1}{\sigma _{n}(A)}\int_{B_{n,\mathbf{i}%
}\cap A}\Delta _{n}(x)dx\text{ and} \\
u_{\mathbf{i},n} &\equiv &\frac{1}{\sqrt{n}}\left\{ \sum_{j=1}^{N}1\left\{
X_{j}\in B_{n,\mathbf{i}}\right\} -nP\left\{ X_{j}\in B_{n,\mathbf{i}%
}\right\} \right\} .
\end{eqnarray*}%
Then, we can write%
\begin{equation*}
S_{n}(A)=\sum_{\mathbf{i}\in \mathcal{I}_{n}}\alpha _{\mathbf{i},n}\text{
and }U_{n}=\sum_{\mathbf{i}\in \mathcal{I}_{n}}u_{\mathbf{i},n}.
\end{equation*}%
Certainly $Var(S_{n}(A))=1$ and it is easy to check that $%
Var(U_{n})=1-\alpha $. Take $\mu _{1},\mu _{2}\in \mathbf{R}$ and let%
\begin{equation*}
y_{\mathbf{i},n}\equiv \mu _{1}\alpha _{\mathbf{i},n}+\mu _{2}u_{\mathbf{i}%
,n}.
\end{equation*}%
From (\ref{SU}),%
\begin{equation*}
Var\left( \sum_{\mathbf{i}\in \mathcal{I}_{n}}y_{\mathbf{i},n}\right)
\rightarrow \mu _{1}^{2}+\mu _{2}^{2}(1-\alpha )\text{ as }n\rightarrow
\infty .
\end{equation*}%
Since $\sigma _{n}^{r}(A)=\sigma ^{r}(A)+o(1),\ r>0$, by Lemma A7 and $%
m(B_{n,\mathbf{i}})\leq Ch^{d}$ for a constant $C>0$, we take $r\in
(2,(2p+2)/p]$ and bound%
\begin{eqnarray*}
&&\sigma _{n}^{r}(A)\sum_{\mathbf{i}\in \mathcal{I}_{n}}\mathbf{E}|\alpha _{%
\mathbf{i},n}|^{r} \\
&\leq &C\sup_{x\in A}\mathbf{E}\left\vert \Delta _{n}(x)\right\vert
^{r}\sum_{\mathbf{i}\in \mathcal{I}_{n}}\left(
\int_{A}\int_{A}\int_{A}1_{B_{n,\mathbf{i}}}(u,v,s)dudvds\right) ^{r/3},
\end{eqnarray*}%
where $1_{B}(u,v,s)\equiv 1\{u\in B\}1\{v\in B\}1\{s\in B\}$. Using Jensen's
inequality, we have%
\begin{eqnarray*}
\sup_{x\in A}\mathbf{E}\left\vert \Delta _{n}(x)\right\vert ^{r} &\leq
&C_{1}n^{rp/2}h^{r(p-1)d/2}\sup_{x\in A}\sum_{j=1}^{J}\mathbf{E}\left\vert
v_{jN}(x)\right\vert ^{rp}w_{j}^{r}(x) \\
&\leq &C_{2}n^{rp/2}h^{r(p-1)d/2}\max_{1\leq j\leq J}\sup_{x\in A\cap 
\mathcal{S}_{j}}\mathbf{E}\left\vert v_{jN}(x)\right\vert ^{rp}
\end{eqnarray*}%
for some $C_{1},C_{2}>0$. As for the last term, we apply Rosenthal's
inequality (see. e.g. Lemma 2.3. of GMZ): for some constant $C>0,$%
\begin{eqnarray*}
&&n^{rp/2}h^{r(p-1)d/2}\sup_{x\in A\cap \mathcal{S}_{j}}\mathbf{E}\left\vert
v_{jN}(x)\right\vert ^{rp} \\
&\leq &Ch^{r(p-1)d/2}\sup_{x\in A\cap \mathcal{S}_{j}}\left( \frac{1}{h^{2d}}%
\mathbf{E}\left[ Y_{ji}^{2}K^{2}\left( \frac{x-X_{i}}{h}\right) \right]
\right) ^{rp/2} \\
&&+Ch^{r(p-1)d/2}\sup_{x\in A\cap \mathcal{S}_{j}}\left( \frac{n}{%
n^{rp/2}h^{rpd}}\mathbf{E}\left[ \left\vert Y_{ji}K\left( \frac{x-X_{i}}{h}%
\right) \right\vert ^{rp}\right] \right) .
\end{eqnarray*}%
By Lemma A4, the first term is $O(h^{-rd/2})$ and the last term is $%
O(n^{1-rp/2}h^{-rdp/2-rd/2+d}).$ Hence we find that 
\begin{eqnarray*}
\sum_{\mathbf{i}\in \mathcal{I}_{n}}\mathbf{E}|\alpha _{\mathbf{i},n}|^{r}
&=&\text{Cardinality of }\mathcal{I}_{n}\times O\left( m(B_{n,\mathbf{i}%
})^{r}h^{-rd/2}\{1+n^{1-rp/2}h^{-rdp/2+d}\}\right)  \\
&=&O\left( h^{rd/2-d}\{1+n^{1-rp/2}h^{-rdp/2+d}\}\right) =o(1)
\end{eqnarray*}%
for any $r\in (2,(2p+2)/p],$ because $n^{-1/2}h^{-d}\rightarrow 0$.
Therefore, as $n\rightarrow \infty ,$%
\begin{equation*}
\sum_{\mathbf{i}\in \mathcal{I}_{n}}\mathbf{E}|\alpha _{\mathbf{i}%
,n}|^{r}\rightarrow 0\text{ for any }r\in (2,(2p+2)/p].
\end{equation*}

Also, arguing similarly as in (6.56) of GMZ, we can show that $\sum_{\mathbf{%
i}\in \mathcal{I}_{n}}\mathbf{E}|u_{\mathbf{i},n}|^{r}\rightarrow 0$ as $%
n\rightarrow \infty $ for any $r\in (2,(2p+2)/p]$. Since $X_{i}$'s are
common across different $j$'s, the sequence $\{y_{\mathbf{i},n}\}_{\mathbf{i}%
\in \mathcal{I}_{n}}$ is a $1$-dependent random field (see Mason and Polonik
(2009)). The desired result of Lemma A9 follows by Theorem 1 of Shergin
(1993) and the Cram\'{e}r-Wold device. $\blacksquare $\bigskip

\noindent \textsc{Lemma A10:}\textbf{\ }\textit{Suppose that the conditions
of Lemma A9\ are satisfied, and let }$A\subset \mathbf{R}^{d}$\textit{\ be a
Borel set in Lemma A9}.\textit{\ Then,}%
\begin{equation*}
\frac{n^{p/2}h^{(p-1)d/2}\left\{ \zeta _{n}(A)-\mathbf{E}\zeta
_{n}(A)\right\} }{\sigma _{n}(A)}\overset{d}{\rightarrow }N(0,1),\text{ 
\textit{as }}n\rightarrow \infty .
\end{equation*}

\noindent \textsc{Proof:}\textbf{\ }The conditional distribution of $%
S_{n}(A) $ given $N=n$ is equal to that of%
\begin{equation*}
\frac{n^{p/2}h^{(p-1)d/2}}{\sigma _{n}(A)}\sum_{j=1}^{J}\int_{A}\left\{
\Lambda _{p}(v_{jn}(x))-\mathbf{E}\Lambda _{p}(v_{jN}(x))\right\} w_{j}(x)dx.
\end{equation*}%
Using Lemma A9 and the de-Poissonization argument of Beirlant and Mason
(1995) (see also Lemma 2.4 of GMZ), this conditional distribution converges
to $N(0,1).$ Now by Lemma A6, it follows that%
\begin{equation*}
n^{p/2}h^{(p-1)d/2}\sum_{j=1}^{J}\int_{A}\left\{ \mathbf{E}\Lambda
_{p}(v_{jN}(x))-\mathbf{E}\Lambda _{p}(v_{jn}(x))\right\}
w_{j}(x)dx\rightarrow 0,
\end{equation*}%
as $n\rightarrow \infty $. This completes the proof. $\blacksquare $\bigskip

\noindent \textsc{Proof of Theorem 1} \textsc{:} Fix $\varepsilon >0$ as in
Assumption 1(iii), and take $n_{0}>0$ such that for all $n\geq n_{0},$%
\begin{equation*}
\{x-uh:x\in \mathcal{S}_{j},u\in \lbrack -1/2,1/2]^{d}\}\subset \mathcal{S}%
_{j}^{\varepsilon }\subset \mathcal{X}\text{ for all }j\in \mathcal{J}.
\end{equation*}%
Since we are considering the least favorable case of the null hypothesis,%
\begin{equation*}
\mathbf{E}[Y_{ji}K((x-X_{i})/h)]/h^{d}=%
\int_{[-1/2,1/2]^{d}}m_{j}(x-uh)K(u)du=0,\text{ for almost all }x\in 
\mathcal{S}_{j}\text{,}
\end{equation*}%
for all $n\geq n_{0}$ and for all $j\in \mathcal{J}$. Therefore, $\hat{g}%
_{jn}(x)=v_{jn}(x)$ for almost all $x\in \mathcal{S}_{j}$, $j\in \mathcal{J}%
, $ and for all $n\geq n_{0}$. From here on, we consider only $n\geq n_{0}$.

We fix $0<\varepsilon _{l}\rightarrow 0$ as $l\rightarrow \infty $ and take
a compact set $\mathcal{W}_{l}\subset \mathcal{S}_{j}$ such that for each $%
j\in \mathcal{J}$, $w_{j}$ is bounded and continuous on $\mathcal{W}_{l}$
and for $s\in \{1,2\},$%
\begin{equation}
\int_{\mathcal{X}\backslash \mathcal{W}_{l}}w_{j}^{s}(x)dx\rightarrow 0\text{
as }l\rightarrow \infty \text{.}  \label{cv7}
\end{equation}%
We can choose such $\mathcal{W}_{l}$ following the arguments in the proof of
Lemma 6.1 of GMZ because $w_{j}^{s}$ is integrable by Assumption 1(ii). Take 
$M_{l,j},v_{l,j}>0,$ $j=1,2,\cdot \cdot \cdot ,J,\ $such that$\ $for$\ 
\mathcal{C}_{l,j}\equiv \lbrack -M_{l,j}+v_{l,j},M_{l,j}-v_{l,j}]^{d},$%
\begin{equation*}
P\left\{ X_{i}\in \mathbf{R}^{d}\backslash \mathcal{C}_{l,j}\right\} >0,
\end{equation*}%
and for some Borel $A_{l,j}\subset \mathcal{C}_{l,j}\cap \mathcal{W}_{l},$ $%
\rho _{j}(\cdot )$ is bounded and continuous on $A_{l,j}$,%
\begin{eqnarray}
\sup_{x\in A_{l,j}}\left\vert \rho _{jn}(x)-\rho _{j}(x)\right\vert
&\rightarrow &0,\text{ as }n\rightarrow \infty ,\text{ and}  \label{cvs} \\
\int_{\mathcal{W}_{l}\backslash A_{l,j}}\rho _{j}(x)w_{j}^{s}(x)dx
&\rightarrow &0,\text{ as }l\rightarrow \infty ,\text{ for }s\in \{1,2\}%
\text{.}  \notag
\end{eqnarray}%
The existence of $M_{l,j},v_{l,j}$ and $\varepsilon _{l}$ and the sets $%
A_{l,j}$ are ensured by Lemma A1. By Assumption 1(i), we find that the
second convergence in (\ref{cvs}) implies that $\int_{\mathcal{W}%
_{l}\backslash A_{l,j}}w_{j}^{s}(x)dx\rightarrow 0$ as $l\rightarrow \infty $%
, for $s\in \{1,2\}$. Now, take\ $A_{l}\equiv \cap _{j=1}^{J}A_{l,j}$ and $%
\mathcal{C}_{l}\equiv \cap _{j=1}^{J}\mathcal{C}_{l,j}$, and observe that
for $s\in \{1,2\}$,%
\begin{equation}
\int_{\mathcal{W}_{l}\backslash A_{l}}w_{j}^{s}(x)dx\leq \sum_{j=1}^{J}\int_{%
\mathcal{W}_{l}\backslash A_{l,j}}w_{j}^{s}(x)dx\rightarrow 0,  \label{cv4}
\end{equation}%
as $l\rightarrow \infty $ for all $j\in \mathcal{J}$.

First, we write%
\begin{eqnarray}
&&\frac{\sum_{j=1}^{J}\left\{ n^{p/2}h^{(p-1)d/2}\Gamma _{j}(\hat{g}%
_{jn})-a_{jn}\right\} }{\sigma _{n}}  \label{dec3} \\
&=&\frac{n^{p/2}h^{(p-1)d/2}}{\sigma _{n}}\left\{ \zeta _{n}(\mathcal{X}%
\backslash \mathcal{W}_{l})-\mathbf{E}\zeta _{n}(\mathcal{X}\backslash 
\mathcal{W}_{l})\right\}  \notag \\
&&+\frac{n^{p/2}h^{(p-1)d/2}}{\sigma _{n}}\left\{ \zeta _{n}(\mathcal{W}%
_{l}\backslash A_{l})-\mathbf{E}\zeta _{n}(\mathcal{W}_{l}\backslash
A_{l})\right\}  \notag \\
&&+\frac{n^{p/2}h^{(p-1)d/2}}{\sigma _{n}}\left\{ \zeta _{n}(A_{l})-\mathbf{E%
}\zeta _{n}(A_{l})\right\} .  \notag
\end{eqnarray}%
Since $\mathcal{X}\backslash A_{l}=(\mathcal{X}\backslash \mathcal{W}%
_{l})\cup (\mathcal{W}_{l}\backslash A_{l})$, by Lemma A8, (\ref{cv7}), and (%
\ref{cv4}),%
\begin{equation}
n^{p/2}h^{(p-1)d/2}\left\{ \zeta _{n}(\mathcal{X}\backslash A_{l})-\mathbf{E}%
\zeta _{n}(\mathcal{X}\backslash A_{l})\right\} \overset{p}{\rightarrow }0,\ 
\text{as }n\rightarrow \infty ,\text{ and }l\rightarrow \infty .  \label{cv3}
\end{equation}%
Furthermore, we write $\left\vert \sigma _{n}^{2}-\sigma
_{n}^{2}(A_{l})\right\vert $ as%
\begin{eqnarray*}
&&\sum_{j=1}^{J}\sum_{k=1}^{J}\int_{\mathcal{X}}q_{jk,p}(x)\left(
1-1_{A_{l}}(x)\right) \rho _{jn}^{p}(x)\rho _{kn}^{p}(x)w_{j}(x)w_{k}(x)dx \\
&\leq &\sum_{j=1}^{J}\sum_{k=1}^{J}\sup_{x\in \mathcal{S}_{j}\cap \mathcal{S}%
_{k}}\left\vert q_{jk,p}(x)\rho _{jn}^{p}(x)\rho _{kn}^{p}(x)\right\vert
\int_{\mathcal{X}}\left( 1-1_{A_{l}}(x)\right) w_{j}(x)w_{k}(x)dx \\
&=&\sum_{j=1}^{J}\sum_{k=1}^{J}\sup_{x\in \mathcal{S}_{j}\cap \mathcal{S}%
_{k}}\left\vert q_{jk,p}(x)\rho _{jn}^{p}(x)\rho _{kn}^{p}(x)\right\vert
\int_{\mathcal{X}\backslash A_{l}}w_{j}(x)w_{k}(x)dx.
\end{eqnarray*}%
Observe that as $l\rightarrow \infty ,$%
\begin{equation*}
\left\vert \int_{\mathcal{X}\backslash A_{l}}w_{j}(x)w_{k}(x)dx\right\vert
^{2}\leq \left( \int_{\mathcal{X}\backslash A_{l}}w_{j}^{2}(x)dx\right)
\left( \int_{\mathcal{X}\backslash A_{l}}w_{k}^{2}(x)dx\right) \rightarrow 0.
\end{equation*}%
From Lemma A4, it follows that%
\begin{equation}
\text{lim}_{l\rightarrow \infty }\text{limsup}_{n\rightarrow \infty
}\left\vert \sigma _{n}^{2}-\sigma _{n}^{2}(A_{l})\right\vert =0.
\label{lim}
\end{equation}%
Furthermore, since $\sigma _{n}^{2}(A_{l})\rightarrow \sigma ^{2}(A_{l})$ as 
$n\rightarrow \infty $ for each $l$ by Lemma A7, and $\sigma
^{2}(A_{l})\rightarrow \sigma ^{2}>0$ as $l\rightarrow \infty $, by
Assumption 1, it follows that for any $\varepsilon _{1}\in (0,\sigma ^{2})$,%
\begin{eqnarray}
0 &<&\sigma ^{2}-\varepsilon _{1}\leq \text{liminf}_{n\rightarrow \infty
}\sigma _{n}^{2}  \label{bdd2} \\
&\leq &\text{limsup}_{n\rightarrow \infty }\sigma _{n}^{2}\leq \sigma
^{2}+\varepsilon _{1}<\infty .  \notag
\end{eqnarray}%
Combining this with (\ref{cv3}), we find that as $n\rightarrow \infty $ and $%
l\rightarrow \infty ,$%
\begin{equation*}
\frac{n^{p/2}h^{(p-1)d/2}}{\sigma _{n}}\left\{ \zeta _{n}(\mathcal{X}%
\backslash A_{l})-\mathbf{E}\zeta _{n}(\mathcal{X}\backslash A_{l})\right\}
=o_{P}(1).
\end{equation*}

As for the last term in (\ref{dec3}), by (\ref{bdd2}) and Lemma A10, as $%
n\rightarrow \infty $ and $l\rightarrow \infty ,$%
\begin{equation*}
n^{p/2}h^{(p-1)d/2}\left\vert \zeta _{n}(A_{l})-\mathbf{E}\zeta
_{n}(A_{l})\right\vert =O_{P}(1).
\end{equation*}%
Therefore, by (\ref{lim}),%
\begin{eqnarray*}
&&\frac{n^{p/2}h^{(p-1)d/2}}{\sigma _{n}}\left\{ \zeta _{n}(A_{l})-\mathbf{E}%
\zeta _{n}(A_{l})\right\} \\
&=&\frac{n^{p/2}h^{(p-1)d/2}}{\sigma _{n}(A_{l})}\left\{ \zeta _{n}(A_{l})-%
\mathbf{E}\zeta _{n}(A_{l})\right\} +o_{P}(1),
\end{eqnarray*}%
where $o_{P}(1)$ is a term that vanishes in probability as $n\rightarrow
\infty $ and $l\rightarrow \infty $. For each $l\geq 1,$ the last term
converges in distribution to $N(0,1)$ by Lemma A10. Since $\sigma
_{n}^{2}(A_{l})\rightarrow \sigma ^{2}$ as $n\rightarrow \infty $ and $%
l\rightarrow \infty $, we conclude that%
\begin{equation*}
\sum_{j=1}^{J}\left\{ n^{p/2}h^{(p-1)d/2}\Gamma _{j}(\hat{g}%
_{jn})-a_{jn}\right\} \overset{d}{\rightarrow }N\left( 0,\sigma ^{2}\right) .
\end{equation*}%
$\blacksquare $

\subsection{Proofs of Other Theorems}

We now give proofs of other theorems in the paper.

\noindent \textsc{Proof of Theorem 2} \textsc{:} We first show that for each 
$j\in \mathcal{J},$%
\begin{eqnarray}
\hat{a}_{jn} &=&a_{jn}+O_{P}(n^{-1/2}h^{-3d/2})\text{ and}  \label{conv5} \\
\hat{\sigma}_{n}^{2} &=&\sigma _{n}^{2}+O_{P}(n^{-1/2}h^{-3d/2}).  \notag
\end{eqnarray}%
For this, we show that for all $j,k=1,\cdot \cdot \cdot ,J,$%
\begin{equation}
\sup_{x\in \mathcal{S}_{j}\cap \mathcal{S}_{k}}|\hat{\rho}_{jk,n}(x)-\rho
_{jk,n}(x)|=O_{P}\left( n^{-1/2}h^{-d}\right) .  \label{conv2}
\end{equation}%
Write $\sup_{x\in \mathcal{S}_{j}\cap \mathcal{S}_{k}}|\hat{\rho}%
_{jk,n}(x)-\rho _{jk,n}(x)|$ as%
\begin{equation*}
\sup_{x\in \mathcal{S}_{j}\cap \mathcal{S}_{k}}\left\vert \frac{1}{nh^{d}}%
\sum_{i=1}^{n}\left\{ Y_{ji}Y_{ki}K^{2}\left( \frac{x-X_{i}}{h}\right) -%
\mathbf{E}\left[ Y_{ji}Y_{ki}K^{2}\left( \frac{x-X_{i}}{h}\right) \right]
\right\} \right\vert .
\end{equation*}%
Let $\varphi _{n,x}(y_{1},y_{2},z)\equiv y_{1}y_{2}K^{2}((x-z)/h)$ and $%
\mathcal{K}_{n}\equiv \{\varphi _{n,x}(\cdot ,\cdot ,\cdot ):x\in \mathcal{S}%
_{j}\cap \mathcal{S}_{k}\}$. We define $N(\varepsilon ,\mathcal{K}%
_{n},L_{2}(Q))$ to be a covering number of $\mathcal{K}_{n}$ with respect to 
$L_{2}(Q)$, i.e., the smallest number of maps $\varphi _{j},\ j=1,\cdot
\cdot \cdot ,N_{1},$ such that for all $\varphi \in \mathcal{K}_{n},$ there
exists $\varphi _{j}$ such that $\int (\varphi _{j}-\varphi )^{2}dQ\leq
\varepsilon ^{2}$. By Assumption 2(b), Lemma 2.6.16 of van der Vaart and
Wellner (1996), and Lemma A.1 of Ghosal, Sen and van der Vaart (2000), we
find that for some $C>0$,%
\begin{equation*}
\sup_{Q}\log N(\varepsilon ,\mathcal{K}_{n},L_{2}(Q))\leq C\log \varepsilon ,
\end{equation*}%
where the supremum is over all discrete probability measures. We take $\bar{%
\varphi}_{n}(y_{1},y_{2},z)\equiv y_{1}y_{2}||K||_{\infty }^{2}$ to be the
envelope of $\mathcal{K}_{n}$. By Theorem 2.14.1 of van der Vaart and
Wellner (1996), we deduce that%
\begin{equation*}
n^{1/2}h^{d}\mathbf{E}\left[ \sup_{x\in \mathcal{S}_{j}\cap \mathcal{S}_{k}}|%
\hat{\rho}_{jk,n}(x)-\rho _{jk,n}(x)|\right] \leq C,
\end{equation*}%
for some positive constant $C$. This yields (\ref{conv2}). In view of the
definitions of $\hat{a}_{jn}$ and $\hat{\sigma}_{n}^{2},$ and Lemma A4, this
completes the proof of (\ref{conv5}).

Since $g_{j}(x)\leq 0$ for all $x\in \mathcal{X}$ under the null hypothesis
and $K$ is nonnegative,%
\begin{eqnarray*}
\sup_{x\in \mathcal{S}_{j}}\mathbf{E}\hat{g}_{jn}(x) &=&\sup_{x\in \mathcal{S%
}_{j}}\int g_{j}(x-uh)K\left( u\right) du\leq \int \sup_{x\in \mathcal{S}%
_{j}}g_{j}(x-uh)K\left( u\right) du \\
&\leq &\int \sup_{x\in \mathcal{X}}g_{j}(x)K\left( u\right) du=\sup_{x\in 
\mathcal{X}}g_{j}(x)\leq 0,
\end{eqnarray*}%
from some large $n$ on. The second inequality follows from Assumption
1(iii). Therefore,%
\begin{equation*}
\int_{\mathcal{X}}\Lambda _{p}(\hat{g}_{jn}(x))w_{j}(x)dx\leq \int_{\mathcal{%
X}}\Lambda _{p}(\hat{g}_{jn}(x)-\mathbf{E}\hat{g}_{jn}(x))w_{j}(x)dx.
\end{equation*}%
Hence by using this and (\ref{conv5}), we bound $P\{\hat{T}_{n}>z_{1-\alpha
}\}$ by%
\begin{equation*}
P\left\{ \frac{1}{\sigma _{n}}\sum_{j=1}^{J}\left\{ n^{p/2}h^{(p-1)d/2}\int_{%
\mathcal{X}}\Lambda _{p}(\hat{g}_{jn}(x)-\mathbf{E}\hat{g}%
_{jn}(x))w_{j}(x)dx-a_{jn}\right\} >z_{1-\alpha }\right\} +o(1).
\end{equation*}%
By Theorem 1, the leading probability converges to $\alpha $ as $%
n\rightarrow \infty $, delivering the desired result. $\blacksquare $\bigskip

\noindent \textsc{Proof of Theorem 3:} Fix $j$ such that $\Gamma
_{j}(g_{j})>0$. We focus on the case with $p>1$. The proof in the case with $%
p=1$ is simpler and hence omitted. Using the triangular inequality, we bound 
$\left\vert \Gamma _{j}(\hat{g}_{jn})-\Gamma _{j}(g_{j})\right\vert $ by%
\begin{eqnarray*}
&&\left\vert \int_{\mathcal{X}}\left\{ \Lambda _{p}(\hat{g}_{jn}(x))-\Lambda
_{p}(\mathbf{E}\hat{g}_{jn}(x))\right\} w_{j}(x)dx\right\vert \\
&&+\left\vert \int_{\mathcal{X}}\left\{ \Lambda _{p}(\mathbf{E}\hat{g}%
_{jn}(x))-\Lambda _{p}(g_{j}(x))\right\} w_{j}(x)dx\right\vert .
\end{eqnarray*}%
There exists $n_{0}$ such that for all $n\geq n_{0}$, $\sup_{x\in \mathcal{S}%
_{j}}|\mathbf{E}\hat{g}_{jn}(x)|<\infty $ by Lemma A4. Also, note that sup$%
_{x\in \mathcal{S}_{j}}|g_{j}(x)|<\infty $ by Assumption 1(i). Hence,
applying Lemma A3, from some large $n$ on, for some $C_{1},C_{2}>0$,%
\begin{eqnarray*}
\left\vert \Gamma _{j}(\hat{g}_{jn})-\Gamma _{j}(g_{j})\right\vert &\leq
&C_{1}\sum_{k=0}^{\lceil p-1\rceil }\int_{\mathcal{X}}|\hat{g}_{jn}(x)-%
\mathbf{E}\hat{g}_{jn}(x)|^{p-kz}w_{j}(x)dx \\
&&+C_{2}\sum_{k=0}^{\lceil p-1\rceil }\int_{\mathcal{X}}|\mathbf{E}\hat{g}%
_{jn}(x))-g_{j}(x)|^{p-kz}w_{j}(x)dx,
\end{eqnarray*}%
where $z=(p-1)/\lceil p-1\rceil $. Observe that $0\leq z\leq 1$.

As for the second integral, take $\varepsilon >0$ and a compact set $%
D\subset \mathbf{R}^{d}$ such that $\int_{\mathcal{X}\backslash
D}w_{j}(x)dx<\varepsilon $ and $g_{j}$ is continuous on $D$. Such a set $D$
exists by Lemma A1. Since $D$ is compact, $g_{j}$ is in fact uniformly
continuous on $D$. By change of variables,%
\begin{eqnarray*}
\mathbf{E}\hat{g}_{jn}(x)-g_{j}(x) &=&\int_{[-1/2,1/2]^{d}}\left\{
g_{j}(x-uh)K(u)-g_{j}(x)\right\} du \\
&=&\int_{[-1/2,1/2]^{d}}\left\{ g_{j}(x-uh)-g_{j}(x)\right\} K(u)du
\end{eqnarray*}%
and obtain that for $k=0,1,\cdot \cdot \cdot ,p-1,$%
\begin{eqnarray*}
&&\int_{\mathcal{X}}\left\vert \mathbf{E}\hat{g}_{jn}(x)-g_{j}(x)\right\vert
^{p-kz}w_{j}(x)dx \\
&=&\int_{D}\left\vert \mathbf{E}\hat{g}_{jn}(x)-g_{j}(x)\right\vert
^{p-kz}w_{j}(x)dx+\int_{\mathcal{X}\backslash D}\left\vert \mathbf{E}\hat{g}%
_{jn}(x)-g_{j}(x)\right\vert ^{p-kz}w_{j}(x)dx \\
&\leq &C_{3}\sup_{u\in \lbrack -1/2,1/2]^{d}}\sup_{x\in D\cap \mathcal{S}%
_{j}}\left\vert g_{j}(x-uh)-g_{j}(x)\right\vert ^{p-kz} \\
&&+C_{4}\int_{\mathcal{X}\backslash D}\int_{[-1/2,1/2]^{d}}\left\vert
g_{j}(x-uh)-g_{j}(x)\right\vert ^{p-kz}w_{j}(x)dudx,
\end{eqnarray*}%
for some positive constants $C_{3}$ and $C_{4}$. Note that the constant $%
C_{4}$ involves $||K||_{\infty }$. The first term is $o(1)$ as $h\rightarrow
0$, because $g_{j}$ is uniformly continuous on $D$. By Assumption 1(i), the
last term is bounded by%
\begin{equation*}
C_{5}\int_{\mathcal{X}\backslash D}w_{j}(x)dx<C_{6}\varepsilon ,\text{ for
some }C_{5},C_{6}>0\text{,}
\end{equation*}%
for some large $n$ on. Since the choice of $\varepsilon $ was arbitrary, we
conclude that as $n\rightarrow \infty ,$%
\begin{equation*}
\left\vert \Gamma _{j}(\hat{g}_{jn})-\Gamma _{j}(g_{j})\right\vert \leq
C_{1}\int_{\mathcal{X}}\left\vert \hat{g}_{jn}(x)-\mathbf{E}\hat{g}%
_{jn}(x)\right\vert ^{p-kz}w_{j}(x)dx+o(1).
\end{equation*}%
As for the leading integral, from the result of Theorem 1 (replacing $%
\Lambda _{p}(\cdot )$ there by $|\cdot |^{p-kz}$), we find that 
\begin{equation*}
\int_{\mathcal{X}}\left\vert \hat{g}_{jn}(x)-\mathbf{E}\hat{g}%
_{jn}(x)\right\vert
^{p-kz}w_{j}(x)dx=O_{P}(n^{-(p-kz)/2}h^{-(p-kz-1)d/2-d/2}).
\end{equation*}%
Since $n^{-1/2}h^{-d/2}\rightarrow 0$ by the condition of the theorem, we
conclude that $\Gamma _{j}(\hat{g}_{jn})\overset{p}{\rightarrow }\Gamma
_{j}(g_{j}).$ Using the similar argument, we can also show that%
\begin{equation*}
\hat{\sigma}_{n}^{2}\overset{p}{\rightarrow }\sigma ^{2}\text{ and }\hat{a}%
_{jn}=O_{P}(h^{-d/2})\text{ for all }j\in \mathcal{J},
\end{equation*}%
where $\sigma ^{2}=\mathbf{1}^{\prime }\Sigma \mathbf{1}>0.$ Hence%
\begin{equation*}
\hat{\sigma}_{n}^{-1}\{\Gamma _{j}(\hat{g}_{jn})-n^{-p/2}h^{-pd/2}h^{d/2}%
\hat{a}_{jn}\}\overset{p}{\rightarrow }\sigma ^{-1}\Gamma _{j}(g_{j})>0.
\end{equation*}%
Therefore,%
\begin{equation*}
P\{\hat{T}_{n}>z_{1-\alpha }\}\geq P\left\{ \sigma ^{-1}\Gamma
_{j}(g_{j})>0\right\} +o(1)\rightarrow 1,
\end{equation*}%
where the inequality holds by the fact that $n^{-1/2}h^{-d/2}\rightarrow 0$
and $\hat{a}_{jn}=O_{P}(h^{-d/2})$. $\blacksquare $\bigskip

\noindent \textsc{Lemma A11:} \textit{Suppose that Assumptions 1-3 hold, }$%
n^{-1/2}h^{-d}\rightarrow 0$, \textit{and that }$\sqrt{n}g_{j}(\cdot
)=\delta _{j}(\cdot ),$ $j\in \mathcal{J}$,\textit{\ for real bounded
functions} $\delta _{j},\ j\in \mathcal{J},$ \textit{for each} $n$. \textit{%
Then,}%
\begin{equation*}
\frac{1}{\sigma _{n}}\sum_{j=1}^{J}\left\{ n^{p/2}h^{(p-1)d/2}\Gamma
_{j,\delta }(\hat{g}_{jn})-\tilde{a}_{jn}\right\} \overset{d}{\rightarrow }%
N(0,1),
\end{equation*}%
\textit{where }$\tilde{a}_{jn}\equiv \int \mathbf{E}\Lambda
_{p}(h^{-d/(2p)}\rho _{jn}(x)\mathbb{Z}_{1}+h^{d(p-1)/(2p)}\delta
_{jn}(x))w_{j}(x)dx$ \textit{and} $\delta _{jn}(x)\equiv \int \delta
_{j}(x-uh)K(u)du$.\bigskip

\noindent \textsc{Proof:} By change of variables,%
\begin{equation*}
\sqrt{n}\mathbf{E}\hat{g}_{jn}(x)=\sqrt{n}\int g_{j}(x-uh)K(u)du=\int \delta
_{j}(x-uh)K(u)du.
\end{equation*}%
Since $\delta _{j}$ is bounded, $\sup_{x\in \mathcal{S}_{j}}\sqrt{n}%
\left\vert \mathbf{E}\hat{g}_{jn}(x)\right\vert =O(1)$. Hence%
\begin{equation}
\frac{\sqrt{nh^{d}}\hat{g}_{jN}(x)}{\rho _{jn}(x)}=\xi _{jn}(x)+\frac{\sqrt{%
nh^{d}}\mathbf{E}\hat{g}_{jn}(x)}{\rho _{jn}(x)}=\xi _{jn}(x)+O(h^{d/2}),
\label{conv}
\end{equation}%
under the local alternatives. Using this and following the proof of Lemma
A7, we find that under the local alternatives, $\sigma _{jk,n}\rightarrow
\sigma _{jk}.$ Also, as in the proof of Theorem 1, we use (\ref{conv}) and
deduce that%
\begin{equation}
\frac{1}{\sigma _{n}}\sum_{j=1}^{J}n^{p/2}h^{(p-1)d/2}\left\{ \Gamma _{j}(%
\hat{g}_{jn})-\mathbf{E}\Gamma _{j}(\hat{g}_{jn})\right\} \overset{d}{%
\rightarrow }N(0,1).  \label{conv6}
\end{equation}%
Now, as for $n^{p/2}h^{(p-1)d/2}\sigma _{n}^{-1}\mathbf{E}\Gamma _{j}(\hat{g}%
_{jn})$, We first note that%
\begin{eqnarray*}
n^{p/2}h^{(p-1)d/2}\Gamma _{j}(\hat{g}_{jn}) &=&h^{-d/2}\Gamma
_{j}(n^{1/2}h^{d/2}\{\hat{g}_{jn}-\mathbf{E}\hat{g}_{jn}\}+n^{1/2}h^{d/2}%
\mathbf{E}\hat{g}_{jn}) \\
&=&\Gamma _{j}(h^{-d/(2p)}\rho _{jn}(x)\xi _{jn}(x)+h^{(p-1)d/(2p)}\delta
_{jn}(x)).
\end{eqnarray*}%
We follow the proof of Lemma A4 and Lemma A6 (applying Lemma A2 with $%
\Lambda _{p}(v)$ in Lemma A6 replaced by $\Lambda
_{p}(v+h^{d(p-1)/(2p)}\delta _{jn}(x)/\rho _{jn}(x))$) to deduce that%
\begin{equation*}
\int \left\{ n^{p/2}h^{(p-1)d/2}\mathbf{E}\Lambda _{p}(\hat{g}_{jn}(x))-%
\mathbf{E}\Lambda _{p}(\bar{Z}_{jn}(x))\right\} w_{j}(x)dx\rightarrow 0,
\end{equation*}%
where $\bar{Z}_{jn}(x)\equiv h^{-d/(2p)}\rho _{jn}(x)\mathbb{Z}%
_{1}+h^{d(p-1)/(2p)}\delta _{jn}(x).$ $\blacksquare $\bigskip

\noindent \textsc{Proof of Theorem 4:} Under the local alternatives, by (\ref%
{conv5}) and (\ref{conv6}),%
\begin{eqnarray}
&&P\{\hat{T}_{n}>z_{1-\alpha }\}  \label{der} \\
&=&P\{\hat{\sigma}_{n}^{-1}\Sigma _{j=1}^{J}\{n^{p/2}h^{(p-1)d/2}\Gamma _{j}(%
\hat{g}_{jn})-\hat{a}_{jn}\}>z_{1-\alpha }\}  \notag \\
&=&P\{\sigma ^{-1}\Sigma _{j=1}^{J}\{n^{p/2}h^{(p-1)d/2}\Gamma _{j}(\hat{g}%
_{jn})-\tilde{a}_{jn}+\tilde{a}_{jn}-\hat{a}_{jn}\}>z_{1-\alpha }\}+o(1) 
\notag \\
&=&P\{\mathbb{Z}_{1}+\sigma ^{-1}\Sigma _{j=1}^{J}\{\tilde{a}%
_{jn}-a_{jn}\})>z_{1-\alpha }\}+o(1).  \notag
\end{eqnarray}%
Fix $\varepsilon >0$ and take a compact set $A_{\varepsilon }\subset 
\mathcal{S}_{j}$ such that $\int_{\mathcal{S}_{j}\backslash A_{\varepsilon
}}w_{j}(x)dx<\varepsilon $. Furthermore, without loss of generality, let $%
A_{\varepsilon }$ be a set on which $\delta _{j}(\cdot )$ and $\rho
_{j}(\cdot )$ are uniformly continuous. Then for any $\varepsilon _{1}>0$,
there exists $\lambda >0$ such that $\sup_{z\in \mathbf{R}%
^{d}:||x-z||<\lambda }|\delta _{j}(z)-\delta _{j}(x)|\leq \varepsilon _{1}$
uniformly over $x\in A_{\varepsilon }$. Hence from some large $n$ on,%
\begin{equation*}
\sup_{x\in A_{\varepsilon }}|\delta _{jn}(x)-\delta _{j}(x)|\leq
\int_{\lbrack -1/2,1/2]^{d}}\sup_{x\in A_{\varepsilon }}\left\vert \delta
_{j}(x-uh)-\delta _{j}(x)\right\vert K(u)du\leq \varepsilon _{1}.
\end{equation*}%
Since the choice of $\varepsilon _{1}$ was arbitrary, we conclude that $%
|\delta _{jn}(x)-\delta _{j}(x)|\rightarrow 0$ uniformly over $x\in
A_{\varepsilon }$. Similarly, we also conclude that $|\rho _{jn}(x)-\rho
_{j}(x)|\rightarrow 0$ uniformly over $x\in A_{\varepsilon }$. Using these
facts, we analyze $\sigma ^{-1}\Sigma _{j=1}^{J}\{\tilde{a}_{jn}-a_{jn}\}$
for each case of $p\in \{1,2\}$.\bigskip

\noindent (i) Suppose $p=1.$\ For $\gamma >0$ and $\mu \in \mathbf{R}$,%
\begin{eqnarray*}
\mathbf{E}\max \{\gamma \mathbb{Z}_{1}+\mu ,0\} &=&\mathbf{E}[\gamma \mathbb{%
Z}_{1}+\mu |\gamma \mathbb{Z}_{1}+\mu >0]P\left\{ \gamma \mathbb{Z}_{1}+\mu
>0\right\} \\
&=&\{\mu +\gamma \phi (-\mu /\gamma )/(1-\Phi (-\mu /\gamma ))\}\left(
1-\Phi (-\mu /\gamma )\right) \\
&=&\mu \left( 1-\Phi (-\mu /\gamma )\right) +\gamma \phi (-\mu /\gamma ) \\
&=&\mu \Phi (\mu /\gamma )+\gamma \phi (\mu /\gamma ).
\end{eqnarray*}

Taking $\gamma _{jn}\equiv h^{-d/2}\rho _{jn}(x)$, we have%
\begin{eqnarray*}
&&\mathbf{E}\max \{\gamma _{jn}\mathbb{Z}_{1}+\delta _{jn}(x),0\}-\mathbf{E}%
\max \{\gamma _{jn}\mathbb{Z}_{1},0\} \\
&=&\delta _{jn}(x)\Phi (\delta _{jn}(x)/\gamma _{jn})+\gamma _{jn}\phi
(\delta _{jn}(x)/\gamma _{jn})-\gamma _{jn}\phi (0) \\
&=&\delta _{jn}(x)\Phi (0)+O(h^{d/2}),
\end{eqnarray*}%
uniformly in $x\in \mathcal{S}_{j}$. Therefore, we can write $%
\lim_{n\rightarrow \infty }\{\tilde{a}_{jn}-a_{jn}\}$ as%
\begin{eqnarray*}
&&\lim_{n\rightarrow \infty }\int_{\mathcal{X}}\mathbf{E}[\Lambda
_{1}(h^{-d/2}\rho _{jn}(x)\mathbb{Z}_{1}+\delta _{jn}(x))-\Lambda
_{1}(h^{-d/2}\rho _{jn}(x)\mathbb{Z}_{1})]w_{j}(x)dx \\
&=&\frac{1}{2}\int_{A_{\varepsilon }}\delta _{j}(x)w_{j}(x)dx+\frac{1}{2}%
\lim_{n\rightarrow \infty }\int_{\mathcal{X}\backslash A_{\varepsilon
}}\delta _{jn}(x)w_{j}(x)dx.
\end{eqnarray*}%
Since $\delta _{jn}$ is uniformly bounded, there exists $C>0$ such that the
last integral is bounded by $C\varepsilon $. Since the choice of $%
\varepsilon >0$ was arbitrary, in view of (\ref{der}), this gives the
desired result.\bigskip

\noindent (ii) Suppose $p=2.$ For $\gamma >0$ and $\mu \in \mathbf{R}$,%
\begin{eqnarray*}
\mathbf{E}\max \{\gamma \mathbb{Z}_{1}+\mu ,0\}^{2} &=&\mathbf{E}[\left(
\gamma \mathbb{Z}_{1}+\mu \right) ^{2}|\gamma \mathbb{Z}_{1}+\mu >0]P\left\{
\gamma \mathbb{Z}_{1}+\mu >0\right\} \\
&=&(\mu ^{2}+\gamma ^{2})\Phi (\mu /\gamma )+\mu \gamma \phi (\mu /\gamma ).
\end{eqnarray*}%
Taking $\gamma _{jn}\equiv h^{-d/4}\rho _{jn}(x)$ and $\mu
_{jn}=h^{d/4}\delta _{jn}(x)$, we have%
\begin{eqnarray*}
&&\mathbf{E}\max \{\gamma _{jn}\mathbb{Z}_{1}+\mu _{jn},0\}^{2}-\mathbf{E}%
\max \{\gamma _{jn}\mathbb{Z}_{1},0\}^{2} \\
&=&\gamma _{jn}^{2}\{\Phi (\mu _{jn}/\gamma _{jn})-\Phi (0)\}+\mu
_{jn}^{2}\Phi (\mu _{jn}/\gamma _{jn})+\mu _{jn}\gamma _{jn}\phi (\mu
_{jn}/\gamma _{jn}) \\
&=&\{\mu _{jn}\gamma _{jn}\phi (0)+O(h^{d/2})\}+O(h^{d/2})+\{\mu _{jn}\gamma
_{jn}\phi (0)+O(h^{d})\} \\
&=&2\phi (0)\delta _{jn}(x)\rho _{jn}(x)+O(h^{d/2}),\text{ uniformly in }%
x\in \mathcal{S}_{j}.
\end{eqnarray*}%
Hence we write $\lim_{n\rightarrow \infty }\{\tilde{a}_{jn}-a_{jn}\}$ as%
\begin{eqnarray*}
&&\lim_{n\rightarrow \infty }\int_{\mathcal{S}_{j}}\mathbf{E}[\Lambda
_{2}(h^{-d/4}\rho _{jn}(x)\mathbb{Z}_{1}+h^{d/4}\delta _{jn}(x))-\Lambda
_{2}(h^{-d/4}\rho _{jn}(x)\mathbb{Z}_{1})]w_{j}(x)dx \\
&=&2\phi (0)\lim_{n\rightarrow \infty }\int_{\mathcal{S}_{j}}\delta
_{jn}(x)\rho _{jn}(x)w_{j}(x)dx+O(h^{d/2}) \\
&=&\sqrt{\frac{2}{\pi }}\lim_{n\rightarrow \infty }\int_{A_{\varepsilon
}}\delta _{jn}(x)\rho _{jn}(x)w_{j}(x)dx+\sqrt{\frac{2}{\pi }}%
\lim_{n\rightarrow \infty }\int_{\mathcal{S}_{j}\backslash A_{\varepsilon
}}\delta _{jn}(x)\rho _{jn}(x)w_{j}(x)dx+O(h^{d/2}).
\end{eqnarray*}%
The second term is bounded by $C\varepsilon $ for some $C>0$, because $%
\delta _{jn}\rho _{jn}$ is bounded. Since the choice of $\varepsilon >0$ was
arbitrary and%
\begin{equation*}
\int_{A_{\varepsilon }}\delta _{jn}(x)\rho _{jn}(x)w_{j}(x)dx\rightarrow
\int_{A_{\varepsilon }}\delta _{j}(x)\rho _{j}(x)w_{j}(x)dx,\text{ as }%
n\rightarrow \infty ,
\end{equation*}%
in view of (\ref{der}), this gives the desired result. $\blacksquare $%
\bigskip

\textsc{Proof of Theorem 4}$^{\ast }$\textsc{:} Let $A_{\varepsilon }\subset 
\mathcal{S}_{j}$ be defined as in the proof of Theorem 4.

\noindent (i) Suppose $p=1$. Under $H_{\delta }^{\ast }$, take $\gamma
\equiv h^{-d/2}\rho _{jn}(x)$ and $\mu =h^{-d/4}\delta _{jn}(x)$ to get%
\begin{eqnarray*}
&&\mathbf{E}\max \{\gamma \mathbb{Z}_{1}+\mu ,0\}-\mathbf{E}\max \{\gamma 
\mathbb{Z}_{1},0\} \\
&=&h^{-d/4}\delta _{jn}(x)\Phi (h^{d/4}\delta _{jn}(x)/\rho
_{jn}(x))+h^{-d/2}\rho _{jn}(x)\left[ \phi (h^{d/4}\delta _{jn}(x)/\rho
_{jn}(x))-\phi (0)\right] \\
&=&h^{-d/4}\delta _{jn}(x)\Phi (0)+\frac{1}{2}\phi (0)\left[ \delta
_{jn}^{2}(x)/\rho _{jn}(x)\right] +O(h^{d/4}),
\end{eqnarray*}%
uniformly in $x\in \mathcal{S}_{j}$. Therefore, if $\eta _{1,0}(w,\delta )=0$
under $H_{\delta }^{\ast },$ we can write $\lim_{n\rightarrow \infty }\{%
\tilde{a}_{jn}-a_{jn}\}$ as%
\begin{eqnarray*}
&&\lim_{n\rightarrow \infty }\int_{\mathcal{X}}\mathbf{E}[\Lambda
_{1}(h^{-d/2}\rho _{jn}(x)\mathbb{Z}_{1}+h^{-d/4}\delta _{jn}(x))-\Lambda
_{1}(h^{-d/2}\rho _{jn}(x)\mathbb{Z}_{1})]w_{j}(x)dx \\
&=&\frac{1}{2}\phi (0)\lim_{n\rightarrow \infty }\int_{\mathcal{S}_{j}}\left[
\delta _{jn}^{2}(x)/\rho _{jn}(x)\right] w_{j}(x)dx \\
&=&\frac{1}{2}\phi (0)\int_{A_{\varepsilon }}\left[ \delta _{j}^{2}(x)/\rho
_{j}(x)\right] w_{j}(x)dx+\frac{1}{2}\phi (0)\lim_{n\rightarrow \infty
}\int_{\mathcal{S}_{j}\backslash A_{\varepsilon }}\left[ \delta
_{jn}^{2}(x)/\rho _{jn}(x)\right] w_{j}(x)dx.
\end{eqnarray*}%
Since $\delta _{jn}^{2}/\rho _{jn}$ is uniformly bounded and the choice of $%
\varepsilon >0$ is arbitrary, we get the desired result.\bigskip

\noindent (ii) Suppose $p=2$. Under $H_{\delta }^{\ast }$, we take $\gamma
\equiv h^{-d/4}\rho _{jn}(x)$ and $\mu =\delta _{jn}(x)$, so that, by a
Taylor expansion,%
\begin{eqnarray*}
&&\mathbf{E}\max \{\gamma \mathbb{Z}_{1}+\mu ,0\}^{2}-\mathbf{E}\max
\{\gamma \mathbb{Z}_{1},0\}^{2} \\
&=&\gamma ^{2}\{\Phi (\mu /\gamma )-\Phi (0)\}+\mu ^{2}\Phi (\mu /\gamma
)+\mu \gamma \phi (\mu /\gamma ) \\
&=&\left\{ \phi (0)\mu \gamma +\frac{1}{6}\phi ^{\prime \prime }(a^{\ast })%
\frac{\mu ^{3}}{\gamma }\right\} +\left\{ \Phi (0)\mu ^{2}+\phi (a^{\ast })%
\frac{\mu ^{3}}{\gamma }\right\} +\mu \gamma \left\{ \phi (0)+\frac{1}{2}%
\phi ^{\prime \prime }(a^{\ast })\frac{\mu ^{3}}{\gamma }\right\} \\
&=&h^{-d/4}\cdot 2\phi (0)\delta _{jn}(x)\rho _{jn}(x)+\frac{1}{2}\delta
_{jn}^{2}(x)+O(h^{d/4}),\text{ }
\end{eqnarray*}%
uniformly in $x\in \mathcal{S}_{j}$, where $a^{\ast }$ denotes a term that
lies between $0$ and $\mu /\gamma .$ Therefore, if $\eta _{1,1}(w,\delta )=0$
under $H_{2\delta },$ then we can write $\lim_{n\rightarrow \infty }\{\tilde{%
a}_{jn}-a_{jn}\}$ as%
\begin{eqnarray*}
&&\lim_{n\rightarrow \infty }\int_{\mathcal{X}}\mathbf{E}[\Lambda
_{2}(h^{-d/2}\rho _{jn}(x)\mathbb{Z}_{1}+h^{-d/4}\delta _{jn}(x))-\Lambda
_{2}(h^{-d/2}\rho _{jn}(x)\mathbb{Z}_{1})]w_{j}(x)dx \\
&=&\frac{1}{2}\lim_{n\rightarrow \infty }\int_{\mathcal{S}_{j}}\delta
_{jn}^{2}(x)w_{j}(x)dx \\
&=&\frac{1}{2}\int_{A_{\varepsilon }}\delta _{j}^{2}(x)w_{j}(x)dx+\frac{1}{2}%
\lim_{n\rightarrow \infty }\int_{\mathcal{S}_{j}\backslash A_{\varepsilon
}}\delta _{jn}^{2}(x)w_{j}(x)dx.
\end{eqnarray*}%
Since $\delta _{jn}^{2}$ is uniformly bounded and the choice of $\varepsilon
>0$ is arbitrary, we get the desired result.\bigskip\ $\blacksquare $\bigskip

\noindent \textsc{Proof of Theorem 5:} Similarly as before, we fix $%
\varepsilon >0$ and take a compact set $A_{\varepsilon }\subset \mathcal{S}%
_{j}$ such that $\int_{\mathcal{S}_{j}\backslash A_{\varepsilon
}}w_{j}(x)dx<\varepsilon $ and $\delta _{j}(\cdot )$ and $\delta _{j}(\cdot
)\rho _{j}^{-1}(\cdot )$ are uniformly continuous on $A_{\varepsilon }$. By
change of variables and uniform continuity,%
\begin{eqnarray*}
\sup_{x\in A_{\varepsilon }}|\delta _{jn}(x)\rho _{jn}^{-1}(x)-\delta
_{j}(x)\rho _{j}^{-1}(x)| &\rightarrow &0\text{ and} \\
\sup_{x\in A_{\varepsilon }}|\delta _{jn}(x)-\delta _{j}(x)| &\rightarrow &0.
\end{eqnarray*}

\noindent (i) Suppose $p=1.$\ For $\gamma >0$ and $\mu \in \mathbf{R}$, 
\begin{equation*}
\mathbf{E}\left\vert \gamma \mathbb{Z}_{1}+\mu \right\vert =2\gamma \phi
(\mu /\gamma )+2\mu \left[ \Phi (\mu /\gamma )-1/2\right] .
\end{equation*}%
With $\gamma _{jn}\equiv h^{-d/2}\rho _{jn}(x)\ $and $\mu
_{jn}=h^{-d/4}\delta _{jn}(x)$, we find that uniformly over $x\in \mathcal{S}%
_{j},$%
\begin{eqnarray*}
&&\mathbf{E}|\gamma _{jn}\mathbb{Z}_{1}+\mu _{jn}|-\mathbf{E}|\gamma _{jn}%
\mathbb{Z}_{1}| \\
&=&2\gamma _{jn}[\phi (\mu _{jn}/\gamma _{jn})-\phi (0)]+2\mu _{jn}[\Phi
(\mu _{jn}/\gamma _{jn})-1/2] \\
&=&\left[ \phi ^{\prime \prime }(0)+2\phi (0)\right] \delta _{jn}^{2}(x)\rho
_{jn}^{-1}(x)+O(h^{d/4}).
\end{eqnarray*}%
Therefore, we write $\lim_{n\rightarrow \infty }\{\tilde{a}_{jn}-a_{jn}\}$ as%
\begin{eqnarray*}
&&\lim_{n\rightarrow \infty }\int_{\mathcal{S}_{j}}\mathbf{E}[\Lambda
_{1}(h^{-d/2}\rho _{jn}(x)\mathbb{Z}_{1}+n^{-d/4}\delta _{jn}(x))-\Lambda
_{1}(h^{-d/2}\rho _{jn}(x)\mathbb{Z}_{1})]w_{j}(x)dx \\
&=&\frac{1}{\sqrt{2\pi }}\lim_{n\rightarrow \infty }\int_{\mathcal{S}%
_{j}}\delta _{jn}^{2}(x)\rho _{jn}^{-1}(x)w_{j}(x)dx+O(h^{d/4}) \\
&=&\frac{1}{\sqrt{2\pi }}\int_{A_{\varepsilon }}\delta _{j}^{2}(x)\rho
_{j}^{-1}(x)w_{j}(x)dx+\frac{1}{\sqrt{2\pi }}\lim_{n\rightarrow \infty
}\int_{\mathcal{S}_{j}\backslash A_{\varepsilon }}\delta _{jn}^{2}(x)\rho
_{jn}^{-1}(x)w_{j}(x)dx+o(1).
\end{eqnarray*}%
By Assumption 4 and Lemma A4, $\delta _{jn}^{2}(x)\rho _{jn}^{-1}(x)$ is
bounded uniformly over $x\in \mathcal{S}_{j}$, enabling us to bound the
second integral by $C\varepsilon $ for some $C>0.$ Since $\varepsilon $ is
arbitrarily chosen, in view of (\ref{der}), this gives the desired
result.\bigskip

\noindent (ii) Suppose $p=2.\ $We have, for each $x\in \mathcal{S}_{j},$%
\begin{equation*}
\mathbf{E}\{h^{-d/4}\rho _{jn}(x)\mathbb{Z}_{1}+\delta _{jn}(x)\}^{2}-%
\mathbf{E}\{h^{-d/4}\rho _{jn}(x)\mathbb{Z}_{1}\}^{2}=\delta _{jn}^{2}(x).
\end{equation*}%
Therefore, we write $\lim_{n\rightarrow \infty }\{\tilde{a}_{jn}-a_{jn}\}$ as%
\begin{eqnarray*}
&&\lim_{n\rightarrow \infty }\int_{\mathcal{X}}\mathbf{E}[\Lambda
_{2}(h^{-d/4}\rho _{jn}(x)\mathbb{Z}_{1}+\delta _{jn}(x))-\Lambda
_{2}(h^{-d/4}\rho _{jn}(x)\mathbb{Z}_{1})]w_{j}(x)dx \\
&=&\int_{A_{\varepsilon }}\delta _{j}^{2}(x)w_{j}(x)dx+\lim_{n\rightarrow
\infty }\int_{\mathcal{S}_{j}\backslash A_{\varepsilon }}\delta
_{jn}(x)w_{j}(x)dx+o(1)
\end{eqnarray*}%
The second integral is bounded by $C\varepsilon $ for some $C>0$, and in
view of (\ref{der}), this gives the desired result. $\blacksquare $\bigskip

\clearpage \newpage 

\clearpage \newpage 
\begin{figure}[htbp]
\caption{Results of Monte Carlo Experiments: $L_1$ test and $\protect\sigma%
(x) = x$}
\label{figure2}
\begin{center}
\makebox{
\includegraphics[origin=bl,scale=.8,angle=180]{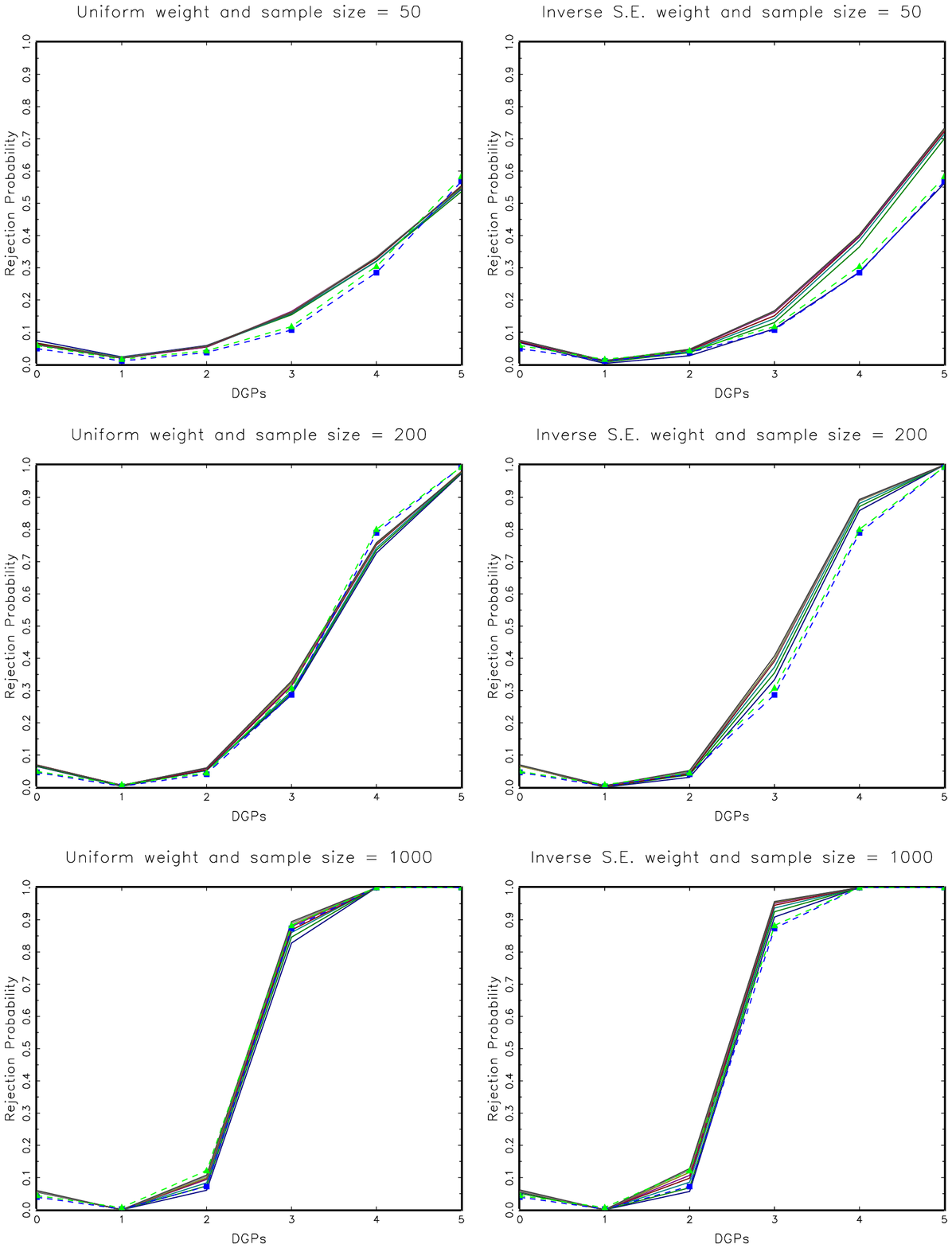}
}
\end{center}
\par
\parbox{5in}{Notes: See notes in Figure \ref{figure1}.}
\end{figure}

\clearpage \newpage 
\begin{figure}[htbp]
\caption{Results of Monte Carlo Experiments: $L_2$ test and $\protect\sigma%
(x) \equiv 1$}
\label{figure3}
\begin{center}
\makebox{
\includegraphics[origin=bl,scale=.8,angle=180]{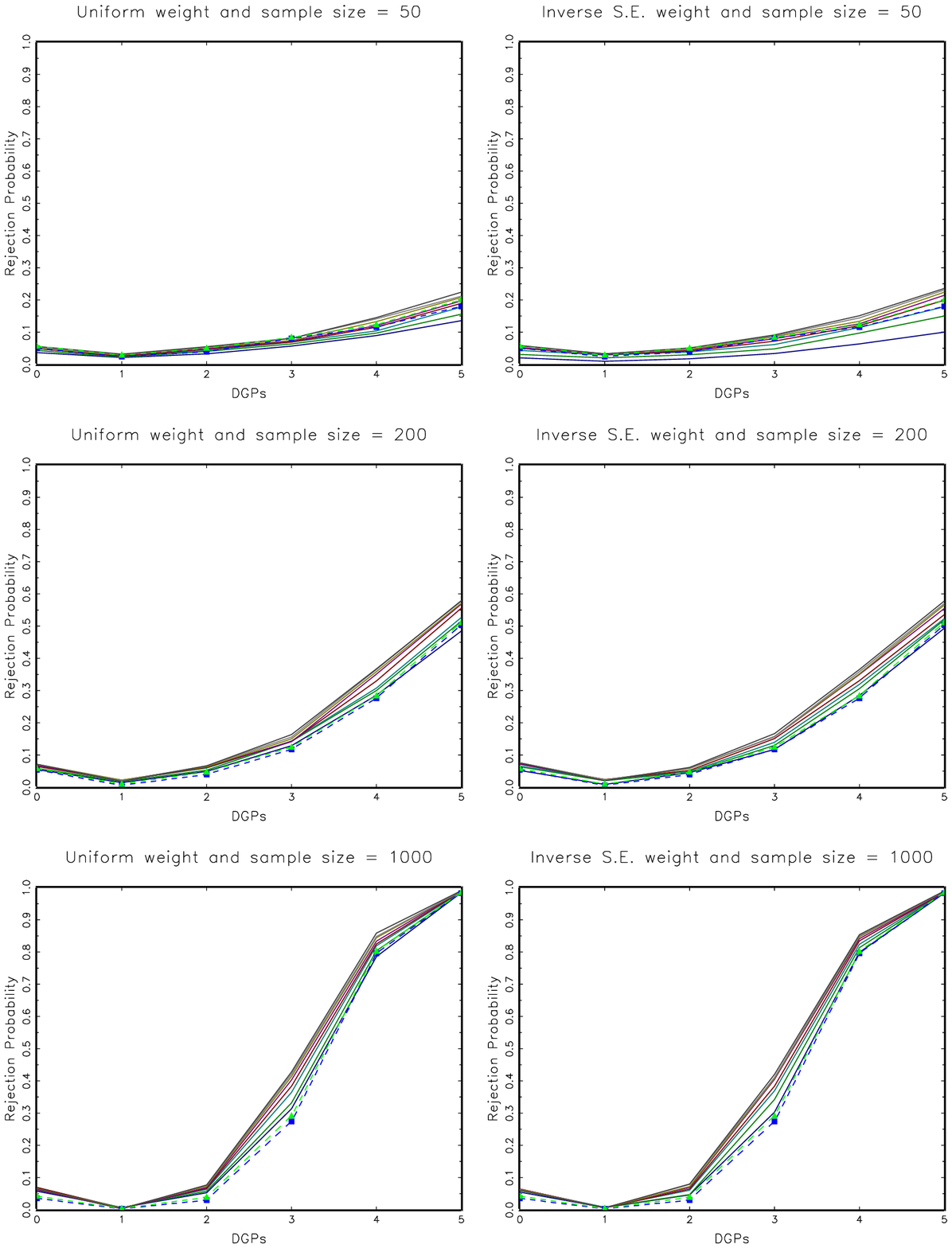}
}
\end{center}
\par
\parbox{5in}{Notes: See notes in Figure \ref{figure1}.}
\end{figure}

\clearpage \newpage 
\begin{figure}[htbp]
\caption{Results of Monte Carlo Experiments: $L_2$ test and $\protect\sigma%
(x) = x$}
\label{figure4}
\begin{center}
\makebox{
\includegraphics[origin=bl,scale=.8,angle=180]{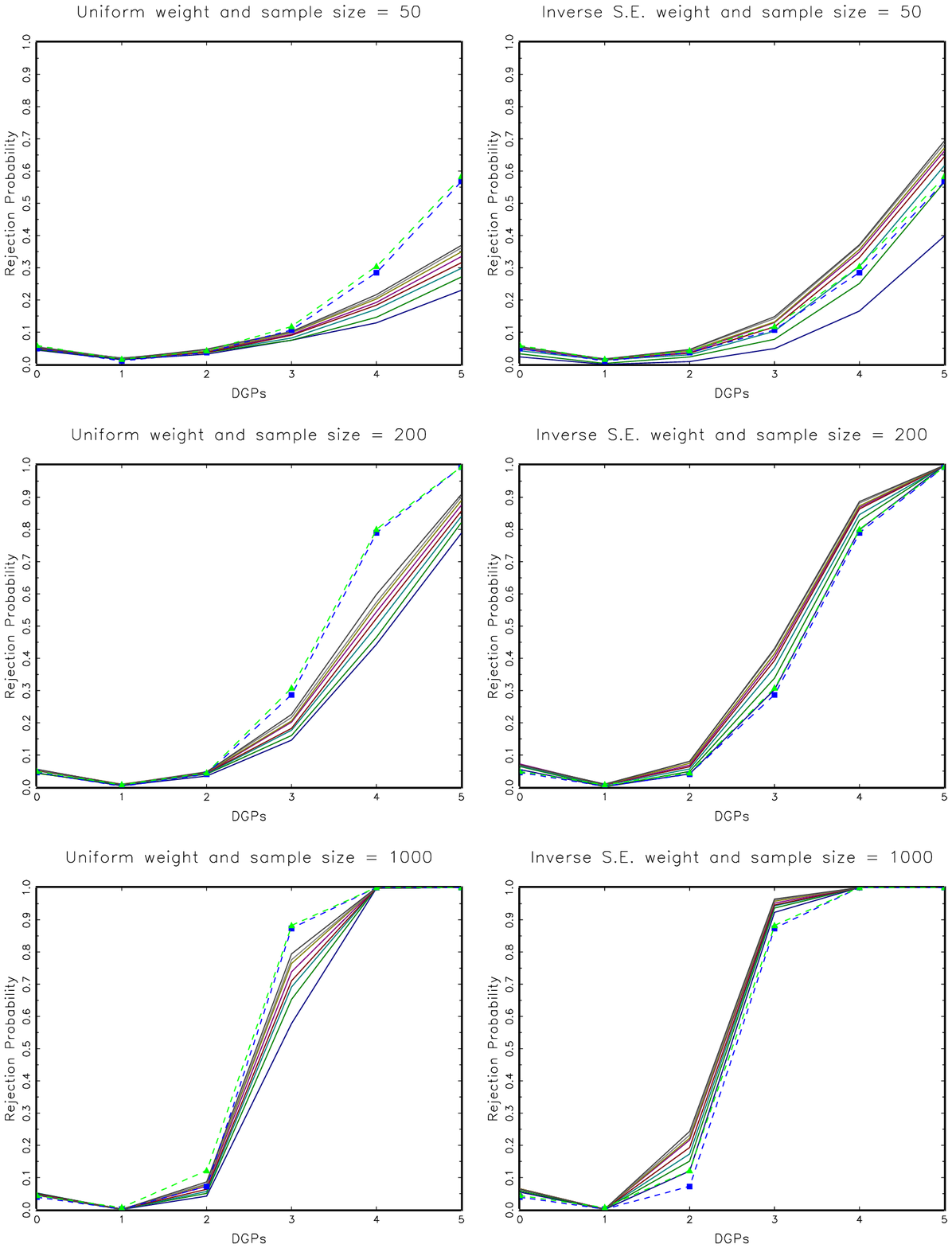}
}
\end{center}
\par
\parbox{5in}{Notes: See notes in Figure \ref{figure1}.}
\end{figure}

\end{document}